
\input amstex
\documentstyle{amsppt}
\magnification=1200
\hsize=16.5truecm
\vsize=23.2truecm
\voffset=0.5truecm
\catcode`\@=11
\redefine\logo@{}
\catcode`\@=13

\define\ms{\operatorname{Ms}}

\define\ml{\operatorname{Exp-Lift}}
\define\dv{\operatorname{Div}}
\define \fm{\frak M}
\define \fn{\frak N}
\define \gq{\Gamma_1(Q)}
\define \bn{\Bbb N}
\define \bz{\Bbb Z}
\define \bq{\Bbb Q}
\define \br{\Bbb R}
\define \bc{\Bbb C}
\define \bh{\Bbb H}
\define \bp{\Bbb P}

\define \M{{\Cal M}}
\define\Ha{{\Cal H}}
\define\La{{\Cal L}}
\define\geg{{\goth g}}
\define\0o{{\overline 0}}
\define\1o{{\overline 1}}
\define\rk{\text{rk}~}
\define\gi{\Gamma_{\infty}}
\define\gm{\Gamma}

\define\mult{\text{mult}}

\define\pd#1#2{\dfrac{\partial#1}{\partial#2}}
\define\Sym{{\text{Sym}}}

\TagsOnRight
\NoBlackBoxes

\document
\document

\topmatter
\title
Automorphic Forms and Lorentzian Kac--Moody algebras. Part II
\endtitle

\author
Valeri A. Gritsenko \footnote{Supported by
RIMS of Kyoto University
\hfill\hfill}
and
Viacheslav V. Nikulin \footnote{Supported by
Grant of Russian Fund of Fundamental Research and RIMS of Kyoto University
\hfill\hfill}
\endauthor

\date{November 25, 1996, ${}\quad{}$ alg-geom/9611028}\enddate

\address
St. Petersburg Department of Steklov Mathematical Institute,
\newline
${}\hskip 8pt $
Fontanka 27, 191011 St. Petersburg,  Russia
\endaddress
\email
gritsenk\@pdmi.ras.ru;
gritsenk\@mpim-bonn.mpg.de
\endemail

\address
Steklov Mathematical Institute,
ul. Vavilova 42, Moscow 117966, GSP-1, Russia
\endaddress
\email
slava\@nikulin.mian.su
\endemail

\abstract
We give variants of lifting  construction, which define new classes
of modular forms on the Siegel upper half-space of complex dimension
$3$ with respect to the full  paramodular  groups  (defining moduli
of Abelian surfaces with arbitrary polarization). The data for these
liftings are Jacobi forms of integral and half-integral indices.
In particular, we get modular forms which are generalizations
of the Dedekind eta-function. Some of these forms define automorphic
corrections of Lorentzian Kac--Moody algebras with hyperbolic
generalized Cartan matrices of rank three classified in Part I of this
paper. We also construct many automorphic forms which give discriminants
of moduli of K3 surfaces with conditions on Picard lattice.
These results  are important for Mirror Symmetry and theory of
Lorentzian Kac--Moody algebras.
\endabstract

\rightheadtext
{Lorentzian Kac--Moody algebras II}
\leftheadtext{V. Gritsenko and  V. Nikulin}
\endtopmatter

\document
\head
\S 0. Introduction
\endhead

In Part I we developed the general theory of reflective
automorphic forms and their particular case of Lie reflective
automorphic forms on Hermitian symmetric domains of type IV.
These automorphic forms are very important
in Mirror Symmetry (for K3's and Calabi--Yau's) and for
Lorentzian Kac--Moody algebras. In Part I, in particular, we
showed that this theory is similar to the theory of hyperbolic
root systems (it is its mirror symmetric variant).
We explained that reflective automorphic forms are very exceptional.
Conjecturally their number is finite similarly to finiteness results
for corresponding hyperbolic root systems with some condition of finiteness
of volume for fundamental polyhedron (i.e., of elliptic or parabolic type).
We believe and hope to show
in further publications that classification of reflective automorphic forms
and corresponding hyperbolic root systems is the key step
in classification of some important class of Calabi--Yau's
(see \cite{GN6}). 
For example, finiteness results for hyperbolic root systems of elliptic 
and parabolic type and for reflective automorphic
forms are related with finiteness of families of these Calabi--Yau's.

We demonstrated  in Part I a general
method of classification of the hyperbolic root systems
on the example of symmetric (and twisted to symmetric)
hyperbolic generalized Cartan matrices of elliptic type of rank 3
and with a lattice Weyl vector. Let us denote
the set of these matrices by $\goth A$. It contains 60 matrices.
In this Part II we consider
methods of construction of reflective automorphic forms. In
particular, for many generalized Cartan matrices
$A\in \goth A$ we find their mirror symmetric objects --- automorphic
forms $F$ defining automorphic corrections of Kac--Moody algebras
$\geg (A)$ corresponding to $A$. These forms $F$
define so-called automorphic Lorentzian Kac--Moody algebras $\geg_F$
containing $\geg (A)$ and having good automorphic properties.
In this paper we consider mainly  $3$-dimensional automorphic forms $F$
with respect to the full paramodular groups,
i.e. automorphic  forms on the Siegel
upper half-space $\bh_2$ of complex dimension $3$
with respect to the paramodular groups
$\Gamma_t\subset Sp_4(\bq)$
(the threefold $\Cal A_t=\Gamma_t\setminus \bh_2$
is the moduli space of Abelian surfaces with polarization
of type $(1,t)$).
We remark that the results  and methods of this paper
can be generalized  to automorphic forms
with  respect to $O(n,2)$ (see \cite{G4}--\cite{G5}).
We do not consider in this paper the generalized Cartan matrices
$A\in {\goth A}$  related with automorphic forms $F$ with respect to some
congruence subgroups of the paramodular groups.
We will treat them later.

\smallpagebreak

First we define the  advance version
of the Maass lifting (see \cite{M1}--\cite{M2})
in the form proposed in \cite{G1}--\cite{G5}.
This new variant (we call it {\it arithmetic lifting}) provides modular
forms with respect to the full paramodular group $\Gamma_t$
with a character $\chi:\Gamma_t\to \Bbb C^*$ of order $Q$
where
$Q$ is necessarily a divisor of $12$, because
$\chi$  is trivial on the commutator subgroup
of $\Gamma_t$.
The datum for  the lifting is a holomorphic Jacobi form
$\phi_{k,t}(\tau, z)$ of weight $k$ of  integral
or {\it half-integral} index $t$ with a character
$v_{SL}\times v_H$ of the  Jacobi group
$\Gamma^J=SL_2(\bz)\ltimes H(\bz)$
induced by a character
$v_{SL}:SL_2(\bz)\to \{\root {12}\of {1}\}$
of order $Q$
and by a binary character
$v_H: H(\bz)\to \{\pm 1\}$ of the  Heisenberg group.
In Theorem 1.12 we define several variants of the arithmetic lifting
numerated by $\mu\in (\bz/Q\bz)^*$:
$$
\hbox{Lift}_\mu(\phi)(Z)
=\sum\Sb m\equiv \mu \, mod \,Q\\
\vspace{0.5\jot} m>0\endSb
m^{2-k}\bigl(\widetilde{\phi}|_k \,T_-^{(Q)}(m)\bigr)(Z)
\tag{0.1}
$$
where
$\widetilde{\phi}_{k,t}(Z)=\phi_{k,t}(\tau,z)
\exp{(2\pi i\,t\omega)}$ is the  modular form with respect to the maximal
parabolic subgroup $\Gamma_{\infty,t}$
associated with $\phi_{k,t}(\tau, z)$
(this parabolic subgroup defines a $1$-dimensional
boundary component of the threefold $\Cal A_t$)
and where the operators $T_-^{(Q)}(m)$ are the ``minus''-embedding
of the standard  elements $T^{(Q)}(m)$ of the Hecke ring
of the principal congruence
subgroup $\Gamma_1(Q)\subset SL_2(\bz)$
in the Hecke ring of the parabolic
subgroup $\Gamma_{\infty,t}$.
The Fourier expansion of the lifting \thetag{0.1}
can be  written in terms of Fourier coefficients
of the Jacobi form used for the lifting.

In \S 1 we construct several series of Jacobi forms which one can use
as data for the arithmetic lifting \thetag{0.1}.
For example, the Jacobi theta-series
$$
\vartheta(\tau,z)=\vartheta_{11}(\tau,z)=
\sum_{m\in \bz}\,\biggl(\frac{-4}{m}\biggr)\, q^{{m^2}/8}\,r^{{m}/2}
$$
is a Jacobi form of weight $\frac{1}2$
and index $\frac 1{2}$.
(For matrix
$\bigl(\smallmatrix \tau&z\\ z&\omega\endsmallmatrix\bigr)
\in \bh_2$ we use  three formal variables
$q=\exp{(2\pi i\, \tau)}$,
$r=\exp{(2\pi i\, z)}$ and
$s=\exp{(2\pi i\,\omega)}$.)
A typical modular form which we construct
is the cusp form
$$
\Delta_1(Z)=\hbox{Lift}_1(\eta(\tau)\vartheta(\tau,z)))
=\sum_{M\ge 1}
\sum\Sb n,\,m >0,\, l\in \bz\\
\vspace{0.5\jot} n,\,m\equiv 1\,mod\,6\\
\vspace{0.5\jot} 4nm-3l^2=M^2\endSb
\biggl(\dsize\frac{-4}{l}\biggr)
\biggl(\dsize\frac{12}{M}\biggr)
\sum\Sb a|(n,l,m)\endSb \biggl(\dsize\frac{6}{a}\biggr)\,
q^{n/6}r^{l/2}s^{m/2}
\tag{0.2}
$$
of weight one
with respect to the paramodular group $\Gamma_3$
with a character of order $6$. This is the lifting
of $\eta(\tau)\vartheta(\tau,z)$ where
$\eta(\tau)$ is the Dedekind $\eta$-function.
For a primitive matrix 
$N=\left(\smallmatrix n&l/2\\ l/2&m\endsmallmatrix\right)$ 
(i.e. $(n,l,m)=1$) the Fourier
coefficient $f(N)$ of $\Delta_1(Z)$ is equal to $\pm 1$ or $0$.  One
can consider $\Delta_1(Z)$ as a variant of {\it Dedekind
$\eta$-function in three variables}.  Below we give some evidence for
this point of view.

The Dedekind $\eta$-function is a character for an affine Kac--Moody
algebra. The cusp form $\Delta_1(Z)$ defines the automorphic
correction $\frak g_{\Delta_1}$ for the symmetric generalized Cartan
matrix $A_{3,II}\in {\goth A}$ (see Theorem 1.3.1 in Part I).
It follows that $\Delta_1(Z)$ is the
denominator function of the generalized Kac--Moody superalgebra
$\frak g_{\Delta_1}$ and it  has an infinite product expansion due to the
Weyl--Kac--Borcherds denominator identity (see \cite{GN1, \S 6} and
\cite{R}). We find this expansion in \S 2.  It gives multiplicities of
roots of $\frak g_{\Delta_1}$.  To describe this infinite product
expansion, let us define a weak Jacobi form
$\phi_{0,3}(\tau,z)=\dsize\sum_{n\ge 0,\,l\in \bz}f_3(n,l)q^nr^l$ of
weight $0$ and index $3$ with integral Fourier coefficients: $$
\phi_{0,3}(\tau ,z)=
\biggl(\frac {\vartheta(\tau ,2z)}{\vartheta(\tau ,z)}\biggr)^2
=r^{-1}
\biggl(\prod_{n\ge 1}(1+q^{n-1}r)(1+q^{n}r^{-1})(1-q^{2n-1}r^2)
(1-q^{2n-1}r^{-2})\biggr)^2.  $$ We have the identity $$
\Delta_1(Z)=
q^{\frac {1}6}r^{\frac {1}2}s^{\frac {1}2}
\prod
\Sb n,\,l,\,m\in \Bbb Z\\
\vspace{0.5\jot}
(n,l,m)>0
\endSb
\bigl(1-q^n r^l s^{3m}\bigr)^{f_{3}(nm,l)}.
\tag{0.3}
$$ This formula is an example of the {\it exponential or Borcherds
lifting}. In general, one can define the exponential lifting of a
Jacobi form $\phi_{0,t}(\tau,z)=\sum_{n,l} f(n,l)q^nr^l$ of weight
zero as follows (compare with \thetag{0.1}): 
$$
\ml(\phi_{0,t})(Z)=\psi(Z)
\exp{\bigl( -\sum_{m\ge 1}
m^{-1}\widetilde{\phi}_{0,t}|_0\, T_-(m)(Z)\bigr)}
\tag{0.4}
$$ 
where $$
\psi(Z)=\eta(\tau)^{f(0,0)}
\prod_{l>0}
\biggl(\frac{\vartheta(\tau,lz)}
{\eta(\tau)}\exp{(\pi i\, l^2\omega)}\biggr)^{f(0,l)}.  $$ The
$T_-(m)$ are  the operators from \thetag{0.1} if we put there $Q=1$.
In Theorem 2.1 we describe properties of the lifting
\thetag{0.4}.
This theorem is a modified version of the Borcherds construction
given in \cite{Bo6, Theorem 10.1}.

The arithmetic and exponential liftings have similar behavior with
respect to the action of Hecke operators, if we modify, in appropriate
way, the action of Hecke operators on \thetag{0.4}.  Let
$X=\sum_{i}\Gamma_t g_i\in H(\Gamma_t)$ be an element of the Hecke
ring of the group $\Gamma_t$. For a modular form $F(Z)$ of weight $k$
with respect to $\Gamma_t$ we define the Hecke product $$
[F]_X(Z)=\prod_{i} (F|_k g_i) (Z).
\tag{0.5}
$$ If $X$ has a good reduction modulo $t$, then we proved in
\cite{GN4, Theorem A.7} that $$ [\ml(\phi_{0,t})]_X=
c\cdot\ml\bigl(\phi_{0,t}|_0\, \Cal J_0^{(t)}(X)\bigr) $$ where ${\Cal
J}^{(t)}_0(X)$ is a natural projection of the Hecke ring $\Cal
H(\Gamma_t)$ on the Hecke-Jacobi ring of the parabolic subgroup
$\Gamma_{\infty}$ and $c$ is a constant.

Another type of Hecke product is the multiplicative symmetrisation
$$
\ms_p:
F(Z)\mapsto \prod_{M_i\in \Gamma_t\cap \Gamma_{tp}
\setminus \Gamma_{tp}}
(F|_kM_i) (Z)
\tag{0.6}
$$ which maps modular forms with respect to $\Gamma_t$ to
modular forms with respect to $\Gamma_{pt}$ for a prime $p$.
Theorem 3.3 says that $$
\ms_p(\hbox{Exp-Lift}(\phi_{0,t}))=
c\cdot\hbox{Exp-Lift}(\phi_{0,t}| T_-(p)) $$ where $T_-(p)$ is the
Hecke operator used in \thetag{0.1} and \thetag{0.4}.  To prove the
relations of the exponential lifting with the defined Hecke
correspondences we use the methods developed in \cite{G2}, \cite{G6}--
\cite{G8} (see also \cite{G11}) 
in the the case of the lifting of type \thetag{0.1}.

In \S 2--\S 4 we find many cases when the arithmetic lifting (or a
finite Hecke product of type \thetag{0.5} or
\thetag{0.6} of lifted modular  forms)
is equal to the exponential lifting.  This is true for the most
fundamental Siegel modular forms $\Delta_5(Z)$ (the product of even
theta-constants), $\Delta_2(Z)$, $\Delta_1(Z)$ and $\Delta_{1/2}(Z)$
(the index denotes the weight) which are the modular forms with
respect to the groups $\Gamma_1$, $\Gamma_2$, $\Gamma_3$ and
$\Gamma_4$ respectively.  The first three functions are cusp forms.
The last one is a Siegel theta-constant (the ``most odd" even
theta-constant) which we consider as a ``trivial" lifting of the
Jacobi theta-series.  This interpretation leads to the fact that this
theta-constant is a modular form with respect to the full paramodular
group $\Gamma_4$ and has infinite product expansion of type
\thetag{0.3} (see Theorem 1.11 where also another singular modular
form with respect to $\Gamma_{36}$ is constructed).

The modular forms $\Delta_k(Z)$ ($k=1/2$, $1$, $2$, $5$) have very
simple divisors.  They
coincide with the irreducible Humbert surfaces $H_1$ of discriminant
$1$ in the corresponding modular threefolds.
The $\Delta_5(Z)$ is associated with the hyperbolic root system
defined by a triangle with vertices at infinity on the hyperbolic
plane (see Fig.1 in \cite{GN5}), the $\Delta_2(Z)$ and $\Delta_1(Z)$
are associated with similar right quadrangle and right hexagon (see
Fig.2 and Fig.3 in \cite{GN5} respectively).
The singular modular form
$\Delta_{1/2}$ is connected with the right $\infty$-gon with vertices at
infinity and defines the generalized Kac--Moody superalgebra with
infinitely many simple real roots (i.e. of parabolic type).  In \S 5.1 we
describe relations between the modular forms constructed in \S 1--\S 4
and generalized Cartan matrices $A \in {\goth A}$.

Many modular forms constructed in \S 1--\S 4 (in particular, all
modular forms of $\Delta$-type) give discriminant automorphic forms
for moduli spaces of algebraic $K3$ surfaces with some special Picard
lattices and give Mirror Symmetry for K3 surfaces in the variant we
considered in \cite{GN3} (see also \cite{GN6}). 
We describe these cases in Sect. 5.2.

Below we summarize properties of the cusp form $\Delta_1(Z)$ (and
other $\Delta$-functions as well)
as the $3$-dimensional generalization
of the Dedekind $\eta$-function:

\smallskip
(a) $\ \Delta_1(Z)$ is the unique (up to a constant) cusp form of
weight one with respect to the full paramodular group $\Gamma_3$.  (We
remark that the weight $\frac{1}2$ is the singular weight for
$Sp_4(\bq)$, thus one is the minimal weight for which a cusp form may
exist.)

\smallskip
(b) $\ \Delta_1(Z)$ is a root of order six from the first cusp form
(with trivial character) with respect to $\Gamma_3$.

\smallskip
(c) $\ \Delta_1(Z)$ satisfies the Euler type identity (the infinite
sum in \thetag{0.2} is equal to the product in \thetag{0.3}) and it is
the denominator function of a generalized Kac--Moody superalgebra.

\smallskip
(d) The divisor of $\ \Delta_1(Z)$ is discriminant of the moduli space
of $K3$ surfaces with the Picard lattice
described in Sect. 5.2.

\smallskip
{\it Acknowledgment. }
This paper was written during our stay at RIMS of Kyoto University.
We are grateful to the Mathematical Institute for hospitality
and its excellent working atmosphere.
\smallskip
The previous  version of the paper was published as the preprint 
RIMS-1122 (1996) (see also alg-geom/9611028).

\head
1. The lifting of Jacobi forms of half-integral indices
\endhead

In this chapter we consider a new variant of the arithmetic lifting
(or the Maass or Saito--Kurokawa lifting) of Jacobi forms to the Siegel
modular forms (see \cite{M1}--\cite{M2} and
\cite{G1}--\cite{G5}).  The Siegel modular forms
under construction are
modular forms with respect to a full paramodular group
(integral symplectic groups of  a skew-symmetric form)
with  a character (or a multiplier system).

In what follows we consider three types of automorphic forms of
integral and half-integral weight:  modular forms for
$SL_2(\bz)$, Siegel modular forms with respect to a paramodular
group and Jacobi forms of integral and half-integral index.

Let
$$
\bh_n=\{Z={}^tZ\in M_n(\bc),\ Z=X+iY,\ Y>0\}
$$
be the Siegel upper-half space of genus $n$.  We denote by $|_k$
($k\in \bz/2$) the
standard slash operator on the space of functions on $\bh_n$:
$$
(F|_kM)(Z):=\hbox{det\,}(CZ+D)^{-k}F(M<Z>)
\tag1.1
$$
where
$$
M=\pmatrix A&B\\C&D\endpmatrix\in Sp_n(\br)\quad\text{ and }
\quad M<Z>=(AZ+B)(CZ+D)^{-1}.
$$
For a half-integral $k$ we choose one of the  holomorphic square roots
by the condition $\sqrt{\hbox{det}(Z/i)}>0$ for $Z=iY\in \bh_n$.

\definition{Definition 1.1}Let $k$ be integral
or half-integral.  Let $\gm\subset Sp_n(\br)$ be a 
subgroup which contains a principal congruence subgroup. 
A modular form of weight $k$ for $\gm$ with a multiplier system 
(or a character)
$v: \gm\to \bc^{\times}$ 
is a holomorphic  function on $\bh_n$ which
satisfies the functional equation
$$
(F|_k M)(Z)=v(M)F(Z) \qquad\qquad \text{for any } M\in \gm.
$$
For $n=1$ we have to add the standard growth condition
(the holomorphicity) at the cusps of the group $\gm$.
We denote by $\fm_k(\gm, v)$ (resp. $\fn_k(\gm, v)$) the space
of all modular (resp. cusp) forms of weight $k$.
\enddefinition

\smallskip
\subhead
1. Jacobi modular forms of half-integral index
\endsubhead
\medskip
We shall consider Jacobi forms of integral and half-integral
weight. To define the corresponding multiplier system
we use the multiplier system of the $\eta$-function.
The Dedekind eta-function
$$
\eta(\tau)=
q^{\frac 1{24}}\,\prod_{n\ge 1} (1-q^n)=
\sum_{n\in \bn} \biggl(\frac{12}n\biggl)q^{{n^2}/{24}},
\tag1.2
$$
where $\tau\in \bh_1$, $q=\exp{(2\pi i \tau)}$ and
$$
\biggl(\frac{12}n\biggl)=
\cases
\hphantom{-}1\  \text{ if }\  n\equiv \pm 1 \operatorname{mod} 12\\
-1\ \text{ if }\  n\equiv \pm 5 \operatorname{mod} 12\\
\hphantom{-}0\  \text{ if }\  (n,12)\ne 1,
\endcases
$$
satisfies the functional equation
$$
\eta(\frac{a\tau+b}{c\tau+d})=
v_\eta(M)(c\tau+d)^{\frac 1{2}}\eta(\tau)
\qquad M=\pmatrix a&b\\c&d\endpmatrix\in SL_2(\bz)
$$
where $v_\eta(M)$ is a $24$th root of unity.
The Dedekind eta-function is a modular form of weight $1/2$.
(One can consider it as a modular form with respect to the double
cover $\widetilde{SL_2(\bz)}$ of $SL_2(\bz)$.  The commutator subgroup
of $\widetilde{SL_2(\bz)}$ has  index $24$.)
For arbitrary even $D$ the $\eta$-multiplier system defines the character
$v_\eta^D$ of $SL_2(\bz)$.

\proclaim{Lemma 1.2}Let $D$ be even integral. We set
$Q=\frac{24}{(24,D)}$ where
$(24,D)=\hbox{\rm g.c.d.\,}(24, D)$.
Then $v_\eta^D$ is a character of $SL_2(\bz)$ and
$$
\hbox{Ker\,}(v_\eta^D)\supset \Gamma_1(Q)=
\{M\in SL_2(\bz)\,|\, M\equiv E_2 \operatorname{mod} Q\}.
\tag{1.3}
$$
\endproclaim
\demo{Proof}The lemma  follows from the fact that
$SL_2(\bz)/[SL_2(\bz), SL_2(\bz)]\cong \bz/12\bz$
or from the exact formula
for the multiplier system $v_\eta(M)$ in terms of the
entries  of the matrix
$M=\pmatrix a&b\\c&d\endpmatrix$.  For even $D$ it is given by
$$
v_\eta^D(M)=\cases
\exp{\bigl(\frac{2\pi i D}{24}((a+d)c-bd(c^2-1)-3c))\bigr)}
&\hskip-2pt c\equiv 1\operatorname{mod}\  2\\
\vspace{2\jot}
\exp{\bigl(\frac{2 \pi i D}{24}((a+d)c-bd(c^2-1)+3(d-cd-1))\bigr)}
&\hskip-2pt d\equiv 1\operatorname{mod}\  2.
\endcases
$$
This formula gives us a better result for special $D$:
if $Q=3$, $6$
$$
\hbox{Ker\,}(v_\eta^D)\supset \Gamma_0^0(Q)=
\{M=\pmatrix a&b\\c&d\endpmatrix\in SL_2(\bz)\,
|\, b,c\equiv 0 \operatorname{mod} Q,
\ a,d\equiv \pm 1 \operatorname{mod} Q\},
$$
if $Q=12$
$$
\hbox{Ker\,}(v_\eta^D)\supset \Gamma_{\{1,5\}}(12)=
\{M=\pmatrix a&b\\c&d\endpmatrix\in SL_2(\bz)\,
|\, b,c\equiv 0 \operatorname{mod} 12,
\ a,d\equiv 1,5 \operatorname{mod} 12\}.
$$
\newline
\qed
\enddemo

Jacobi forms of half-integral indices  appear naturally
as Fourier-Jacobi coefficients of the theta-constants or
well known $Sp_4(\bz)$-cusp form  $\Delta_5(Z)$ of weight $5$
(see \thetag{1.24} and Example 1.14 below).
Let us introduce a maximal parabolic subgroup of the integral
symplectic group
$$
\gi=
\left\{\pmatrix *&0&*&*\\
                *&*&*&*\\ *&0&*&*\\ 0&0&0&* \endpmatrix
\in Sp_4(\Bbb Z)\right\}.
\tag{1.4}
$$
The Jacobi group is  defined by
$
\Gamma^J=\gi/\{\pm E_4\} \cong SL_2(\bz)\ltimes H(\bz),
$
where $H(\bz)$ is the integral Heisenberg group, which is 
the central extension
$$
0\to \bz \to H(\bz)\to \bz \times \bz \to 0.
$$
We use the following embeddings ``tilde" of $SL_2(\bz)$ and
$[\dots]$ of $H(\bz)$ in
$\Gamma^J$
$$
\hskip-1pt \widetilde{
\bigl(\smallmatrix a&b\\c&d\endsmallmatrix\bigr)}:=
\Biggl(\smallmatrix
a&0&b&0\\
0&1&0&0\\ c&0&d&0\\ 0&0&0&1
\endsmallmatrix\Biggr),\qquad
H(\Bbb Z)\cong\left\{\lambda,\mu,\kappa\in \bz\ |\
 [\lambda,\mu; \kappa]=
\Biggl(\smallmatrix
1&0&0&\hphantom{-}\mu\\
\lambda&1&\mu&\hphantom{-}\kappa\\
0&0&1&{-}\lambda\\
0&0&0&\hphantom{-}1
\endsmallmatrix\Biggr)
\right\}.
\tag{1.5}
$$
There exists the unique binary character $v_H$ of
$H(\bz)$
$$
v_H([\lambda,\mu;\kappa]):=(-1)^{\lambda+\mu+\lambda\mu+\kappa}
\tag1.6
$$
which can be extended to the Jacobi group.
\proclaim{Lemma 1.3}Let $\chi$ be a character of $SL_2(\bz)$.
Then
$$
(\chi\times v_H) (\gamma \cdot h ):= \chi(\gamma)\cdot v_H(h),
\qquad
(\chi\times \hbox{id}_H) (\gamma \cdot h ):= \chi(\gamma)
$$
are correctly defined
characters of the Jacobi group
where $\gamma\in SL_2(\bz)$ and $h\in H(\bz)$, .
\endproclaim
\demo{Proof}It is easy to see that
$$
v_H(\gamma^{-1}h\gamma)=v_H(h)
$$
for any $\gamma\in SL_2(\bz)$.
\newline
\qed
\enddemo
In what follows we fix the notations
$$
Z=\pmatrix \tau &z\\z&\omega\endpmatrix\in \bh_2,\
e(z)=\exp{(2\pi i z)},\
q=e(\tau), \  r=e(z),\  s=e(\omega).
\tag1.7
$$
\definition{Definition 1.4}Let $t$ be an integral  or
{\it half-integral} positive number.
A holomorphic  function
$\phi(\tau ,z)$ on $\bh_1\times \bc$ is called
{\it a  Jacobi form of
weight $k$ and index $t$ with a multiplier system  (or a character)}
$v:\Gamma^J\to \bc^{\times}$ if the function
$$
\widetilde{\phi}(Z):=\phi(\tau ,z)\exp{(2\pi i t\omega )},
\qquad Z\in \bh_2,
$$
is a $\gi$-modular form of weight $k$
with the multiplier system $v$,
i.e if it satisfies the functional equation
$$
(\widetilde{\phi}|_k M)(Z)=v(M)\widetilde{\phi}(Z)
\qquad\qquad \text{for any } M\in \Gamma^J
$$
and has a Fourier expansion of type
$$
\phi(\tau ,z)=\sum\Sb n,\,l\\
\vspace{0.5\jot} 4nt-l^2\ge 0\endSb
f(n,l) \exp{(2\pi i(n\tau +lz))}
$$
where the summation is taken over $n$ and $l$ from some
free $\bz$-modules. The condition $f(n,l)=0$ unless $4tn-l^2\ge 0$
is equivalent to the holomorphicity of $\phi$ at infinity.
The form $\phi(\tau ,z)$ is called a Jacobi cusp form if 
$f(n,l)=0$ unless $4tn-l^2>0$.
We call a holomorphic function $\phi(\tau,z)$
{\it a nearly holomorphic} Jacobi form of weight $k$
and index $t$ if it satisfies the functional equation of the definition
and there exists $n\in \bn$ such that
$\Delta(\tau)^n\phi(\tau,z)$ is a  Jacobi form.
We call the form  $\phi(\tau,z)$
{\it a weak} Jacobi  form if 
$f(n,l)\ne 0$ only for $n\ge 0$ in its Fourier expansion.
\enddefinition

We denote the space of all Jacobi forms (resp. cusp forms, weak forms
or nearly holomorphic Jacobi forms) with the multiplier system $v$
by
$J_{k,t}(v)$ (resp. $J_{k,t}^{cusp}(v)$, $J_{k,t}^{weak}(v)$ 
or $J_{k,t}^{nh}(v)$).
We remark that $J_{k,t}(v)$ contains only zero if $k\le 0$, but
there exist weak Jacobi forms of non-positive weights.
One can give  a similar definition for any congruence-subgroup of
the Jacobi group, but with the definition given above we would like
to emphasize
that there exists a large class of Jacobi forms of half-integral indices
with respect to  the full Jacobi group.
In what follows the two classical examples are
very important.

\example{Example 1.5} {\it Jacobi theta-series, the  Jacobi triple product
and  the quintiple product.}
The Jacobi theta-series is defined as
$$
\vartheta(\tau ,z)=\hskip-2pt\sum\Sb n\equiv 1\, mod\, 2 \endSb
\,(-1)^{\frac{n-1}2}
\exp{(\frac{\pi i\, n^2}{4} \tau +\pi i\,n z)}=
\sum_{m\in \bz}\,\biggl(\frac{-4}{m}\biggr)\, q^{{m^2}/8}\,r^{{m}/2}
\tag1.8
$$
where
$$
\biggl(\frac{-4}{m}\biggr)=\cases \pm 1 &\text{if }
m\equiv \pm 1\ \hbox{mod}\ 4\\
\hphantom{\pm}0 &\text{if }
m\equiv \hphantom{\pm} 0 \ \hbox{mod}\ 2.
\endcases
$$
It is known that $\vartheta(\tau ,z)$ satisfies the following
transformations formulae with respect to the standard generators of
$SL_2(\bz)$ and $H(\bz)$:
$$
\align
\vartheta(\tau ,z+\lambda \tau +\mu)&=
(-1)^{\lambda+\mu}\exp{(-\pi i\,(\lambda^2 \tau  +2\lambda z))}\,
\vartheta(\tau , z)\\
\vartheta(-\frac 1{\tau },\, \frac {z}{\tau })&=
\exp{(-\frac{3\pi i}4)}\sqrt{\tau }\,\exp{(\pi i\,\frac {z^2}{\tau })}\,
\vartheta(\tau ,z)
\\
\vartheta(\tau +1,\,z)&
=\exp{(\frac {\pi i}{4})}\,\vartheta(\tau ,\,z).
\endalign
$$
This is a Jacobi form of weight $1/2$ and index $1/2$ with a
multiplier system $v_\vartheta$.  By definition, $v_\vartheta$ is the
multiplier system of the function $\vartheta(\tau , z)\exp{(\pi i \omega )}$.
The first functional equation for $\vartheta(\tau , z)$ is
equivalent to
$$
\vartheta(\tau ,z)\exp{(\pi i {\omega })}|_{\frac{1}2}
[\lambda,\mu;\kappa]=
(-1)^{\lambda+\mu+\lambda\mu+\kappa}\vartheta(\tau ,z)\,
\exp{(\pi i {\omega })}.
$$
Thus the restriction of $v_\vartheta$ to $H(\bz)$ is the character $v_H$
defined in \thetag{1.6}.  Moreover we have
$$
\frac{\partial\vartheta(\tau ,z)}{\partial z}\big|_{z=0}=
2\pi i \sum\Sb n\equiv 1\, mod \, 4\\
\vspace{0.5\jot}n>0 \endSb
\biggl(\frac{-4}{n}\biggr)nq^{n^2/8}
=2\pi i
 \,\eta(\tau )^3.
$$
Thus
$$
v_\vartheta(\widetilde M)=v_\eta(M)^3 \qquad \text{for }\ M\in SL_2(\bz).
$$
The formula for
$\eta(\tau)^3$ is due to Jacobi and follows from his famous triple
formula
$$
\prod_{n\ge 1}(1-q^{n-1}r)(1-q^n r^{-1})(1-q^n)=
\sum_{m\in \bz}(-1)^m q^{\frac 1{2}m(m-1)}r^m
$$
or equivalently
$$
\vartheta(\tau ,\,z)=
-q^{1/8}r^{-1/2}\prod_{n\ge 1}\,(1-q^{n-1} r)(1-q^n r^{-1})(1-q^n).
\tag1.9
$$
We recall   the  quintiple product formula. We shall use it
in the  form
$$
\multline
\sum_{n\in \bz}
\biggl(\frac{12}n\biggl)\,q^{{n^2}/{24}}r^{{n}/2}\\
=
q^{\frac 1{24}}r^{-\frac 1{2}}
\prod_{n\ge 1}(1+q^{n-1}r)(1+q^{n}r^{-1})(1-q^{2n-1}r^2)
(1-q^{2n-1}r^{-2})(1-q^n)
\endmultline
$$
where the generalized Kronecker symbol was defined in
\thetag{1.2}.
The left hand side of the last identity is a Jacobi modular form.
\proclaim{Lemma 1.6}The function
$$
\vartheta_{3/2}(\tau ,z)=\sum_{n\in \bz}
\biggl(\frac{12}n\biggl)\,q^{{n^2}/{24}}r^{{n}/2}
\in J_{\frac 1{2}, \frac 3{2}}(v_\eta\times v_H)
$$
is a Jacobi modular form of weight $1/2$ and index $3/2$ with
the multiplier system $v_\eta\times v_H$.
\endproclaim
\demo{Proof}The statement follows from the identity
$$
\vartheta_{3/2}(\tau ,z)
=\frac{\eta(\tau )\vartheta(\tau , 2z)}
{\vartheta(\tau , z)}
$$
and 
$
\vartheta(\tau , 2z)\in
J_{\frac 1{2}, 2}(v_\eta^3\times \hbox{id}_H)
$.
One can prove this directly or using the Hecke operator
$\Lambda_2$  which will be defined in the next section.
\newline\qed\enddemo
\endexample

\subhead
1.2. Hecke operators
\endsubhead
\medskip
We shall use some Hecke operators on the space of Jacobi modular forms.
Let us define the Hecke ring of the Jacobi
group $\Gamma^J$ (equivalently, the parabolic subgroup $\gi$ defined
in \thetag{1.4})
and its congruence subgroup
$$
\gm^J(Q)=\gq\ltimes H(\bz)
$$
where $\gq$ is the principal congruence subgroup of $SL_2(\bz)$
(see \thetag{1.3}). We denote the corresponding
subgroup of $\Gamma_\infty$ by  $\Gamma_\infty(Q)$.
An element of the Hecke ring is  a formal finite sum
of double left cosets with respect to $\Gamma_\infty$
(about the Hecke rings of this type see 
\cite{G6}--\cite{G8}, \cite{G11}). If
$$
X=\sum_i a_i\gi N_i\gi=\sum_j b_j \gi M_j\in H(\gi),
$$
we define in a standard way
its action on the space of all $\gi$-modular forms 
$$
(F|_k X)(Z):=
\cases
\sum_j  \nu(M_j)^{2k-3} b_j (F|_k M_j)(Z)&\quad\text{if } \ k\ne 0\\
\sum_j b_j (F|_0 M_j)(Z)&  \quad\text{if } \ k= 0,
\endcases
$$
where $\nu(M_j)$ is the degree of the simplectic similitude $M_j$.
We use the normalizing factor $\nu(M)^{2k-3}$ connected with
$Sp_4$ since  the corresponding Hecke operators are also used
in the construction of some $L$-functions
(e.g. see \cite{G3}, \cite{G11} for  the $Spin$-$L$-function of Siegel 
modular forms and \cite{G8}--\cite{G9} for the skew-symmetric square of 
the standard $L$-function of $SU(2,2)$-modular forms).

Let us recall the definition of the Hecke element $T^{(Q)}(m)$ of
the ring $H(\gq, M_2^+(Q))$
where
$$
M_2^+(Q)=\{\pmatrix a&b\\c&d\endpmatrix\in M_2^+(\bz)\,|
\pmatrix a&b\\c&d\endpmatrix\equiv
\pmatrix 1&0\\0&1\endpmatrix \ \hbox{mod Q}\}.
$$
For  $(m,Q)=1$ we have
$$
T^{(Q)}(m)=\sum\Sb ad=m\\ \vspace{ 0.5\jot} b\, mod\, d\endSb
\gq \sigma_a \pmatrix a&Qb\\0&d\endpmatrix,
$$
where $a>0$ and  $\sigma_a\in SL_2(\bz)$ such that
$\sigma_a\equiv
\left(\smallmatrix a^{-1}&0\\0&a\endsmallmatrix\right)
\hbox{ mod}\ Q$.

We shall use two different types of Hecke operator on the space of
Jacobi forms $J_{k,t}(\chi)$:
$$
{}\hskip-5pt
\Lambda_n=\gi(Q)\,\hbox{diag}(1,n,1,n^{-1})\,\gi(Q)=
\gi(Q)\,\hbox{diag}(1,n,1,n^{-1}) \in H(\gi(Q))
\tag1.10
$$
and
$$
T^{(Q)}_-(m)=\sum\Sb ad=m\\
\vspace{ 0.5\jot} b\, mod\, d\endSb
\gi(Q)\, \widetilde{\sigma_a}
\pmatrix a&0&Qb&0\\0&m&0&0\\
0&0&d&0\\0&0&0&1\endpmatrix\in H(\gi(Q))
\tag1.11
$$
where $(m,Q)=1$.
These elements realize the embeddings of the
semigroup $\bn^{\times}$ (equivalently, the Hecke ring
$H(\{1\}, \bn^{\times})$ of the
trivial group) and the Hecke ring of $\gq$ into
the Hecke
ring of the parabolic subgroups $\gi(Q)\subset \gi$.
One can consider $H(\gi)$ as a non-commutative extension of
the usual commutative Hecke ring $H(Sp_4(\bz))$
 (see the papers of the first
author mentioned above for a general case of
 of the Hecke rings of parabolic subgroups of the classical groups
over local fields and \cite{G2} for the case of the paramodular groups).

Let us suppose that the $SL_2(\bz)$-part of the character $\chi$
of Jacobi group  has the conductor $Q$.
Then element \thetag{1.11} defines the following operator on the
spaces of Jacobi forms of integral weight $k$ and index $t$
($t\in \bz/2$) with
the character $\chi$:
$$
\widetilde{\phi}|_k T^{(Q)}_-(m)(Z)=
m^{2k-3}\sum\Sb ad=m\\ \vspace {0.5\jot} b\, mod\, d\endSb
d^{-k}\chi(\sigma_a)\phi(\frac{a\tau +bQ}d,\,az)
\exp{(2\pi i mt\omega)}.
\tag{1.12}
$$
The element $\Lambda_n$ defines a Hecke  operator for
an arbitrary (integral
or half-integral) weight $k$
$$
(\widetilde{\phi}|_k \Lambda_n)(Z)=n^{-k}\phi(\tau , nz)
\exp{(2\pi i\, n^2 t\omega )}.
\tag1.13
$$
Usually  we shall omit the variable $\omega$ in the action of
Hecke operators.

\proclaim{Lemma 1.7}Let
$\phi(\tau ,z)\in J_{k,t}(\chi\times v_H^\varepsilon)$, where $\chi$
is a character of $SL_2(\bz)$ and $\varepsilon =0$ or $=1$.  Then
$$
\phi|_k \Lambda_n\in J_{k,tn^2}(\chi\times v_H^{n\varepsilon}).
$$
We assume that $k$ is integral, $\Gamma_1(Q)\subset
\hbox{Ker}\,(\chi)$ and  $(m,2^{\varepsilon}Q)=1$.
Then
$$
(\phi|_k\, T^{(Q)}_-(m))(\tau ,z)\in J_{k,mt}(\chi_m\times
v_H^{\varepsilon})
$$
where $\chi_m$ is a character of $SL_2(\bz)$ defined by
$$
\chi_{m}(\alpha):=\chi(\alpha_m)
$$
with $\alpha_m\in SL_2(\bz)$ such that
$
\alpha_m\equiv \bigl(\smallmatrix 1&0\\0&m\endsmallmatrix\bigr)
\alpha
\bigl(\smallmatrix 1&0\\0&m^{-1}\endsmallmatrix\bigr)
\operatorname{mod} Q
$.
\endproclaim
\remark{Remark}It is clear from the definition that
the space $J_{k,t}(\chi\times v_H^\varepsilon)$ is not empty only if
$2t\equiv \varepsilon \mod 2$. 
\endremark
\demo{Proof}It is easy to check that
$$
(\phi|_k \Lambda_n)|_k[\lambda,\mu;\kappa]=
(-1)^{m\lambda+m\mu+m^2\lambda\mu+m^2\kappa}\phi|_k \Lambda_n.
$$
This proves the first statement.
Any
representative from the system indicated in \thetag{1.11} acts on
$\omega $ by multiplication on $m$.
Thus the form
$$
\widetilde\psi(Z)=\psi(\tau ,z)e^{2\pi i t m\omega }=
(\widetilde\phi|_k T^{(Q)}_-(m))(Z)
$$
is a Jacobi form of index  $mt$.  Let us find its character.
For any $\alpha\in SL_2(\bz)$
$$
\widetilde{\psi}|_k\widetilde\alpha
=\sum_{M\in T^{(Q)}_-(m)}
\widetilde{\phi}|_k(\widetilde{\alpha_M}\cdot M)
$$
(see \thetag{1.5}).
By definition of the system of representatives $\{M\}$
given in \thetag{1.11},
$$
\alpha_M\bigl(\smallmatrix 1&0\\0&m\endsmallmatrix\bigr)\equiv
\bigl(\smallmatrix 1&0\\0&m\endsmallmatrix\bigr)
\alpha\operatorname{\,mod\,} Q.
$$
This gives the formula for the character $\chi_m(\alpha)$.
If index $t$ is half-integral, then for arbitrary
$[\lambda,\mu;\kappa]\in H(\bz)$ we have
$$
\psi|_k\widetilde{\sigma_a}
\Biggl(\smallmatrix a&0&Qb&0\\0&m&0&0\\
0&0&d&0\\0&0&0&1\endsmallmatrix\Biggr)
[\lambda,\mu;\kappa]=
(-1)^{d\lambda+a\mu+m\lambda\mu -Qb\lambda(d\lambda+1)+m\kappa}
\phi|_k \widetilde{\sigma_a}
\Biggl(\smallmatrix a&0&Qb&0\\0&m&0&0\\
0&0&d&0\\0&0&0&1\endsmallmatrix\Biggr) .
$$
This finishes the proof.
\newline
\qed
\enddemo
Using the exact formula for $v_\eta^D$ from Lemma 1.2, we obtain

\proclaim{Lemma 1.8}Let $\chi$ be a character of $\gi$ of
type $v_\eta^D\times v_H^\varepsilon$, where $D$ is even.
Let $Q=24/(D,24)$ and $m\equiv -1\operatorname{mod}\, Q$.
Then
$$
\chi_{m}=\overline{\chi}
$$
where bar denotes complex conjugation.
\endproclaim
\smallskip
\subhead
1.3. Paramodular groups and Humbert modular surfaces
\endsubhead
\medskip

It is well known that the  symplectic group of rank four over a field
is isomorphic to  an  orthogonal group of signature $(3,2)$.
We recall the corresponding construction over $\bz$.
We fix a $\bz$-module
$L=\Bbb Z e_1\oplus \Bbb Z e_2 \oplus \Bbb Z e_3 \oplus \Bbb Z e_4$.
Arbitrary $\bz$-linear map $g:\,L\to L$ induces the linear map
${\bigwedge}^2 g: L\wedge L\to L\wedge L$ of the $\bz$-module
$L\wedge L$
of integral bivectors.
$L\wedge L$ is  isomorphic to the module
of integral skew-symmetric matrices ($e_i\wedge e_j$ corresponds
to the elementary skew-symmetric matrix $E_{ij}$ having only two
non-zero elements $e_{ij}=1$ and $e_{ji}=-1$).
The scalar product $(u,v)$ on $L\wedge L$ is defined by
$u\wedge v=(u,v)e_1\wedge e_2\wedge e_3\wedge e_4\in \wedge^4 L$. This is
an even unimodular integral symmetric bilinear
form of signature $(3,3)$ on $L\wedge L$
which
is invariant with respect to the action of $SL(L)$ on
$L\wedge L$. We recall that
if the matrix $X$ corresponds to the bivector
$x=\sum_{i<j}x_{ij}e_i\wedge e_j$, then
$(x,x)=2\,\hbox{Pf}\,(X)$, where $\hbox{Pf}\,(X)$
is the Pfaffian of $X$ and
$\hbox{Pf}\,(MX{}^tM)=\hbox{det}(M)\,\hbox{Pf}\,(X)$.

We fix a skew-symmetric form $J_t$ on $L$ by the property:
$$
J_t(x,y)e_1\wedge e_2\wedge e_3\wedge e_4=
-x\wedge y\wedge w_t,\qquad w_t=te_1\wedge e_3+ e_2\wedge e_4.
$$
The group
$$
\widetilde{\Gamma}_t = \{g:L\to L\ |\  J_t(gx,gy)=J_t(x,y)\}
$$
is called {\it the integral paramodular group of level} $t$.
The lattice $L_t=w_t^{\perp}$ consisting  of all elements of $L\wedge L$
orthogonal to $w_t$ has the basis
$$
f_1=e_1\wedge e_2,\
f_2= e_2\wedge e_3,\
f_3= te_1\wedge e_3 -e_2\wedge e_4,\
f_{-2}=e_4\wedge e_1,\
f_{-1}=e_4\wedge e_3.
\tag{1.14}
$$
The Pfaffian  defines a  quadratic form
$S_t$ of signature $(3,2)$ on the lattice $L_t$ which has
the following matrix
$$
S_t=\left(\smallmatrix
0&0&0&0&-1\\
0&0&0&-1&0\\
0&0&2t&0&0\\
0&-1&0&0&0\\
-1&0&0&0&0
\endsmallmatrix\right)
$$
in the basis $\{f_i\}$.
It gives a homomorphism from the integral paramodular  group
to the orthogonal group  of  the lattice $L_t$
$$
{\tsize\bigwedge}^2:\, \widetilde{\Gamma}_t \to \operatorname{O}(L_t).
$$
The paramodular group $\widetilde{\Gamma}_t$ is conjugate
to a subgroup of the usual rational symplectic group $Sp_4(\bq)$:
$$
\Gamma_t:=I_t \widetilde{\Gamma}_t I_t^{-1}=\left\{\pmatrix
*    &   t*   &   *   &   *\\
*   &    *   &  *   &   *t^{-1}\\
{*}   &  t *   &   *   &   *\\
t{*}  &  t*& t*   &   *
\endpmatrix
 \in Sp_4 (\bq)\,|\, \text{ all $*$ are integral}\, \right\},
\tag{1.15}
$$
where
$I_t=\hbox{diag}\,(1,1,1,t)$. Therefore we get a homomorphism
$$
\Phi: \Gamma_t\to \operatorname{O}(L_t),\qquad
\Phi(g)={\tsize\bigwedge}^2(I_t^{-1}\, g \,I_t).
\tag{1.16}
$$
One can check, that
$$
\Phi(\Gamma_t)\subset \widehat{\hbox{\rm SO}}(L_t)
=\widehat{\hbox{\rm O}}(L_t)\cap \operatorname{SO }(L_t),\qquad
\hbox{\rm Ker}\,\Phi=\{\pm E_4\},
$$
where
$$
\widehat{\hbox{O}}(L_t)=\{g\in \operatorname{O}(L_t)\,|\
\forall \,\ell\in \widehat L_t\quad g\ell-\ell\in L_t\}
$$
($\widehat L_t$ denotes the dual lattice)
is the subgroup of the orthogonal group consisting of elements
acting identically on the discriminant group
$$
A_{L_t}:=\widehat L_t/L_t = (2t)^{-1} \bz/\bz \cong \bz/2t \bz
$$
equipped with the corresponding  finite quadratic form $q_{L_t}$
(see \cite{N1}).
The next lemma is well known (see, for example, \cite{GH1})

\proclaim{Lemma 1.9}The homomorphism $\Phi$ defines an isomorphism
between $\Gamma_t$ and
$\widehat {SO}^+(L_t)={SO}^+(L_t)\cap \widehat{SO}(L_t)$,
where ${SO}^+(L_t)$ denote the subgroup of elements with real spinor
norm equals one.
\endproclaim

The paramodular group has  normal extensions generated by the
$\Phi$-preimages of the elements of $\operatorname{SO}^+(L_t)$.
{}From the result of \cite {N1}, it  follows that
$SO^+(L_t)/
\widehat{SO}^+(L_t)
\cong O(q_{L_t})$
where $O(q_{L_t})$
is the finite orthogonal group of  the discriminant group
$(A_{L_t}, q_{L_t})$.
For the  case of the lattice $L_t$
$$
O(q_{L_t})=
\{b\ \hbox{mod}\  2t\,|\, b^2\equiv 1\  \hbox{mod}\  4t\}
\cong(\bz/2\bz)^{\nu(t)}
\tag{1.17}
$$
where  $\nu(t)$ is the number of prime divisors of $t$.
Using \thetag{1.16}, one can define a normal extension
$\Gamma_t^*$ of the paramodular group
$$
\Gamma_t^*=\ <\Gamma_t, V_d,\  \text{ where }\  d||t >
\ \subset Sp_4(\Bbb R),
$$
where for  a strict divisor $d$ of $t$
($d||t$ means that $d|t$ and $(d, \frac{t}d)=1$)
we set
$$
V_d=(\sqrt{t})^{-1}{\widetilde V}_d,\quad
{\widetilde V}_d=\Biggl(\smallmatrix dx&-t&0&0\\
         -y&d&0&0\\
          0&0&d&y\\
          0&0&t&xd\endsmallmatrix
\Biggr)\in GSp_4(\bz)
\quad
(xd-y\frac td=1).
\tag{1.18}
$$
Then we have  $\Gamma_t^*/\Gamma_t\cong O(q_{L_t})$.
For a square free $t$ the group $\Gamma_t^*$ is a maximal discrete
group acting on $\bh_2$.
The homomorphism ${\tsize\bigwedge}^2$ defines the isomorphism
$$
\Phi: \operatorname{P}\Gamma_t^* \to PO^+(L_t).
$$
We shall use also a double normal extension
$$
\Gamma_t^+=\Gamma_t\cup \Gamma_tV_t, \qquad
V_t={\tsize\frac 1{\sqrt{t}}}\left(\smallmatrix
0   &   t   &   0   &   0\\
1  & 0 &   0   &   0\\
0   &   0   &   0   &   1\\
0   &   0   &   t   &   0
\endsmallmatrix\right).
\tag{1.19}
$$
One can prove (see \cite{G1, Lemma 2.2})
\proclaim{Lemma 1.10}
The group $\Gamma_t^+$ is generated by the maximal
parabolic subgroup $\Gamma_{\infty, t}=\Gamma_t\cap \gi(\bq)$
and $V_t$.
\endproclaim

We shall use  the exact images of the generators of $\Gamma_t^+$
with respect to $\Phi$:
$$
\aligned
\Phi(\widetilde{\left(\smallmatrix
 a&b\\c&d\endsmallmatrix\right)})&=
\left(\smallmatrix
a&-b&0&0&0\\
-c&d&0&0&0\\
0&0&1&0&0\\
0&0&0&a&b\\
0&0&0&c&d
\endsmallmatrix\right),\\
\vspace{2\jot}
\Phi(V_t)&=\left(\smallmatrix
-1&0&0&0&0\\
 0&0&0&-1&0\\
 0&0&-1&0&0\\
 0&-1&0&0&0\\
 0&0&0&0&-1
\endsmallmatrix\right),
\endaligned
\quad
\aligned
\Phi([\lambda,\mu;\frac{\kappa}t])&=
\left(\smallmatrix
1&0&2t\mu    &t\lambda\mu-\kappa&t\mu^2\\
0&1&2t\lambda& t\lambda^2       &t\lambda\mu+\kappa\\
0&0&1&\lambda&\mu\\
0&0&0&1&0\\
0&0&0&0&1
\endsmallmatrix\right),\\
\vspace{2\jot}
\Phi(J_t)&=\left(\smallmatrix
 0&0&0&0&-1\\
 0&0&0&-1&0\\
 0&0&1&0&0\\
 0&-1&0&0&0\\
-1&0&0&0&0
\endsmallmatrix\right).
\endaligned
$$
The involution $\Phi(V_t)$ acts on the discriminant group
$A_{L_t}$ as  multiplication on $-1$, thus we have the isomorphism
$$
\Phi: \hbox{P}\Gamma_t^+\to P\widehat{O}^+(L_t).
$$
We remark that for the case of a perfect square $t=d^2$  the group
$\Gamma_{d^2}$ is conjugate (by the matrix
$\hbox{diag}(1,d,1,d^{-1})$) to a subgroup of $\Gamma_1=Sp_4(\bz)$.
In this realization we have
$$
\Gamma_{d^2}\cong{\Gamma}_{d,1}=\{M\in Sp_4(\bz)\,|\,M\equiv
\Biggl(\smallmatrix
*   &   0   &   *   &   0\\
0  & * &   0   &   *\\
*   &   0   &   *   &   0\\
0   &   *   &   0  &   *
\endsmallmatrix\Biggr)\ \hbox{mod}\  d \}, \
V_{d^2}\cong V_{d,1}=\Biggl(\smallmatrix
0   &   1   &   0   &   0\\
1  & 0 &   0   &   0\\
0   &   0   &   0   &   1\\
0   &   0   &   1  &   0
\endsmallmatrix\Biggr).
\tag{1.20}
$$

The real orthogonal group $O^+(L_t\otimes \br)$ acts on 
the homogeneous domain $\Omega_t^+$
of type IV corresponding to the lattice $L_t$,
which is a subdomain of a quadric in the projective space
$\Bbb P^4$
$$
\align
\Omega_t^{+}&=
\{Z\in \bp(L_t\otimes \bc)\,|\ (Z,Z)=0,\ (Z,\overline Z)<0\}^+\\
{}&=
\{Z={}^t\bigl((t z_2^2-z_1z_3),\,z_3,\,z_2,\,z_1,\,1\bigr)
\cdot z_0\
|\ \hbox{Im}\,(z_1z_3-t z_2^2)>0,\ \hbox{Im}\,(z_1)>0\}\\
{}&\cong\{\frak z={}^t(z_3,z_2,z_1)\in \bc^3\,|\,
\hbox{\,Im}\,(z_1z_3-t z_2^2)>0,\ \hbox{Im}\,(z_1)>0\}=
\bh_t^+,
\endalign
$$
where $\bh_t^+$ is the $3$-dimensional tube domain of type IV
(see the part I of this paper for a general definition).
Let us define an isomorphism $\phi_t$ of the Siegel
upper half-plane $\bh_2$ with $\bh_t^+$ as
$$
\psi_t(\pmatrix \tau &z\\z&\omega \endpmatrix)=
{}^t\bigl(t\omega ,\,z,\,\tau\bigr)\in \bh_t^+.
$$
Using the formulae for $\Phi$-images of
generators of the symplectic
group given above we  obtain the commutative diagram
$$
\CD
\bh_2@>g>>\bh_2\\
@V\psi_t VV @V\psi_t  VV\\
\bh_t^+@>\Phi(g)>>\bh_t^+
\endCD
\tag{1.21}
$$
where  $g\in Sp_4(\br)$.

Let us compare the Fourier expansions of
modular forms with respect to symplectic and orthogonal
groups. For $F(Z)\in \frak M_{k}(\Gamma_t,\chi)$ we have
$$
F(Z)=\sum
\Sb N=\bigl(\smallmatrix n&l/2\\l/2&mt\endsmallmatrix
\bigr)\ge 0\endSb
a(N)\exp{(2\pi i\,\hbox{tr}(NZ))}=
\sum
\Sb N=\bigl(\smallmatrix n&l/2\\l/2&mt\endsmallmatrix
\bigr)\ge 0\endSb a(N)\,q^nr^ls^{tm}
$$
where $n$, $l$ and $m$ run over some free $\bz$-modules
depending on the character $\chi$.
If we consider $F(\phi_t(Z))$
as a modular form with respect
to the orthogonal group $\widehat {SO}^+(L_t)$
then we can rewrite the  Fourier expansion above as
$$
F(\frak z)=F(z_1f_2+z_2f_3+z_3f_3)=\hskip-2pt
\sum\Sb \ell = nf_2-l\widehat{f}_3+mf_{-2}\\
\vspace{0.5\jot} -(\ell, \ell) \ge 0\endSb
\hskip-2pt
a(\ell)\exp{\bigl(-2\pi i\,(\ell, \frak z)_{L_t}\bigr)}=
\sum\Sb \ell
\endSb
a(\ell)q^nr^ls^{m}
$$
where the summation is taken over the same $n$, $l$,
$m$ like in the symplectic Fourier expansion above,
and
$\widehat{f}_3=\frac 1{2t} f_3$,
$q=\exp{(2\pi i z_1)}$,
$r=\exp{(2\pi i z_2)}$,
$s=\exp{(2\pi i z_3)}$. Thus
$$
a(\pmatrix n&l/2\\l/2&m\endpmatrix)=
a(nf_2-l\widehat{f}_3+mf_{-2}).
$$
A primitive
$\ell=(e, a, -\frac{b}{2t}, c, f)\in \widehat{L}_t$
(here primitive means $(e,a,b,c,f)=1$) determines the quadratic rational
divisor
$$
\Cal H_{\ell}=\{Z\in \Omega_t^+ \ |\ (\ell,Z)=0\}\cong
\{
\pmatrix \tau&z\\ z&\omega \endpmatrix\in \bh_2\ |\
tf(z^2-\tau\omega)+tc\omega+bz+a\tau+e=0\,\}
\tag{1.22}
$$
of the discriminant
$D(\ell)=2t(\ell,\ell)=b^2-4tef-4tac$.
Due to the isomorphism $\Phi$, the groups $\Gamma_t$, $\Gamma_t^+$
and $\Gamma_t^*$ act
on the set of all rational quadratic divisors with
a fixed discriminant. We define a Humbert modular surface $H_\ell$
in the Siegel modular threefold
$\Cal A_t=\Gamma_t\setminus \bh_2$
($\Cal A_t^+=\Gamma_t^+\setminus \bh_2$ or
$\Cal A_t^*=\Gamma_t^*\setminus \bh_2$ respectively)
by
$$
H_\ell=\pi_t\,\bigl(\bigcup
\Sb g\in \gm_t  \endSb
g^*(\Cal H_\ell) \bigr),
$$
where $\pi_t: \bh_2\to \Cal A_t$
is the natural projection.
We remark that according to our definition a Humbert modular surface
is irreducible.
The $\widehat{SO}(L_t)$-orbit of a primitive
vector $\ell\in \widehat{L}_t$ depends only on the norm
of $\ell$ and its image in the discriminant group.
Thus any Humbert surface  in $\Cal A_t$
of discriminant $D$
can be represented in the form
$$
H_D(b)=\pi_t(\{Z\in \bh_2\,|\,a\tau+bz+t\omega=0\})
\tag{1.23}
$$
where  $a,b\in \bz$, $D=b^2-4ta$
and $b\,\hbox{mod}\, 2t$. From that follows that
the number of  Humbert surfaces in $\Cal A_t$
of discriminant
$D$ is equal to
$
\#\,\{\ b\,\operatorname{mod }\  2t\,
|\ b^2\equiv D\,\operatorname{mod }\  4t\}.
$
The involution $V_t$
acts on the discriminant group
by multiplication on $-1$.
Therefore the number of the Humbert surfaces of
discriminant $D$ in $\Cal A_t^+$ is equal to
$
\#\,\{\ \pm b\,\operatorname{mod }\  2t\,
|\ b^2\equiv D\,\operatorname{mod }\  4t\}
$.
About the theory of  Humbert surfaces for the case of
a non-principle polarization see \cite{vdG} and \cite{GH1}.

\medskip
\subhead
1.4. The ``trivial" lifting of the Jacobi theta-function
\endsubhead
\smallskip
We define the arithmetic lifting of
Jacobi forms as an action of a formal operator  $L$-function
on a given Jacobi modular form.
We may illustrate this approach with a classical example.
Let us consider the even theta-series
$$
\vartheta(\tau)=\sum_{n\in \bz}\exp{(2\pi i n^2\tau)}.
$$
The periodic function
$\exp{(2\pi i\, \tau)}$ is invariant with respect to the parabolic subgroup
$\Gamma_0
=\{\bigl(\smallmatrix \pm 1&m\\0&\pm 1\endsmallmatrix\bigr)
\,|\,m\in \bz\}\subset SL_2(\bz)$.
The element
$[n^{-1}]=\Gamma_0
\bigl(\smallmatrix n&0\\0&n^{-1}\endsmallmatrix\bigr)
\Gamma_0$  of the Hecke ring
$H(\Gamma_0)$
acts on the function $\exp{(2\pi i \tau)}$ in a very simple
way:
$\exp{(2\pi i \tau)}|[n^{-1}]=\exp{(2\pi i n^2\tau)}$.
Thus
$$
\vartheta(\tau)=1+2\sum_{n\in \bn} \exp{(2\pi i \tau)}\big|[n^{-1}]=
1+2\exp{(2\pi i \tau)}\big|\bigl(\sum_{n\in \bn}\, [n^{-1}]\bigl)
$$
where we consider $[\zeta(1)]=\sum_{n\ge 1} [n^{-1}]$ as
a formal Dirichlet series
over the Hecke ring $H(\Gamma_0)$ which has an expansion
in the formal infinite product $[\zeta(1)]=\prod_{p} (1-[p^{-1}])^{-1}$
over all primes.
We can interpret
the last representation of the theta-series as a
lifting of  $\Gamma_0$-modular form $\exp{(2\pi i\, \tau)}$
defined by  the formal zeta-function  $[\zeta(1)]$.
(See \cite{G6} for some applications
of this representation  of $\vartheta(\tau)$.)

In this section we use this kind of a ``trivial" lifting
to define two special
singular modular forms  with respect to the full paramodular group
$\Gamma_{4}$ and $\Gamma_{36}$.
The first form is one of the classical Siegel theta-constants.
The behavior of the theta-constants with respect to
the congruence-subgroups of Hecke type
is well known (see for example \cite{Fr2}).
Bellow we make an observation that
the ``most odd'' even  theta-constant
$$
\Theta_{\bold 1,\bold 1}(Z)
=
\frac {1}2\sum_{l_1,l_2\in \bz}
\exp{\bigl(\pi i\, (Z
\bigl[
\smallmatrix
l_1+\frac{1}2\\
l_2+\frac{1}2
\endsmallmatrix
\bigr]+
l_1+l_2)\bigr)}
=\frac{1}2\sum\Sb n,\,m\in \bz \endSb
\,\biggl(\dsize\frac{-4}{n}\biggr)\biggl(\dsize\frac{-4}{m}\biggr)
q^{n^2/8}r^{nm/4}s^{m^2/8},
$$
where $Z[M]={}^tMZM$,
is a modular form with respect to the subgroup
${\Gamma}_{2,1}\subset Sp_4(\bz)$ conjugated to the paramodular
group
$\Gamma_4$ (see \thetag{1.20}).
Our  second example is a  modular form of weight $1/2$
with respect to the full paramodular group $\Gamma_{36}$.
\proclaim{Theorem 1.11}
1. The function
$$
\Delta_{1/2}(Z)=
\frac{1}2\sum\Sb n,\,m\in \bz \endSb
\,\biggl(\dsize\frac{-4}{n}\biggr)\biggl(\dsize\frac{-4}{m}\biggr)
q^{n^2/8}r^{nm/2}s^{m^2/2}
$$
is a modular form of weight $1/2$
with respect to the double extension
$\Gamma_4^+$ of the paramodular group
$\Gamma_4$(see \thetag{1.19})
with a multiplier system
$v_{8}:\Gamma_4^+\to \{\root{8}\of{1}\}$
which is induced by $v_\eta^3\times v_H$. It means, that
$$
v_{8}|_{SL_2(\bz)}=v_\eta^3,\quad
v_{8}|_{H(\bz)}=v_H,\quad
v_{8}([0,0;\frac{\kappa}4])=\exp{(\frac{\pi i \kappa}4)}
\ \  (\kappa\in \bz)
$$
and $\Delta_{1/2}(V_4(Z))=\Delta_{1/2}(Z)$.
Moreover the divisor of  $\Delta_{1/2}(Z)$ on
$\Cal A_4^+={\Gamma}_4^+\setminus \bh_2$ is exactly equal  to
the  Humbert modular surface of discriminant $1$
$$
\hbox{Div}_{\Cal A_4^+}(\Delta_{1/2}(Z))=H_1=
\pi_4^+(\{Z\in \bh_2\,|\,z=0\}).
$$
2. The function
$$
D_{1/2}(Z)=\frac{1}2 \sum_{m,\, n\in \bz}
\biggl(\frac{12}n\biggl)\biggl(\frac{12}m\biggl)
\,q^{{n^2}/{24}}r^{{nm}/2}s^{{3m^2}/{2}}
$$
is a modular form of weight $1/2$
with respect to the normal extension
$\Gamma_{36}^*=[\Gamma_{36}, V_{4}, V_{9}, V_{36}]$
of order $4$ of $\Gamma_{36}$ with a multiplier system
$v_{24}:\Gamma_{36}^*\to
\{\root {24} \of {1}\}$  induced by $v_\eta\times v_H$, i.e.
$$
v_{24}|_{SL_2(\bz)}=v_\eta,\quad
v_{24}|_{H(\bz)}=v_H,\quad
v_{24}([0,0;\frac{\kappa}{36}])=\exp{(\frac{\pi i \kappa}{12})},
\ (\kappa\in \bz)
$$
and
$$
D_{1/2}(V_{36}(Z))=D_{1/2}(V_{9}(Z))
=D_{1/2}(V_{4}(Z))=D_{1/2}(Z).
$$
Moreover  the divisor of $D_{1/2}(Z)$
on  the threefold
$\Cal A_{36}^*={\Gamma}_{36}^*\setminus \bh_2$
consists of the  three irreducible components with multiplicity one
$$
\hbox{Div}_{\Cal A_{36}^*}(D_{1/2})=
H_4^*+H_9^*(27)+H_{16}^*(32),
$$
where $H_4^*=\pi_{36}^*(\{Z\in \bz_2\,|\,2z-1=0\})\ $ and
$$
H_9^*(27)=\pi_{36}^*(\{5\tau+27z+36\omega=0\}),\qquad
H_{16}^*(32)=\pi_{36}^*(\{7\tau+32z+36\omega=0\}).
$$
\endproclaim
\demo{Proof}Let us define a lifting
of the Jacobi theta-series
$\widetilde{\vartheta}(Z)=\vartheta(\tau,z)\exp{(\pi i \omega)}$
using a formal  Dirichlet $L$-series over the Hecke ring
$H(\gi)$
$$
L_-(\frac 1{2}, \biggl(\frac{-4}{\cdot}\biggr))
=\sum_{m\ge 1} m^{-\frac {1}2}\biggl(\frac{-4}{m}\biggr)\,\Lambda_m
$$
(see \thetag{1.10}, \thetag{1.13}).
After action of the formal operator Dirichlet series
on $\widetilde\vartheta(Z)$, we get
$$
\Delta_{1/{2}}(Z)=
\sum_{m>0} m^{-\frac1{2}}\biggl(\dsize\frac{-4}{m}\biggr)
\,\widetilde{\vartheta}\big|_{\frac 1{2}}\Lambda_m(Z)
=\sum\Sb n\in \bz,\, m\in \bn \endSb
\,\biggl(\dsize\frac{-4}{n}\biggr)\biggl(\dsize\frac{-4}{m}\biggr)
q^{n^2/8}r^{nm/2}s^{m^2/2}.
$$
Thus $\Delta_{1/{2}}(Z)$
is invariant with  weight $k$ and with the  multiplier
system mentioned in the theorem with respect to the parabolic
subgroup
$\Gamma_{\infty,4}=\Gamma_4\cap \gi(\bq)$ of $\Gamma_4$.
Moreover we see from this representation that
$\Delta_{1/{2}}(Z)$ is invariant
with respect to the transformation $q\to s^4$, $s\to q^{\frac 1{4}}$
(i.e.  $\tau \to 4\omega $, $\omega \to {\tau }/4$)
defined by $V_4$ (see \thetag{1.19}).
According to Lemma 1.10 the $\Delta_{1/{2}}(Z)$ is a modular form
with respect to
$\Gamma_4^+$. The last group is conjugate to
$\Gamma_{2,1}^+$ (see \thetag{1.20}). Thus the function
$$
\Theta_{\bold 1,\bold 1}(Z)=\Delta_{\frac 1{2}}
(\pmatrix \tau&z/2\\ z/2&
{\omega}/4\endpmatrix)
$$
is a modular form with respect to the subgroup $\Gamma_{2,1}$
of $Sp_4(\bz)$.
{}From the construction of $\Delta_{1/2}(Z)$
as a lifting  of the triple product, it
follows that its divisor contains the Humbert surfaces $H_1$.
To prove that this is  the full  divisor of
$\Delta_{1/2}(Z)$, let us consider a
$Sp_4(\bz)$-modular form
$$
F_5(Z)=\prod_{g\in \Gamma_{2,1}^+\setminus \Gamma_1}
(\Theta_{\bold 1,\bold 1}|_{\frac 1{2}} g) (Z).
\tag{1.24}
$$
Since $[\Gamma_1:\Gamma_{p,1}^+]=\frac{p^4+p^2}2$, the modular form
$F_5(Z)$ has weight $5$. This is a cusp form, because
$F_5(\bigl(
\smallmatrix \tau &z\\z&\omega\endsmallmatrix
\bigr))\equiv 0$ for $z=0$.
Up to a constant there exists only one $Sp_4(\bz)$-cusp form
of weight $5$ whose divisor is equal to $H_1$
(see \cite{Fr2}).
From this fact, it easily follows that the divisor of
$\Delta_{1/2}(Z)$ on $\Cal A_4^+$ is exactly $H_1$.

The construction of  $D_{1/2}(Z)$ is similar. The form
$$
D_{1/2}(Z)
=
\sum_{m>0} m^{-\frac 1{2}}\biggl(\dsize\frac{12}{m}\biggr)
\,\widetilde{\vartheta}_{3/2}|_{\frac 1{2}}\Lambda_m (Z)=
\sum_{n\in \bz,\, m\in \bn}
\biggl(\frac{12}n\biggl)\biggl(\frac{12}m\biggl)
\,q^{{n^2}/{24}}r^{{nm}/2}s^{{3m^2}/{2}}
$$
is modular with respect to the parabolic subgroup
$\Gamma_{\infty,36}=\Gamma_{36}\cap \gi(\bq)$ with a
multiplier system $v_{\eta}\times v_H$
and is invariant with respect to the action of $V_{36}$.
The involutions $V_9$ and $V_4$ satisfy the relation
$\Gamma_{36}(V_9V_4)=\Gamma_{36}V_{36}$.
Taking $V_9$ in the form \thetag{1.18} we
see that the involution $V_9$
does not change the Fourier expansion. Thus $D_{1/2}(Z)$ is
a modular form with respect to $\Gamma_{36}^*$.

The threefold $\Cal A_{36}^*$ contains the only  Humbert surface
$H^*_4$ of discriminant $4$ which is a part of the divisor
of $D_{1/2}(Z)$ since $\vartheta_{3/2}{(\tau, z)}\equiv 0$ if
$z=\tsize\frac 1{2}$.
We remark that $\Cal A_{36}^+=\Gamma_{36}^+\setminus \bh_2$
contains two Humbert surfaces
of discriminant $4$:
$$
H_4^+(2)=\pi_{36}^+(\{2z-1=0\})\quad\text{and}\quad
H_4^+(34)=\pi_{36}^+(\{8\tau+34z+36\omega=0\}).
$$
Let us prove that $D_{1/2}(Z)$ is anti-invariant with respect
to the involutions which define the surfaces
$H_9^*(27)$ and $H_{16}^*(32)$. For that we consider
this function
as a modular form with respect to the orthogonal group
$PO(L_{36})$ (see Sect. 1.3)
$$
D_{1/2}(\frak z)=
\frac 1{2}\sum_{n,m\in \bz}
\biggl(\frac{12}n\biggl)\biggl(\frac{12}m\biggl)
\exp{\bigl(2\pi i\, (\frac {n^2 }{24}z_1 +\frac{nm }2 z_2
+\frac{m^2 }{24} z_3))}.
$$
The Humbert surface $H_{9}^*(27)$ is defined by the
element
$l_{9}=(0, 5, \frac{27}{72},1, 0)\in \widehat{L}_{36}$
(see \thetag{1.21} and \thetag{1.23}). In the basis
$(f_2,f_3,f_{-2})$ defined in \thetag{1.14},
the involution $\sigma_{l_9}$ has the matrix
$$
(\sigma_{l_9})=\pmatrix
81&-24\cdot 90& 400\\
6&-161&30\\
16&-24\cdot 18&81\endpmatrix.
$$
Thus
$$
\gather
 D_{1/2}(\sigma_{l_9}(\frak z))=\\
=\sum_{n\in \bz,\, m\in \bn}
\biggl(\frac{12}n\biggl)\biggl(\frac{12}m\biggl)
\exp{\bigl(2\pi i (\tsize\frac {(20m-9n)^2}{24} z_1+
\tsize\frac{(20m-9n)(9m-4n) }2 z_2
+\tsize\frac{(9m-4n)^2 }{24} z_3))}\\
=-D_{1/2}(\frak z)
\endgather
$$
and $D_{1/2}(Z)$ has zero along the Humbert surface
$H_{9}^*(27)$.
For the surface $H_{16}^*(32)$ the proof is similar.

Using the  ``projection" of type \thetag{1.24} of the modular form
$D_{1/2}
(\bigl(\smallmatrix \tau& z/6\\ {z}/6&{\omega}/36
\endsmallmatrix\bigr))$ with respect to the group
$\Gamma_{6,1}$
to the space of modular forms with respect
to $Sp_4(\bz)$, we can prove that its  divisor consists
of three Humbert surfaces $H_4^*$, $H_{9}^*(27)$
and $H_{16}^*(32)$ with multiplicities one,
because the weight of the  $Sp_4(\bz)$-modular form,
whose divisor is exactly  $H_{d^2}$,
is known (see \cite{GN4} for
the  construction of such functions). In \S 2 we shall prove
the statement about the divisor of $D_{1/2}(Z)$ using another
method, which also  gives us
infinite product expansions of $\Delta_{1/2}(Z)$
and $D_{1/2}(Z)$ (see \thetag{2.11} and \thetag{2.14}).
\newline
\qed
\enddemo

\medskip
\subhead
1.5.
The arithmetic  lifting of Jacobi forms of half-integral index
\endsubhead
\smallskip
The next theorem is a generalization of  a construction
of cusp  forms with respect to the paramodular groups
proposed in \cite{G1} and \cite{G3}. This is an advanced
version of Maass lifting (see \cite{M1}--\cite{M2}).
The main construction is valid for $SO(n,2)$-group
(see \cite{G4}--\cite{G5}).

\proclaim{Theorem 1.12} Let $k$ be integral, $t$ be integral
or half-integral, $D$ be an even divisor of $24$. We take
the conductor $Q={24}/D$ and  $\mu\in (\bz/Q\bz)^*$.
Let us consider a Jacobi
cusp form
$\phi\in J_{k,t}^{cusp}(v_\eta^D\times v_H^{\varepsilon})$,
where
$\varepsilon=0$ or $\varepsilon =1$.
Then the function
$$
F_\phi(Z)=\hbox{Lift}_\mu(\phi)(Z)
=\sum\Sb m\equiv \mu \, mod \,Q\\
\vspace{0.5\jot} m>0\endSb
m^{2-k}\bigl(\widetilde{\phi}|_k T_-^{(Q)}(m)\bigr)(Z)
$$
(see \thetag{1.12})
is a cusp form with respect to a paramodular group,
and
$$
\align
\hbox{Lift}_\mu(\phi)(Z)&\in
\frak N_k(\Gamma_{Qt}^+, \chi_{D,\varepsilon,\mu})
\qquad(\text{if }\  tQ\in \bz)\\
\hbox{Lift}_\mu(\phi)(Z)&\in
\frak N_k({\Gamma'}_{4Qt}, \chi_{D,1,\mu})
\quad(\text {if }\  tQ \text{ is half-integral)}
\endalign
$$
where
${\Gamma'}_{4Qt}=\delta_2{\Gamma}_{4Qt}^+\delta_2^{-1}$
$(\delta_2=\hbox{diag}(1,2,1,2^{-1}))$.
If $\hbox{Lift}_\mu(\phi)\not\equiv 0$, then
the $\chi_{D,\varepsilon,\mu}$ is a character of order $Q$ or $2Q$
of the group $\Gamma_{Qt}^+$ or ${\Gamma'}_{4Qt}$ respectively.
This character is induced by
$v_{\eta, \mu}^D\times v_H^{\varepsilon}$,
where $v_{\eta,\mu}^D$ is a character of $SL_2(\bz)$
$\mu$-conjugated to $v_{\eta}^D$ (see Lemma 1.7), and by the relations
$$
\chi_{D,\varepsilon,\mu}(V_{Qt})=(-1)^k,\qquad
\chi_{D,\varepsilon,\mu}([0,0;\frac{\kappa}{Qt}])=
\exp{(2\pi i\, \frac{\mu\kappa}Q)}\quad (\kappa\in \bz).
$$
If $\mu=1$, then $\hbox{Lift}_1(\phi)(Z)\not\equiv 0$
for $\phi\not\equiv 0$, i.e. we have an embedding of the space
$J_{k,t}^{cusp}(v_\eta^D\times v_H^{\varepsilon})$ into
the space of Siegel modular forms.
\endproclaim
\remark{Remark 1}We denote $\hbox{Lift}_1$ by $\hbox{Lift}$.
If $\mu\ne 1$,  the lifting
can be zero for some non-zero  Jacobi forms.
\endremark
\remark{Remark 2}It is easy to give a variant of this theorem for
non-cusp Jacobi forms. One should only add  a natural condition on
weights and  characters of $SL_2$-Eisenstein series in order to get
convergence of the lifting  (compare with \cite{G4}).
It gives interesting examples of Eisenstein series
with respect to the paramodular groups. We hope to describe
them somewhere else.
\endremark
\remark{Remark 3}In the notation
$\frak N_k(\Gamma_t, \chi_{D,\varepsilon,\mu})$
we will put the corresponding
character $v_{\eta,\mu}^D\times v_H^{\varepsilon}$
of the Jacobi group  instead of the induced character
of $\Gamma_{Qt}^+$.
\endremark

\demo{Proof}
The convergence of the series
defining $\hbox{Lift}_\mu(\phi)$ follows from
the  upper bound
of  Jacobi cusp forms
of weight $k$ and index $t$ on $\Bbb H_1\times \Bbb C\,$:
$$
|\phi(\tau, z)|< C {v}^{-\frac k2}\exp (2\pi ty^2/v),
$$
where
$v=\hbox{Im}\,\tau>0$, $y=\hbox{Im}\, z$
and  the constant $C$ does not depend on $\tau$ and $z$.
To prove the last inequality, we take the function
$$
\phi^*(\tau,z)=v^{\frac k2}\hbox{exp}(-2\pi ty^2/v)|\phi(\tau,z)|
$$
which is $\Gamma^J$-invariant and
is bounded on any compact subset in $\Bbb H_1\times
\Bbb C$.
If we take the following realization of
the fundamental domain $D$
of $\Gamma^J$ on $\Bbb H_1\times \Bbb C$
$$
\Cal D=\{(\tau, \alpha\tau+\beta)\,|\,
-1\le \alpha, \beta\le 1,\  \tau\in SL_2(\Bbb Z)\setminus \Bbb H_1\},
$$
then the function $\phi^*$ is bounded on the set
$\{(\tau,z)\in \Cal D, \ \hbox{Im\,}\tau>C\}$
because
$\phi^*(\tau,z)\to 0$ as $v\to \infty$
for any cusp form $\phi(\tau,z)$.

Let us consider the Fourier expansion of the Jacobi form
$\phi(\tau,z)$ with the character $v_\eta^D\times v_H^{\varepsilon}$
$$
\align
\phi(\tau  ,\,z )&=\sum\Sb n\equiv D\,mod\,24\\
\vspace{0.5\jot} l\equiv \varepsilon \,mod\,2 \\
\vspace{0.5\jot}n> 0,\,  4nt> l^2\endSb
f(n,l)\,\exp{\bigl(2\pi i (\frac n{24}\tau +\frac l{2}z)\bigr)}\\
{}&=\sum\Sb N\equiv 1\,mod\,Q\\
\vspace{0.5\jot} l\equiv \varepsilon \,mod\,2 \\
\vspace{0.5\jot}N> 0,\,  4NDt> l^2\endSb
f(ND,l)\,\exp{\bigl(2\pi i (\frac N{Q}\tau +\frac l{2}z)\bigr)}.
\endalign
$$
By the definition \thetag{1.12} of the Hecke operators, we have
$$
\gather
m^{2-k}\bigl(\widetilde{\phi}|_k \,T_{-}^{(Q)}(m)
\bigr)(Z)=
\\
m^{k-1}\sum_{ad=m}d^{-k}
v_\eta^D(\sigma_a)
\sum
\Sb N\equiv 1\,mod\,Q\\
\vspace{0.5\jot} l\equiv \varepsilon \,mod\,2 \\
\vspace{0.8\jot} b\,mod\, d
\endSb
f(DN,l)\,
\exp{\bigl(2\pi i (\frac {N(a\tau +bQ)}{Qd}+\frac {al}{2}z+mt\omega )
\bigr)}\\
=
\sum_{ad=m}a^{k-1}v_\eta^D(\sigma_a)
\sum
\Sb dN_1\equiv 1\,mod\,Q\ \\
\vspace{0.5\jot} l\equiv \varepsilon \,mod\,2
\endSb
f(dN_1D,l)\,\exp{\bigl(2\pi i
(\frac{aN_1}{Q}\tau +\frac {al}{2}z+mt\omega )\bigr)}.
\endgather
$$
If $m\equiv \mu\ \hbox{mod }Q$, then
the condition $dN_1\equiv 1\,\operatorname{mod} Q$
($N=dN_1$) is equivalent
to $aN_1\equiv \mu\,\operatorname{mod} Q$, because
$ad\equiv \mu\,\operatorname{mod} Q$,
$Q$ is a divisor of $24$
and for arbitrary $d$ with  $(d,24)=1$
it is true that $d^2\equiv 1\,\operatorname{mod} 24$.
Hence taking the summation over all
$m\equiv \mu\,\operatorname{mod} Q$,
we get
$$
\gather
\sum_{m\equiv \mu\,mod\,Q}m^{2-k}(\widetilde{\phi}|_k
T_-^{(Q)}(m))(Z)=
\\
\sum
\Sb
a>0,\  d>0\\ \vspace{0.5\jot} ad\equiv \mu\,mod\,Q
\endSb
a^{k-1}v_\eta^D(\sigma_a)
\sum
\Sb aN_1\equiv \mu\,mod\,Q,\ \\
\vspace{0.5\jot} l\equiv \varepsilon \,mod\,2
\endSb
f(dN_1D,\,l)\,\exp{(2\pi i (\frac {aN_1}{Q}\tau +\frac {al}{2}z+adt\omega ))}
\\
=\sum
\Sb N,\,M>0\\
\vspace{0.5\jot} N,M\equiv \mu\,mod\,Q\\
\vspace{0.5\jot} L\equiv \varepsilon \,mod\,2
\endSb
\biggl( \sum_{a|(N,L,M)}
a^{k-1}v_\eta^D(\sigma_a)
f(\frac{NMD}{a^2},\, \frac{L}{a})\biggr)
\exp{(2\pi i (\frac N{Q}\tau +\frac {L}{2}z+Mt\omega ))}.
\tag{1.26}
\endgather
$$
Let us suppose that $Qt$ is integral.
Then either $\varepsilon=0$ or $\varepsilon=1$ and $Q$ is even.
In both cases
$\widetilde{\phi}|_k\,T_-^{(Q)}(m)$ is a $\gi$-modular form with
character $v_{\eta,\mu}^D\times v_H^\varepsilon$ for all
$m\equiv \mu\,\operatorname{mod} Q$ (see Lemma 1.8).
Moreover
$$
\widetilde{\phi}|_k\,T_-^{(Q)}(m)\big|_k[0,0;\frac \kappa{Qt}]=
\exp{(2\pi i \frac {\kappa\mu}{Q})}\,\widetilde{\phi}|_kT_-^{(Q)}(m).
$$
Thus the action of the center of the parabolic subgroup
$\Gamma_{\infty, Qt}=\Gamma_{Qt}\cap \gi(\bq)$
on $F_\phi(Z)$ is given by the character
$$
\nu_Q ([0,0;\frac{\kappa}{Qt}])=\exp{(2\pi i\,
\frac{\mu\kappa}Q)}
$$
of order $Q$. Therefore the lifting $F_\phi(Z)$ is
a $\Gamma_{\infty,Qt}$-modular form with character
$v_{\eta,\mu}^D\times v_H^\varepsilon\times \nu_Q$.
The Fourier expansion calculated above shows  that
$F_\phi(Z)$ is invariant
under the  change of the variables $\{\tau \to Qt\omega ,\
z\to z,\ \omega \to (Qt)^{-1}\tau \}$. This transformation is made
by the element  $V_{Qt}$ (see \thetag{1.19}). Thus
$$
(F_{\phi}|_k\, V_{Qt}) (Z)=(-1)^kF_{\phi}(Z).
$$
Moreover we have  $F_{\phi}|_k\, J_{Qt}=F_{\phi}$, where
$$
J_{Qt}=\left(\smallmatrix 0&0&1&0\\
             0&0&0&(Qt)^{-1}\\
             -1&0&0&0\\
             0&-Qt&0&0
\endsmallmatrix\right)
$$
is the standard element from the  group  $\Gamma_{Qt}$, since
$$
V_{Qt}IV_{Qt}I=J_{Qt} \qquad \text{where}\ \
I=\left(\smallmatrix 0&0&1&0\\
           0&1&0&0\\
           -1&0&0&0\\
             0&0&0&1
\endsmallmatrix\right)\in \gi.
$$
This proves that for arbitrary
$M\in \Gamma_{Qt}^+=\Gamma_{Qt}\cup \Gamma_{Qt}V_{Qt}$
$$
(F_\phi|_k M ) (Z)=\chi(M) F(Z)
$$
where $\chi(M)^{Q}=1$ for $M\in \Gamma_t$.
If $F_\phi(Z)\not \equiv 0$,
then we have a character
$$
\chi: \Gamma_{tQ}^+\to \bc^*
$$
of the paramodular group extending the character
$v_{\eta,\mu}^D\times v_H^\varepsilon$ because for $k\in \bz$
the operator $\  |_k$ defines the action of the group
on the space of functions.
We denote this character $\chi_{D,\varepsilon,\mu}$.
Thus $F_\phi(Z)$ is a $\Gamma_{Qt}^+$-modular form with this character.
If $\mu=1$, then the lifting is an injective embedding
of the space of Jacobi forms into the corresponding space
of $\Gamma_t$-modular forms, because the first
Fourier-Jacobi coefficient
of $F_\phi(Z)$ is $\phi(\tau,z)\not\equiv 0$.
This finishes the proof in the case $tQ\in \bz$.

Let us suppose that $tQ$ is half-integral. Then $\varepsilon=1$
and $Q$ is odd, i.e. $Q=3$. Let us consider the Jacobi form
$$
\psi(\tau,z)=\phi(\tau ,2z)=
2^{-{k}}(\phi|_k \Lambda_2)(\tau ,z)
\in J_{k,4t}(v^D_{\eta}\times 1_H).
$$
Its lifting $F_{\psi}(Z)$ is a modular form of weight $k$
with respect to the group $\Gamma_{4Qt}$ with a character induced by
the character $v^D_{\eta,\mu}\times \hbox{id}_H$.
The operators $T_-^{(Q)}(m)$ and $\Lambda_2$
commute. Thus
$$
F_{\psi}(Z)=2^{k}\sum_{m\equiv \mu \, mod \,Q}
m^{2-k}(\widetilde{\phi}|_k T_-(m)|_k\Lambda_2)(Z)=
(F_{\phi}|_k\delta_2)(Z)
$$
where $\delta_2=\hbox{diag}(1,2,1,2^{-1})$.
It proves that $F_{\phi}$ is a modular form with respect to
$\delta_2^{-1}\Gamma_{4Qt}\delta_2$ with a character of order $3$
if $Qt$ is half-integral.

The proof that the lifting is a cusp form
is the same as for $Q=1$ (see \cite{G1}).
\newline
\qed
\enddemo

\subhead
1.6. Examples of the arithmetic lifting
\endsubhead
\smallskip

Below we give series of examples  of liftings of Jacobi forms.
It provides us with many important Siegel modular forms
with respect to the paramodular groups.
These modular forms have interesting
applications to the algebraic geometry and to the theory
of Lorentzian Kac--Moody algebras.
It was proved in \cite{G1} that the field of rational
functions with respect to $\Gamma_t$
(or, equivalently, the moduli space of Abelian surfaces with
polarization of type $(1,t)$) might be rational
only for the twenty exceptional polarizations
$t=1,\dots, 12$, $14$, $15$, $16$, $18$, $20$,
$24$, $30$, $36$.
This list partly explains our interest to the modular forms
of small weight with respect to $\Gamma_t$.
Another starting  point for the construction of these
modular forms is the
classification  of
hyperbolic generalized Cartan matrices in the part I
of this paper (see also \cite{GN4}).
The Siegel modular forms obtained as the arithmetic lifting
of holomorphic Jacobi forms will define  automorphic
Lorentzian Kac--Moody algebras corresponding
to these generalized Cartan matrices.
One of the advantages of the arithmetic lifting
is that we have  formula \thetag{1.26} for all  Fourier coefficients
of the lifted form in terms of Fourier coefficients
of the corresponding Jacobi form.
The only factor we have to calculate is $v_\eta^{D}(\sigma_a)=\pm 1$.
According to Lemma 1.2,
 we have
$$
v_\eta^{D}(\sigma_a)
=\cases {}\quad 1 &\text{if }\  Q_D=1,\,2,\, 3,\,6\\
\biggl(\dsize\frac{-4}a\biggr) &\text{if }\  Q_D=4,\,12
\endcases
\tag{1.27}
$$
for $Q_D=24/(24,D)$.
\example{Example  1.13} {\it The case of trivial character.}
If $\chi=\hbox{id}_{SL_2}\times \hbox{id}_H$,
we get the lifting constructed in \cite{G1} and \cite{G4}
$$
\hbox{Lift}: J_{k,t}^{cusp}\to \frak N_k(\Gamma_t).
$$
The case of $t=1$ is the original Maass lifting defined
in \cite{M1}.
It is known that $\hbox{dim} (J_{1,t})=0$ (see \cite{Sk}).
Thus the arithmetic lifting with trivial character gives us
modular forms of weight $k\ge 2$. For $k=3$ we obtain
canonical differential forms on the moduli spaces of
Abelian surfaces with $(1,t)$-polarization (see \cite{G1},
\cite{G2}) or on some finite quotients of such threefolds
(see \cite{GH1}). In the examples below we shall see that
if we admit Jacobi forms with commutator characters
(characters of type $v_\eta^D\times v_H^\varepsilon$), we can
construct roots of order $d$ ($d$ is a divisor of $12$)
from some  $\Gamma_t$-modular forms.
\endexample

\example{Example  1.14} {\it Three dimensional variants of Dedekind
$\eta$-function:} $\Delta_{1}(Z)$, $\Delta_{2}(Z)$
{\it and}
$\Delta_{5}(Z)$.
We construct a series of examples of cusp   forms using
the Jacobi forms of the type
$\eta^d(\tau)\vartheta(\tau,z)
\in J_{\frac {d+1}2, \frac 1{2}}^{cusp}(v_\eta^{d+3}\times v_H)$
where $d= 1,\  3,\  9, \ 23$.
For any $d$ mentioned above  we obtain the cusp forms
$$
\hbox{Lift}_1(\eta^d(\tau)\vartheta(\tau,z))\in
\frak N_{\frac {d+1}2}(\Gamma_{t_d}, v_\eta^{d+3}\times v_H)
$$
for the paramodular groups with $t_1=3$, $t_3=2$, $t_9=1$
and $t_{23}=6$.
Let us put
$$
\eta(\tau)^d=\sum_{n\in \bn} \tau_d(n)q^{n/{24}}.
$$
Then according to \thetag{1.26} we get an exact formula
for the Fourier coefficients of the lifting
$
\hbox{Lift}_1(\eta^d(\tau)\vartheta(\tau,z))(Z)
$
in terms of $\tau_d(n)$.
For all cases when we know an elementary expression for  $\tau_d(n)$,
we get an elementary formula for the Fourier
coefficients of the lifting.
In particular for $d=1$ and $d=3$
we obtain two very nice  cusp forms of weight $1$ and $2$
with respect to the group $\Gamma_3$ and $\Gamma_2$ respectively
with the following Fourier expansions
$$
\Delta_1(Z)=\sum_{M\ge 1}
\sum\Sb n,\,m >0,\,l\in \bz\\
\vspace{0.5\jot} n,\,m\equiv 1\,mod\,6\\
\vspace{0.5\jot} 4nm-3l^2=M^2\endSb
\hskip-4pt
\biggl(\dsize\frac{-4}{l}\biggr)
\biggl(\dsize\frac{12}{M}\biggr)
\sum\Sb a|(n,l,m)\endSb \biggl(\dsize\frac{6}{a}\biggr)
q^{n/6}r^{l/2}s^{m/2}
\in \frak N_1(\Gamma_3,v_\eta^4\times v_H)
$$
(where we use notation \thetag{1.7}) and
$$
\Delta_2(Z)=\sum_{N\ge 1}\
\sum\Sb
 n,\,m >0,\,l\in \bz\\
\vspace{0.5\jot} n,\,m\equiv 1\,mod\,4\\
\vspace{0.5\jot} 2nm-l^2=N^2
\endSb
\hskip-4pt
N\biggl(\frac {-4}{Nl}\biggr)
\sum_{a\,|\,(n,l,m)} \biggl(\frac {-4}{a}\biggr)
\, q^{n/4} r^{l/2} s^{m/2}
\in \frak N_2(\Gamma_2,v_\eta^6\times v_H).
$$
Hence  all Fourier coefficients $a(N)$ of the {\it cusp form
of weight one}
$\Delta_1(Z)$
corresponding to the primitive matrices $N$ are equal to
$\pm 1$ or $0$!
We may say that this function is the simplest Siegel cusp form.
$\Delta_1(Z)$ and $\Delta_2(Z)$ have properties similar to
the Dedekind $\eta$-function.
For example, $\Delta_1(Z)$ is a root of order $6$ of the unique
up to a constant cusp form (with trivial character) of weight
$6$ for $\Gamma_3$ and
$\Delta_2(Z)$ is  a root of order $4$ from the unique
cusp form of weight $8$ for $\Gamma_2$.
Both functions  have infinite product
expansions and they are the discriminants of  moduli spaces of
$K3$ surfaces of special types.
(See \S 5 below  and
\cite{GN2}--\cite{GN3} where  $\Delta_2$ was used.)

The modular forms $\Delta_1(Z)$ and $\Delta_2(Z)$
are  also connected with
the theory  of moduli spaces $\Cal A_t=\Gamma_t\setminus \bh_2$
of Abelian surfaces with polarization of type $(1,t)$
(without level structure).
Using  $\Delta_1(Z)$ and $\Delta_2(Z)$ we can prove the rationality
of the moduli spaces $\Cal A_3$ and $\Cal A_2$
respectively.

We would like to remark that
$\Delta_1(Z)^3dZ=\Delta_1(Z)^3d\tau \wedge dz \wedge d\omega $
defines a canonical differential form on the double cover
$$
{}^{2}\Cal A_3={}^{2}\Gamma_3\setminus \bh_2
\overset{2:1}\to \longrightarrow \Cal A_3.
$$
where ${}^{2}\Gamma_3=\hbox{Ker}(\chi_6^3)$ is a subgroup
of index $2$ of $\Gamma_3$ ($\chi_6$ is the character of
$\Delta_1(Z)$ of order $6$). For some other applications of constructed
modular forms  to the moduli spaces
of Abelian and Kummer surfaces see  \cite{G2}, \cite{GH1}--\cite{GH2}.

For $d=9$ the lifting of $\eta(\tau)^9\vartheta(\tau,z)$
is equal to  the well known Siegel cusp  form
$\Delta_5(Z)$ of weight $5$
with respect to the full
Siegel modular group $\Gamma_1=\operatorname{Sp}_4(\bz)$
with the non-trivial binary character
$$
\Delta_5(Z)=\hbox{Lift}(\eta(\tau )^{9}\vartheta(\tau ,z))=
\frac 1{64}\prod_{(a,b)}\vartheta_{a,b}(Z)
$$
where the product is taken over all even Siegel theta-constants
(compare this definition with  construction \thetag{1.24}).
This function was constructed as a lifting by Maass in \cite{M2}.
It has the Fourier expansion
$$
\Delta_5(Z)=\sum\Sb n,l,m \equiv 1\operatorname{mod}2\\
\vspace{0.5\jot} n,m>0\endSb \hskip10pt
\sum_{d|(n,l,m)} (-1)^{\frac{l+d+2}2}\,d^{4}
\tau_9\,(\frac{4nm-l^2}{d^2})\,
q^{n/2}\,r^{l/2}\,p^{m/2}.
$$
It  was proved in \cite{GN1}--\cite{GN3} that the functions $\Delta_5$
and  $\Delta_2$ are Weyl--Kac denominator functions for some
generalized Kac--Moody superalgebras.
\endexample

\example{Example  1.15}
{\it The function $\Delta_{11}(Z)$ and the $\mu$-Lifting for $\mu =-1$}.
The functions $\Delta_{1}(Z)$, $\Delta_{2}(Z)$
and $\Delta_{5}(Z)$ are connected with the first
case of Theorem 1.12 when $Qt$ is even.
Now we consider the second one when $Qt$ is odd.
Let us take Jacobi cusp forms
$$
\eta^5(\tau)\vartheta(\tau, 2z)\in
J_{3, {2}}(v_\eta^8\times \hbox{id}_H),
\qquad
\eta^{21}\vartheta(\tau, 2z)
\in J_{11, {2}}
$$
where the last Jacobi form has trivial character.
Then $Q=3$ for the first  and $Q=1$ for the second Jacobi form,
thus
$$\align
\Delta_{11}(Z)&=\hbox{Lift}(\eta(\tau )^{21}\vartheta(\tau ,2z))
\in \frak N_{11}(\Gamma_{2}),
\tag{1.28}\\
F_3^{(6)}(Z)&=\hbox{Lift}(\eta(\tau )^5\vartheta(\tau ,2z))
\in \frak N_3(\Gamma_6, v_\eta^8\times \hbox{id}_H).
\endalign
$$
(We denote by $F_{k}^{(t)}(Z)$ a modular form of weight $k$
with respect to the  paramodular group $\Gamma_t$.)
According to Theorem 1.12 and Lemma 1.8 for $Q=3$
($\mu\equiv 2\equiv -1 \hbox{\, mod 3}$)
$$
\hbox{Lift}_{2}(\eta^5(\tau)\vartheta(\tau,2z))\in
\frak N_3(\Gamma_6, v_\eta^{16}\times \hbox{id}_H),
$$
if it does not vanish identically. To prove this,  we calculate
$$
\bigl(\eta^5(\tau)\vartheta(\tau,2z)|_3 \,T_-^{(3)}(2)\bigr)(Z)=
q^{\frac 2{3}}r^{-3}(1-r^2)(1-r)^4+\dots \not\equiv 0.
$$
In fact we can prove that
$$
\bigl(\eta^5(\tau)\vartheta(\tau,2z)|_3 \,T_-^{(3)}(2)\bigr)(Z)=
\eta(\tau)\vartheta^4(\tau,z)
\vartheta(\tau,2z)\in J_{3,4}(v_\eta^{16}\times \hbox{id}_H).
$$
The product of two ``conjugated" liftings is a cusp form with
trivial character
$$
F_6^{(6)}(Z)=\hbox{Lift}_1\bigl(\eta(\tau)^5\vartheta(\tau,2z)\bigr)
\cdot
\hbox{Lift}_{2}\bigl(\eta(\tau)^5\vartheta(\tau,2z)\bigr)\in
\frak N_6(\Gamma_6).
$$
\endexample

The construction of the arithmetic lifting provides us
with some information about divisor of lifted forms.
Using Theorem 1.12 and \thetag{1.12} one can prove

\proclaim{Lemma 1.16}Let $\hbox{Lift}_\mu(\phi)$ be the lifting
defined in Theorem 1.12. Let us assume that
$\phi(\tau,z)\big|_{z=0}\equiv 0$, and this zero has
order $m$. Then $\hbox{Lift}_\mu(\phi)$
has zero of order at least  $m$ along the Humbert
surface $H_1(0)$,
where $H_1(0)=\pi_{Qt}^+(\{Z\,|\,z=0\})\subset \Cal A_{Qt}^+$
(see Sect. 1.3).

More generally, let $M\subset\bq/\bz$
be a subset which is invariant with respect to multiplication
on arbitrary $a\in \bz$ such that $(a,Q)=1$. Let us assume that
$\phi(\tau,z)\big|_{z=\alpha}\equiv 0$ of order $m$
for all $\alpha\in M$.
Then $\hbox{Lift}_\mu(\phi)$
has zero of order at least  $m$ along the Humbert
surface $H(\alpha)=\pi_{Qt}^+(\{Z\,|\,z=\alpha\})
\subset \Cal A_{Qt}^+$ for any $\alpha\in M$.
\endproclaim

In particular, for the form $F_6^{(6)}(Z)$ given above
$$
\hbox{Div}_{\Cal A_6^+}(F_6(Z))\supset 2H_1+2H_4.
$$
\example{Examples 1.17} {\it The quintiple product.}
Let us consider the quintiple product
$$
\vartheta_{3/2}(\tau,z)=
\sum_{n\in \bz}
\biggl(\frac{12}n\biggl)\,q^{{n^2}/{24}}r^{{n}/2}
\in J_{\frac 1{2}, \frac 3{2}}(v_\eta\times v_H)
$$
(see Lemma 1.6).
Then for $d=1,\  3,\  5,\  11$
we can construct the following Siegel cusp forms
$$
\operatorname{Lift}\,(\eta(\tau)^{d}\vartheta_{3/2}(\tau,z))\in
\frak N_{\frac {d+1}2}(\Gamma_{t_d}, v_\eta^{d+1}\times v_H)
$$
for the paramodular groups of level
$t_1=18$, $t_3=9$, $t_5=6$ and $t_{11}=3$ respectively.
For $d=7$ and $23$ we define
$$
\align
\operatorname{Lift}(\eta(\tau)^7\vartheta_{3/2}(\tau ,2z))
&\in
\frak N_{4}(\Gamma_{18}, v_\eta^{8}\times {\hbox{id}}_H),\\
\operatorname{Lift}(\eta(\tau)^{23}\vartheta_{3/2}(\tau ,2z))
&\in
\frak N_{12}(\Gamma_{6}).
\endalign
$$
Like in Example  1.14 we may calculate the exact form of the Fourier
expansion of such forms. The  case of
$d=1$ and $d=3$, when we get the cusp forms of weight $1$ and $2$,
are of the  special interest.
They are
$$
D_1(Z)=\sum_{M\ge 1}
\sum\Sb n,\,m >0,\,l\in \bz\\
\vspace{0.5\jot} n,\,m\equiv 1\,mod\,12\\
\vspace{0.5\jot} 2nm-l^2=M^2\endSb
\hskip-4pt
\biggl(\dsize\frac{12}{Ml}\biggr)
\sum\Sb a|(n,l,m)\endSb \biggl(\dsize\frac{-4}{a}\biggr)
q^{n/12}r^{l/2}s^{3m/2}
\in
\frak N_{1}(\Gamma_{18}, v_\eta^{2}\times v_H),
\tag{1.29}
$$
$$
D_2(Z)=\sum_{N\ge 1}
\sum\Sb n,\,m >0,\,l\in \bz\\
\vspace{0.5\jot} n,\,m\equiv 1\,mod\,6\\
\vspace{0.5\jot} 4nm-l^2=N^2\endSb
\hskip-4pt
N\biggl(\dsize\frac{-4}{N}\biggr)\biggl(\dsize\frac{12}{l}\biggr)
\sum\Sb a|(n,l,m)\endSb \biggl(\dsize\frac{6}{a}\biggr)
q^{n/6}r^{l/2}s^{3m/2}
\in
\frak N_{2}(\Gamma_{9}, v_\eta^{4}\times v_H).
\tag{1.30}
$$
\endexample

Next two examples show  how one can use the Hecke operators
$\Lambda_n$ to construct new Jacobi cusp  forms of half-integral index.

\proclaim{Lemma 1.18} Let $a, b\in \bn$ be such that
$(a,b)=1$. Then
$$
\align
\vartheta(\tau ,az)\vartheta(\tau ,bz)
&\in J_{1, \frac 1{2} ({a^2+b^2})}^{cusp}(v_\eta^6\times v_H),
\qquad\  \text{if }\  ab \text{ is even,}\\
\vartheta_{3/2}(\tau ,az)\vartheta_{3/2}(\tau ,bz)
&\in J_{1, \frac{3}{2}(a^2+b^2)}^{cusp}(v_\eta^2\times v_H^{a+b}),
\quad\  \text{ if }\  (ab, 6)\ne 1,\\
\vartheta(\tau ,az)\vartheta_{3/2}(\tau ,bz)&
\in J_{1, \frac 1{2}(a^2+3b^2)}^{cusp}(v_\eta^4\times v_H^{a+b}),
\quad \text{ if }\  (a, 3)=1 \vee (a,2)=2 \vee (b,6)\ne 1.
\endalign
$$
\endproclaim
\demo{Proof}
Let us consider the Fourier expansion
$$
\vartheta(\tau ,az)\vartheta(\tau ,bz)=
\sum_{n,m\in \bz}
\biggl(\dsize\frac{-4}{n}\biggr)\biggl(\dsize\frac{-4}{m}\biggr)
q^{\frac{n^2+m^2}8} r^{\frac{an+bm}2}.
$$
For an arbitrary Fourier coefficient $f(N,L)$ the norm of
its index
$$
2(a^2+b^2)N-L^2=\frac 1{4}(am-bn)^2\ge 0
$$
can be zero only if $n$ or $m$ is even,
because $ab$ is even.
Since
$\biggl(\dsize\frac{-4}{2n}\biggr)=0$, we get
the first statement of the lemma.
The proof of the second and third statement is similar.
\enddemo

\example{Example  1.19}{\it Cusp forms of weight one.}
The  Jacobi cusp forms of weight $1$
from the last lemma generate a series
of Siegel cusp  forms of weight one:
$$
\align
\hbox{Lift}(\vartheta(\tau ,az)\vartheta(\tau ,bz))
&\in \frak N_{1}(\Gamma_{2(a^2+b^2)},\,v_{\eta}^6\times v_{H}),\\
\hbox{Lift}(\vartheta_{3/2}(\tau ,az)\vartheta_{3/2}(\tau ,bz))
&\in \frak N_{1}(\Gamma_{18(a^2+b^2)},\,v_{\eta}^2\times v_{H}^{a+b}),\\
\hbox{Lift}(\vartheta(\tau ,az)\vartheta_{3/2}(\tau ,bz))
&\in \frak N_{1}(\Gamma_{3(a^2+3b^2)},\,v_{\eta}^4\times v_{H}^{a+b}).
\endalign
$$
\endexample
\example{Example  1.20}
{\it $\vartheta(\tau,z) \vartheta(\tau,2z)$-lifting series}.
We can use Jacobi cusp forms of weight one
to produce Siegel modular forms
of small weights with respect to $\Gamma_{t}$
of small level $t$. We give some examples
for the groups $\Gamma_4$ -- $\Gamma_{10}$
in the case of $\vartheta(\tau,z) \vartheta(\tau,2z)$.
First we have
$$
\hbox{Lift\,}(\vartheta(\tau,z) \vartheta(\tau,2z))
\in \frak N_{1}(\Gamma_{10}, v_\eta^6\times v_H).
\tag{1.31}
$$
We may combine the last Jacobi cusp form  together with
$\eta(\tau)$ or $\vartheta(\tau,z)$.
For example we can get
the cusp form  of minimal weight
and  with trivial character for $\Gamma_5$
$$
F_5^{(5)}(Z)=\hbox{Lift\,}(\eta(\tau)^3
\vartheta(\tau,z)^6\vartheta(\tau,2z))
\in \frak N_{5}(\Gamma_{5}).
\tag{1.32}
$$
According to Lemma 1.16
$$
\hbox{Div}_{\Cal A_5^+}(F_5(Z))\supset 7H_1+H_4
$$
where $H_d$ are the corresponding Humbert surfaces
in the threefold
$\Cal A_5^+=\Gamma^+_5\setminus \bh_2$.
We also obtain
$$
\aligned
\hbox{Lift\,}(\eta^6\vartheta(\tau,z) \vartheta(\tau,2z))&
\in \frak N_{4}(\Gamma_{5},v_\eta^{12}\times v_H),\\
\hbox{Lift\,}(\vartheta^3\cdot \vartheta(\tau,2z))&
\in \frak N_{2}(\Gamma_{7},v_\eta^{12}\times v_H),\\
\hbox{Lift\,}(\eta^9\vartheta^4\cdot
\vartheta(\tau,2z))&
\in \frak N_{7}(\Gamma_{4}),
\endaligned
\quad
\aligned
\hbox{Lift\,}(\eta^3\vartheta^2\cdot \vartheta(\tau,2z))&
\in \frak N_{3}(\Gamma_{6},v_\eta^{12}\times \hbox{id}_H),\\
\hbox{Lift\,}(\eta^3\vartheta\cdot \vartheta(\tau,2z)^2)&
\in \frak N_{3}(\Gamma_{9},v_\eta^{12}\times v_H),\\
\hbox{Lift\,}(\vartheta\cdot \vartheta(\tau,2z)^3)&
\in \frak N_{2}(\Gamma_{13},v_\eta^{12}\times v_H).
\endaligned
$$
\endexample

In the lemma below we consider a construction which is similar
to the quintiple product (see Lemma 1.6).
\proclaim{Lemma  1.21}The   functions
$$
\align
\phi_{1,4}(\tau ,z)&=
\eta(\tau )^2\,\frac{\vartheta(\tau ,3z)}{\vartheta(\tau ,z)}
\in J_{1,4}^{cusp}(v_\eta^2\times \hbox{id}_H),\\
\psi(\tau ,z)&=
\eta(\tau )^3\,\biggl(\frac{\vartheta(\tau ,3z)}
{\vartheta(\tau ,z)}\biggr)^2
\in J_{\frac 3{2},8}(v_\eta^3\times \hbox{id}_H)
\endalign
$$
are holomorphic Jacobi forms.
\endproclaim
\demo{Proof}Let us consider a weak Jacobi form of weight $0$ and
index $4$
$$
\phi_{0,4}(\tau ,z)=
\frac{\vartheta(\tau ,3z)}{\vartheta(\tau ,z)}
=r+1+r^{-1}+q(-r^4-r^3+r+2+r^{-1}-r^{-3}-r^{-4})+\dots .
\tag{1.33}
$$
The Fourier coefficient $f(n,l)$ of $\phi_{0,4}$ depends only
on the ``norm" $16n-l^2$ and $\pm l \hbox{ mod } 8$. Moreover
$16n-l^2\ge -16$ if $f(n,l)\ne 0$. From the exact form of
coefficients with $q^0$ it follows that
$16n-l^2\ge -1$ if $f(n,l)\ne 0$, and
$f(n,l)=1$ if $16n-l^2=-1$.
Thus $\phi(\tau ,z)$ is a cusp form.
The condition $16n-l^2\ge -1$ implies  holomorphicity
of $\psi(\tau ,z)$.
\newline\qed
\enddemo

\example{Example  1.22}The lifting of the Jacobi cusp form $\phi_{1,4}$
considered in the last lemma gives us  Siegel cusp forms
of weight $1$, $2$ and $4$
$$
\gather
\hbox{Lift\,}(\phi_{1,4})\in
\frak N_{1} (\Gamma_{48}, v_\eta^2\times \hbox{id}_H),\qquad
\hbox{Lift\,}(\eta^2\phi_{1,4})\in
\frak N_{2} (\Gamma_{24}, v_\eta^4\times \hbox{id}_H),\\
\hbox{Lift\,}(\eta^6\phi_{1,4})\in
\frak N_{4} (\Gamma_{12}, v_\eta^6\times \hbox{id}_H).
\endgather
$$
All these modular forms have zero along the Humbert surface $H_9$
in the corresponding moduli spaces $\Cal A_t^+$.
\endexample

One can use  some differential operators to construct new Jacobi forms
with a commutator character.
The following lemma is a direct reformulation of
Eichler--Zagier construction.
\proclaim{Lemma 1.23}(See \cite{EZ, Theorem 9.5}.)
Let
$$
\phi_1\in J_{k_1,m_1}(v_\eta^{d_1}\times v_H^{\varepsilon_1}),
\qquad
\phi_2\in J_{k_2,m_2}(v_\eta^{d_2}\times v_H^{\varepsilon_2})
$$
be two Jacobi forms
of integral or half-integral indices,
where $\varepsilon_i=0$ or $1$. Then one can define the Jacobi form
$$
[\phi_1,\phi_2]=
\frac{1}{2\pi i}(m_2\phi_1'\phi_2-m_1\phi_1\phi_2')
\in J_{k_1+k_2+1,m_1+m_2}(v_\eta^{d_1+d_2}\times
v_H^{\varepsilon_1+\varepsilon_2}),
$$
where
$$
\phi'(\tau ,z)=\pd {\phi(\tau ,z)}{z}.
$$
\endproclaim
\demo{Proof}One should consider the Jacobi function
$\psi=\phi_1^{m_2}/\phi_2^{m_1}$ of index zero. Then
$$
\left(\frac{\phi_2^{m_1+1}}{\phi_1^{m_2-1}}\right)
\pd{\psi(\tau ,z)}{z}=
(m_2\phi_1'\phi_2-m_1\phi_1\phi_2')\in
J_{k_1+k_2+1,m_1+m_2}(v_\eta^{d_1+d_2}\times v_H^{\varepsilon(m_1+m_2)}).
$$
\newline\qed
\enddemo
We can combine two Jacobi theta-functions
$\vartheta(\tau,z)$ and $\vartheta_{3/2}(\tau,z)$ using
the differential  operator of Lemma 1.23.

\proclaim{Lemma 1.24}
The Jacobi  form
$$
\phi_{2,2}(\tau ,z)=2[\vartheta(\tau ,z), \vartheta_{3/2}(\tau ,z)]
\in J_{2, {2}}(v_\eta^4\times \hbox{id}_H)
$$
is a cusp form with integral Fourier coefficients and
$$
\frac {\phi_{2,2}(\tau ,z)}{\eta(\tau )}\in
 J_{\frac 3{2}, {2}}(v_\eta^3\times \hbox{id}_H).
$$
\endproclaim
\demo{Proof}
We  can calculate its  Fourier expansion
$$
\phi_{2,2}(\tau ,z)=\frac {1}2
\sum_{m\,,n\in \bz}
(3m-n)\biggl(\frac{-4}m\biggl)\biggl(\frac{12}n\biggl)\,
q^{\frac{3m^2+n^2}{24}}r^{\frac{m+n}2}.
\tag{1.34}
$$
For all Fourier coefficients $f(N,L)$ of this Jacobi form
we have $8N-L^2=\frac {1}{12}(3m-n)^2\ge \frac{1}3$.
The Fourier expansion shows that
$\phi_{2,2}/\eta$ is still holomorphic.
\enddemo
\remark{Remark}One can check that $[\phi_1,\phi_2]$ in Lemma 1.23
is a Jacobi cusp form
(maybe, identical to zero in some cases).
\endremark

The  Jacobi form constructed in the last lemma, produces
liftings
$$
\gather
\operatorname{Lift}(\phi_{2,2})\in
\frak N_{2}(\Gamma_{12}, v_\eta^{4}\times {\hbox{id}}_H),\qquad
\operatorname{Lift}(\eta^2\phi_{2,2})\in
\frak N_{3}(\Gamma_{8}, v_\eta^{6}\times {\hbox{id}}_H),\\
\operatorname{Lift}(\eta^4\phi_{2,2})\in
\frak N_{4}(\Gamma_{6}, v_\eta^{8}\times {\hbox{id}}_H),
\qquad
\operatorname{Lift}(\eta^8\phi_{2,2})\in
\frak N_{6}(\Gamma_{4}, v_\eta^{12}\times {\hbox{id}}_H),\\
\operatorname{Lift}(\eta^{20}\phi_{2,2})\in
\frak N_{12}(\Gamma_{2}).
\endgather
$$
One can rewrite the definition of $\phi_{2,2}$ using only the
Jacobi theta-series
$\vartheta(\tau,z)$
$$
\phi_{2,2}(\tau ,z)=2\eta(\tau )
\frac
{[\vartheta(\tau ,z),\vartheta(\tau ,2z)]}
{\vartheta(\tau ,z)}.
$$
It gives another formula for the non-cusp form
${\phi_{2,2}(\tau ,z)}/{\eta(\tau )}$.

\example{Example 1.25} {\it The Jacobi form} $\phi_{3,1}(\tau,z)$.
One can construct Jacobi forms of half-integral indices
using also Jacobi forms of integral indices.
For example, let us consider the quotient
$$
\phi_{3,1}(\tau ,z)=\frac {\phi_{12,1}(\tau ,z)}{\eta(\tau )^{18}}
\in  J_{3, {1}}(v_\eta^6\times \hbox{id}_H)
\tag{1.35}
$$
where
$$
\phi_{12,1}(\tau ,z)=(r+10+r^{-1})q
+(10r^2-88r-132-88r^{-1}+10r^{-2})q^2+\dots
$$
is the unique Jacobi cusp form of weight $12$ and index $1$ with integral
coprime Fourier coefficients (see \cite{EZ}).
This quotient is again a holomorphic Jacobi form since
$4n-l^2\ge 3$ for all Fourier coefficients of $\phi_{12,1}$.
The Jacobi form $\phi_{3,1}(\tau ,z)$ is   very useful
in many questions.
We have the following  lifting for $d=2$ and $d=6$
$$
\operatorname{Lift}(\eta^2\phi_{3,1})\in
\frak N_{4}(\Gamma_{3}, v_\eta^{6}\times {\hbox{id}}_H),
\qquad
\operatorname{Lift}(\eta^6\phi_{3,1})\in
\frak N_{6}(\Gamma_{2}, v_\eta^{12}\times {\hbox{id}}_H).
$$
One can use these modular forms and the modular forms from
the first example to construct all generators of the  graded rings
of symmetric  modular forms for $\Gamma_2$ and $\Gamma_3$.
We are planning to consider these type of questions in a publication
which follows.
\endexample

\head
\S~2. Infinite product expansion
\endhead
In this section we consider an exponential or Borcherds lifting
of Jacobi forms of weight $0$ in the space of  meromorphic
Siegel modular forms with respect to a paramodular group.
It will give us infinite product expansion of the modular forms
constructed in \S~1 and description of their divisors.

\subhead
2.1 Borcherds lifting
\endsubhead
We prove below a variant of the  Borcherds construction
of automorphic forms as infinite products
(see \cite{Bo6, Theorem 10.1}). We remark that   Theorem 10.1
in \cite{Bo6}
was formulated only for an unimodular lattice. The statement
and the main part of the proof are  valid for an arbitrary lattice
$L$ of signature $(n+2,2)$ with two orthogonal unimodular
isotropic planes if   one replaces in the theorem
a nearly holomorphic modular form $f(\tau)$ with
a  nearly holomorphic Jacobi form $\phi(\tau,\frak z)$ of weight
$0$ and index $1$ with respect to an anisotropic lattice
of rank $n$, and if one uses  results about generators of the orthogonal
group $\widehat{O}(L)$ proved in \cite{G4}.
In this section
we consider  the case of the full paramodular
group $\Gamma_t$.

Let
$$
\phi_{0,t}(\tau,z)=\sum_{n,l\in \bz} f(n,l)q^nr^l\in
J_{0,t}^{nh}
\tag{2.1}
$$
be a nearly holomorphic Jacobi form of weight $0$
and index $t$ (i.e. $n$ might be negative in the Fourier expansion).
We recall  the notations
$q=\exp{(2\pi i \tau)}$, $r=\exp{(2\pi i z)}$,
$s=\exp{(2\pi i \omega)}$
and
$\widetilde{\phi}_{0,t}(Z)=\phi_{0,t}(\tau,z)\exp{(2\pi i t\omega)}$.
Let
$$
\phi^{(0)}_{0,t}(z)=\sum_{l\in \bz} f(0,l)\,r^l
\tag{2.2}
$$
be the $q^0$-part of $\phi_{0,t}(\tau,z)$.
The Fourier coefficient $f(n,l)$ of $\phi_{0,t}$ depends only on
the norm  $4tn-l^2$ of $(n,l)$ and $l$ mod $2t$.
{}From the definition of nearly holomorphic forms, it follows
that the norm of indices  of non-zero Fourier coefficients
are bounded from bellow.

\proclaim{Theorem 2.1}Assume that
the Fourier coefficients of Jacobi form $\phi_{0,t}$
from \thetag{2.1} are integral.
Then  the product
$$
\ml(\phi_{0,t})(Z)=B_\phi(Z)=
q^{A}r^Bs^{C}
\prod\Sb n,l,m\in \bz\\
\vspace{0.5\jot}
(n,l,m)>0\endSb
 (1-q^nr^ls^{tm})^{f(nm,l)},
\tag{2.3}
$$
where
$$
A=\frac{1}{24}\sum_{l}f(0,l),\quad
B=\frac{1}{2}\sum_{l>0}lf(0,l),\quad
C=\frac{1}{4}\sum_{l}l^2f(0,l),
$$
and $(n,l,m)>0$ means that 
if $m> 0$, then $l$ and $n$ are arbitrary integers,
if $m=0$, then $n>0$ and $l\in \bz$
or  $l<0$ if $n=m=0$,
defines  a meromorphic modular form of weight
$\frac{f(0,0)}2$ with respect to $\Gamma_t^+$ with a character
(or a multiplier system if the weight is half-integral)
induced by $v_\eta^{24A}\times v_H^{2B}$. All divisors
of $\ml(\phi_{0,t})(Z)$ on $\Cal A_t^+$ are the Humbert modular
surfaces $H_D(b)$ of discriminant $D=b^2-4ta$ (see \thetag{1.23})
with  multiplicities
$$
m_{D,b}=\sum_{n>0}f(n^2a,nb).
\tag{2.4}
$$
Moreover
$$
B_\phi(V_t(Z))=(-1)^{D} B_\phi(Z)
\qquad\text{with}\quad
D=\sum\Sb n<0 \\
\vspace{0.5\jot} l \in \bz \endSb
\sigma_1(-n)f(n,l)
\tag{2.5}
$$
where
$ \sigma_1(n)=\sum_{d|n}d$.
\endproclaim
\remark{Remark}One obtains the same statement for any
Jacobi form  $\phi_{0,t}^{(S)}(\tau,z_1,\dots,z_n)$ of
weight $0$ and  index $t$
where $S$ is an even anisotropic quadratic form of rank $n$.
Then the function $\ml(\phi^{(S)}_{0,t})$ is a meromorphic automorphic
form  with respect to the orthogonal  group
$\widehat{SO}^+(U^2\oplus S(t))$ of signature $(n+2,2)$
where $U$ is the unimodular isotropic plane and $S(t)$ is
the  lattice of rank $n$
with the quadratic form $tS$.  
\endremark
\demo{Proof}The product \thetag{2.3}
is a particular
case of automorphic products considered in
\cite{Bo6, Theorem 5.1}. It converges for
$\hbox{det\,}(\hbox{Im}(Z))>c$, where $c$
is sufficiently large, and it can be
extended to a multi-valued meromorphic function on $\bh_2$
whose singularities, including zeros,
lie on rational quadratic divisors $\Cal H_\ell$ of type
\thetag{1.22}.

Let us decompose the  product of Theorem 2.1
in two factors
$$
B_\phi(Z)=
q^{A}r^Bs^{C}
\prod\Sb (n,l,0)>0\endSb
(1-q^nr^l)^{f(0,l)}\times
\prod\Sb n,l,m\in \bz\\
\vspace{0.5\jot}m>0
\endSb
(1-q^nr^ls^{tm})^{f(nm,l)}
\tag{2.6}
$$
and let us calculate the Fourier expansion of the logarithm
of the second factor:
$$
\gather
\hbox{log\,}\biggl(
\prod\Sb n,l,m\in \bz\\
\vspace{0.5\jot} m> 0\endSb
(1-q^nr^ls^{tm})^{f(nm,l)}\biggr)=
-\sum\Sb n,l\in \bz\,,\,m> 0\endSb
f(nm,l)\ \sum_{e\ge 1}\frac 1{e}\,q^{en}r^{el}s^{emt}\\
{}= -
\sum\Sb a,b,c\in \bz\\
\vspace{0.5\jot}c> 0\endSb
\ \bigl(\sum\Sb d|(a,b,c)\endSb
d^{-1}f(\frac{ac}{d^2},\,\frac{b}{d})\bigr)
\,q^ar^bs^{tc}.
\endgather
$$
The last sum can be written as the action of
the formal Dirichlet series
$\sum_{m\ge 1}m^{-1}T_{-}(m)$ on the Jacobi form
$\phi_{0,t}$ (see Theorem 1.12). Thus we obtain
$$
\hbox{log\,}\big(
\prod\Sb n,l,m\in \Bbb Z
\vspace{0.5\jot} m> 0\endSb
\dots\bigr)=-\sum_{m\ge 1}
m^{-1}\bigl(\widetilde{\phi}_{0,t}\,|\,T_-(m)\bigr)(Z).
$$
This expansion shows us that  the second factor in \thetag{2.6}
is invariant with respect to the action of the parabolic subgroup
$\Gamma_{\infty,t}$ whenever the product converges.

It is easy to see that the first factor in \thetag{2.6}
is equal to a product of  Jacobi theta-series and Dedekind
eta-functions
$$
q^{A}r^Bs^{C}
\prod_{(n,l,0)>0}
(1-q^nr^l)^{f(0,l)}=
\eta(\tau)^{f(0,0)}
\prod_{l>0}
\biggl(\frac{\vartheta(\tau,lz)e^{\pi i l^2\omega}}
{\eta(\tau)}\biggr)^{f(0,l)}.
$$
The last  identity explains the form of the factor
$q^{A}r^Bs^{C}$ in the definition
of the function of the theorem.
Thus we proved that
$$
\multline
q^{A}r^Bs^{C}
\prod\Sb n,l,m\in \bz\\
\vspace{0.5\jot}(n,l,m)>0\endSb
 (1-q^nr^ls^m)^{f(nm,l)}\\
=
\eta(\tau)^{f(0,0)}
\prod_{l>0}
\biggl(\frac{\vartheta(\tau,lz)e^{\pi i l^2\omega}}
{\eta(\tau)}\biggr)^{f(0,l)}
\exp{\bigl( -\sum_{m\ge 1}
m^{-1}\widetilde{\phi}_{0,t}| T_-(m)(Z)\bigr)}
\endmultline
\tag{2.7}
$$
whenever the product converges.  Thus $B_\phi(Z)$ transforms
like a $\Gamma_{\infty,t}$-modular form of weight
$\frac{f(0,0)}2$ with the multiplier system of the theorem.
It is useful to write down the whole product $B_\phi(Z)$
in terms of Hecke operators $T_-(m)$.
We can  get such expression using
 the involution $V_t$ (see $\thetag{1.19}$)
$$\multline
\ml(\phi_{0,t})(Z)=
B_\phi(Z)=\\
=q^{A}r^Bs^{C}
\exp{\bigl(-\sum_{m\ge 1}
m^{-1}\widetilde{\phi}_{0,t}| T_-(m)(Z)\bigr)}
\exp{\bigl(-\sum_{m\ge 1}
m^{-1}(\widetilde{\phi}_{0,t}^{(0)}+{\phi}_{0,t}^{(0)})| T_-(m)|V_t(Z)\bigr)}.
\endmultline
\tag{2.8}
$$
The functions
${\phi}_{0,t}^{(0)}(z)=\sum_{l}f(0,l)r^l$ and
$\widetilde{\phi}_{0,t}^{(0)}(Z)=\sum_{l}f(0,l)r^ls^t$
are not Jacobi forms, and we fix the standard system of representatives
\thetag{1.11} in $T_-(m)$ to define the corresponding formal action.
The exponent of  the function
$\widetilde{\phi}_{0,t}^{(0)}$  in \thetag{2.8}
defines the subproduct over all $(n,l,0)$ with $n>0$ in \thetag{2.3}.
The exponent with the function  ${\phi}_{0,t}^{(0)}$,
which does not depend on $\tau$ and $\omega$,
defines the finite subproduct
over $(0,l,0)$ with $l>0$.
The representation $\thetag{2.8}$ shows us analogy between
the exponential lifting and the lifting of holomorphic forms
defined in Theorem 1.12. In \S 3 we shall use $\thetag{2.8}$
to prove that the exponential lifting commutes
with some Hecke correspondence.
(In fact one can consider the factor before the exponent in \thetag{2.7}
as the Hecke operator $T(0)$.)

Let us check the behavior of $B_\phi(Z)$ with respect to
the involution
$
V_t: (q,r,s)\mapsto (s^t,r, q^{1/t})
$.
In the product only the terms $(1-q^nr^ls^{tm})^{f(nm,l)}$
with $n<0$ do not have
$V_t$-pairing terms.
Thus
$$
\frac{B_\phi(V_t(Z))}{B_\phi(Z)}=
(sq^{1/t})^{{tA}-tD-{C}}
\prod\Sb(n,l,m)>0\\
\vspace{0.5\jot} n<0\endSb
\frac{(s^{-tn}-r^lq^m)^{f(nm,l)}}
{(q^{-n}r^l-s^{tm})^{f(nm,l)}}
=(-1)^D(sq^{1/t})^{{tA}-tD-{C}}
$$
where $D$ was defined in \thetag{2.5}.
The finite product in the last formula is equal to
$(-1)^D$. The first factor is equal to one according
to
\enddemo
\proclaim{Lemma 2.2}For arbitrary nearly holomorphic Jacobi
form  $\phi_{0,t}=\sum_{n,l}f(n,l)q^nr^l$ the identity
$$
24({tA}-tD-{C})=t\sum_{l}f(0,l)-24t\sum\Sb n<0, l\endSb \sigma_1(n)f(n,l)-
6\sum_{l}l^2f(0,l)=0
$$
is valid.
\endproclaim
\demo{Proof of the lemma} We use
the differential operator
$$
L_k=8\pi i \pd{}{\tau}-\dfrac{\partial^2{}}{{\partial z}^2}
-\left(\frac{2k-1}{z}\right)\dfrac{\partial{}}{\partial z}
$$
defined in \cite{EZ, \S 3}.
If the Jacobi form
$\phi_{12,t}(\tau,z)=\Delta(\tau)\phi_{0,t}(\tau,z)$
is weak holomorphic
(i.e. when $f(n,l)\ne 0$ implies $n\ge -1$), then
$t\sum_{l}f(0,l)-24tf(-1,l)-6\sum_{l}l^2f(0,l)$
is the constant term of
$(\Cal D_{2}\phi_{12,t})(\tau)$ (see \cite{EZ, Theorem 3.1})
which is a cusp form of weight $14$
for $SL_2(\bz)$. Thus it is equal to zero.
For arbitrary $\phi_{0,t}$ let us consider
$$
f_2(\tau)=
\Delta(\tau)^{-1}
L_{12}\bigl(\Delta(\tau)\phi_{0,t}(\tau,z)\bigr)|_{z=0}
$$
which
is a nearly holomorphic $SL_2(\bz)$-form of weight $2$.
Let
$
\phi_{0,t}(\tau,z)=\sum_{\nu\ge 0}\chi_{\nu}(\tau)z^n
$
is a Taylor expansion around $z=0$ where
$
\chi_{\nu}(\tau)
=\frac{(2\pi i)^{\nu}}{\nu!}\sum_{n}\bigl(\sum_{l} l^\nu f(n,l)\bigr)q^n
$.
Thus
$$
f_2(\tau)=8\pi i\biggl(\frac {\Delta'(\tau)}{\Delta(\tau)}-
\chi_0'(\tau)\biggr)-48\chi_2(\tau)
$$
where
$
(2\pi i)^{-1}\frac {\Delta'(\tau)}{\Delta(\tau)}=
E_2(\tau)=1-24\sum_{n\ge 1}\sigma_1(n)q^n
$.
It shows that the sum in the right hand side of the  identity
of the lemma is (up to the constant $(4\pi i)^2$) the constant term
of $f_2(\tau)$. But the constant term of
any nearly holomorphic modular form of weight $2$ is equal
to zero (see \cite{Bo6, Lemma 9.2}).
The lemma is proved.
\enddemo

The involution $V_t$ and the group
$\Gamma_{\infty,t}$ generate the double extension $\Gamma_t^+$
of the paramodular group.
Hence we have proved that  $B_\phi(Z)$ transforms like
a Siegel modular form of weight $\frac{f(0,0)}{2}$
with a character (or a multiplier system)
$v:\Gamma_t^+\to \bc^*$ induced by
the character of the product of the Jacobi forms
and the Dedekind $\eta$-functions in  \thetag{2.7}
together with  the relation \thetag{2.5}.

By \cite{Bo6} (see the proof of Theorem 5.1 and Theorem 10.1 there)
any two branches of the analytical continuation of the modular
product $B_\phi$ differ by multiplication on a non-zero constant.
Any Humbert surface in $\Cal A_t^+$ can be
represented in the form
$H_\ell$ with  a primitive $\ell=(0,a,-\frac{b}{2t},1,0)\in L_t^*$ 
(see Sect. 1.3).
For such $\ell$ with $D=2t(\ell,\ell)=b^2-4at>0$,
the divisor  $\Cal H_\ell$  has
non-trivial intersection with a neighborhood  of $\bh_2$
at infinity
where the product
$B_\phi(Z)$ converges if  it does not coincide with zero
of a factor of type  $(1-q^{na}r^{nb}s^{tn})$
with $f(n^2a,nb)\ne 0$
of the product.   Therefore $B_\phi(Z)$
is holomorphic univalent along any  quadratic divisor
$H_D(b)$  or it has zero (or pole)
along this divisor of order
$m_{D,b}=\sum_{n>0}f(n^2a,nb)$.
The theorem is proved.

\subhead
2.2 The basic Siegel modular forms
\endsubhead
In the rest of the paper we show
that many of  modular forms of small weights constructed
as the arithmetic lifting of Jacobi forms of half-integral indices have
infinite product expansion,
i.e. they can be represented as an exponential
lifting. In this subsection we consider the most fundamental examples
of Siegel modular forms $\Delta_{1/2}(Z)$,
$\Delta_{1}(Z)$, $\Delta_{2}(Z)$ and
$\Delta_{5}(Z)$ with the divisor equals to  Humbert
modular surface $H_1$ in the corresponding threefold $\Cal A_t^+$
(for $t=4$, $3$, $2$, $1$ respectively).
These modular forms  define  the automorphic corrections of
the  most important  Lorentzian Kac--Moody algebras of signature $(2,1)$.
An automorphic correction
of a Lorentzian Kac--Moody algebra is defined by the Fourier
coefficients of a modular form with appropriate character
 with respect to the Weyl group of a reflective
hyperbolic lattice. According to the Weyl--Kac--Borcherds denominator
formula for generalized Kac--Moody superalgebras,
the multiplicities of the infinite product expansion define
the multiplicities of positive roots of the corresponding
algebra.

We start with some general remarks.
If one has an automorphic product of type considered above,
then  it is an important and difficult problem to calculate 
an exact form of the Fourier coefficients of this product.
For example, this coefficients defines the set of  imaginary simple
roots of automorphic generalized Kac--Moody algebra
(see \cite{Bo1}--\cite{Bo5} and \cite{GN1}--\cite{GN4}).
Our approach to this problem is to find identities between 
arithmetic and exponential liftings.
If an arithmetic lifting  $\hbox{Lift}(\phi_{k,d})$ would
coincide with exponential one, then
there exists a natural formula for the Jacobi form of weight zero
$\phi_{0,t}$
which is the datum for the exponential lifting
$$
\phi_{0,t}=\frac{m^{2-k}\phi_{k,d}|_k T_-(m_2)}{\phi_{k,d}}
\tag{2.9}
$$
where $T_-(m_2)$ is the second Hecke operator of the sequence
of operators used to produce the arithmetic lifting of $\phi_{k,d}$.
Thus one needs to check
Fourier coefficients with negative norm of
${\phi_{k,d}|_kT_-(m_2)}/{\phi_{k,d}}$.
If it is good enough, then using information about divisors
or arguments with dimension of the space of modular forms one can
prove an identity between liftings. We really do this in many cases
below.
To illustrate this method, we start with an  infinite product expansion
of  modular forms of the singular weight constructed in Theorem 1.11.

\example{Example 2.3}{\it Infinite product expansion of theta-functions.}
Let us analyze the Jacobi  form $\thetag{2.9}$
for the case of two liftings of singular Jacobi forms
(see Theorem 1.11).
For the case of the theta-function $\Delta_{1/2}(Z)$
we get  the weak Jacobi form defined in \thetag{1.33}
$$
\gather
\phi_{0,4}(\tau ,z)=
\frac {\vartheta(\tau ,3z)}{\vartheta(\tau ,z)}=
\sum_{n\ge 0,\, l\in \bz}f_4(n,l)q^nr^l  \\
=r^{-1}\prod_{m\ge 1}(1+q^{m-1}r+q^{2m-2}r^{2})
(1+q^{m}r^{-1}+q^{2m}r^{-2})
\prod\Sb n\equiv 1,2 \,mod\,3\\ n\ge 1\endSb
(1-q^{n}r^3) (1-q^{n}r^{-3})
\\
=(r+1+r^{-1})-q(r^4+r^3-r+2-r^{-1}+r^{-3}+r^{-4})
+q^2(\dots)
\tag{2.10}
\endgather
$$
where all Fourier coefficients $f_4(n,l)$
of the weak Jacobi form are integral
(in fact they are Fourier coefficients
of automorphic forms of weight $-1/2$).
Thus according to Theorem 2.1,
$\hbox{Exp-Lift}(\phi_{0,4})$  is a modular form of weight $1/2$
with respect to the paramodular group $\Gamma_4^+$ having
irreducible Humbert modular surface  $H_1$ as its divisor.
It implies that the quotient
$\Delta_{1/2}(Z)/\hbox{Exp-Lift}(\phi_{0,4})(Z)$
is  a holomorphic automorphic function invariant  with respect to
$\Gamma_4^+$, thus it is a constant.
Moreover we get the following infinite product expansion
of $\Delta_{1/2}(Z)$:
$$
\frac{1}2\sum\Sb n,m\in \bz \endSb
\,\biggl(\dsize\frac{-4}{n}\biggr)\biggl(\dsize\frac{-4}{m}\biggr)
q^{n^2/8}r^{nm/2}s^{m^2/2}=
q^{1/8}r^{1/2}s^{1/2}
\prod_{(n,l,m)> 0}(1-q^nr^ls^{4m})^{f_4(nm,l)}.
\tag{2.11}
$$

Let us consider the weak Jacobi form \thetag{2.9}
connected with $D_{1/2}(Z)$
$$
\multline
\phi_{0,36}(\tau,z)=\frac{\vartheta_{3/2}(\tau,5z)}
{\vartheta_{3/2}(\tau,z)}=
\frac{\vartheta(\tau,10z)\vartheta(\tau,z)}
{\vartheta(\tau,5z)\vartheta(\tau,2z)}=
\sum_{n\ge 0,\ l\in \bz}
f_{36}(n,l)\,q^nr^l\\
=r^{-2}\prod_{n\ge 1}
\frac{(1+q^{n-1}r^5)(1+q^nr^{-5})
(1-q^{2n-1}r^{10})(1-q^{2n-1}r^{-10})}
{(1+q^{n-1}r)(1+q^nr^{-1})
(1-q^{2n-1}r^{2})(1-q^{2n-1}r^{-2})}\\
=(r^2-r^1+1-r^{-1}+r^{-2})+q^2(-r^{17}+\dots)
+q^5(r^{27}+\dots)\\
+q^7(r^{32}+\dots)+q^8(r^{34}+\dots)+\dots,
\endmultline
\tag{2.12}
$$
where we include in the last formula only summands $q^nr^l$ with 
the negative norm $144n-l^2$:
$144\cdot 2-{17}^2=-1$, $144\cdot 5-{27}^2=-9$,
$144\cdot 7-{32}^2=-16$,
$144\cdot 8-{34}^2=-4$.

In calculations with weak Jacobi forms we shall often use
the following simple considerations. Let
$$
\phi_{0,t}(\tau,z)=\sum_{n\ge 0, \, l\in \bz}f_t(n,l)q^nr^l\in
J_{0,t}^{weak}.
$$
Then
$f_t(n,l)$ depends only on the norm  $4nt-l^2$ and
$\pm l\operatorname{mod}2t$;
$f_t(n,l)=0$  if $4nt-l^2<-t^2$. Moreover,
to calculate all Fourier coefficients with negative norm,
it is enough to find
$f_t(n,l)$ for
$$
n\le \bigl\lfloor \frac{t}4\bigl\rfloor \
(t\not\equiv 0\,\operatorname{mod}\, 4)\quad
\text{ or } \quad  n \le \frac{t}4 -1
\ (t\equiv 0\,\operatorname{mod}\, 4).
\tag{2.13}
$$
These arguments imply the fact that  the last  part of   \thetag{2.12} 
contains all orbits of  possible Fourier coefficients with negative norm.
(See \thetag{4.5} below for  another formula for $\phi_{0,36}(\tau,z)$.)
Using Theorem 2.1,  we get
$$
\hbox{Div}_{\Cal A_{36}^+}\bigl(\hbox{Exp-Lift}(\phi_{0,36})\bigr)=
H_4(2)+H_4(34)+H_9(27)+H_{16}(32).
$$
(See the notation \thetag{1.23}). The surfaces
$H_4(2)$ and $H_4(34)$ are equivalent with respect
to the group $\Gamma_{36}^*$ (see Theorem 1.11).
This finishes the proof of the statement of  Theorem 1.11
about the divisor of $D_{1/2}(Z)$ and gives us the
infinite product expansion of $D_{1/2}(Z)$:
$$
\frac{1}2\sum_{m,n\in \bz}
\biggl(\frac{12}n\biggl)\biggl(\frac{12}m\biggl)
\,q^{{n^2}/{24}}r^{{nm}/2}s^{{3m^2}/{2}}=
q^{1/24}r^{1/2}s^{3/2}
\prod_{(n,l,m)> 0}(1-q^nr^ls^{36m})^{f_{36}(nm,l)}.
\tag{2.14}
$$
\endexample

\example{Example 2.4}{\it The modular form $\Delta_5(Z)$.}
We recall that there exists a weak Jacobi form of weight zero
and  index one
$$
\phi_{0,1}(\tau,z)=\frac{\phi_{12,1}(\tau,z)}{\eta(\tau)^{24}}
=(r+10+r^{-1})+q(10r^{-2}-64r^{-1}+108-64r+10r^2)+q^2(\dots )
\tag{2.15}
$$
which is one of  the standard generators of the  graded
ring of weak Jacobi forms.
There is a formula for  Fourier coefficients
of $\phi_{12,1}(\tau,z)$ in terms of H. Cohen's numbers
(see \cite{EZ}, \S 9).
The  function $\phi_{0,1}$ gives
the following  result from \cite{GN1}
about the product of even theta-constants:
$$
\align
\Delta_5(Z)
&=\sum\Sb n,l,m \equiv 1\operatorname{mod}2\\
\vspace{0.5\jot} n,m>0\endSb \hskip10pt
\sum_{a|(n,l,m)} (-1)^{\frac{l+a+2}2}\,a^{4}
\tau_9\,(\frac{4nm-l^2}{a^2})\,
q^{n/2}\,r^{l/2}\,s^{m/2}\\
\vspace{2\jot}
{}&=(qrs)^{1/2} \prod
\Sb n,\,l,\,m\in \Bbb Z\\
\vspace{0.5\jot}
(n,l,m)>0\,\endSb
\bigl(1-q^n r^l s^m\bigr)^{f_{1}(nm,l)}
\in \frak N_5(\Gamma_1, \chi_2)
\tag{2.16}
\endalign
$$
where $f_1(n,l)$ are the Fourier coefficients of $\phi_{0,1}(\tau,z)$.
We note  that \thetag{2.9} provides us with another formula
for $\phi_{0,1}(\tau,z)$:
$$
\phi_{0,1}(\tau,z)=-
\frac{(\eta(\tau)^9\vartheta(\tau,z))|_5 T_-(3)}
{\eta(\tau)^9\vartheta(\tau,z)}
\tag{2.17}
$$
which represents $\phi_{0,1}(\tau,z)$ as a sum of four
infinite products.
\endexample

Our next purpose is to  construct the infinite product
expansion of the modular forms $\Delta_1$ and $\Delta_2$
defined in Example 1.14.
For a large $m$ ($m=5$ and $7$ in the case under consideration)
the calculation of
$\phi_{k,d}|T_-(m)/{\phi_{k,d}}$ from $\thetag{2.9}$
might be rather long, thus we prefer to construct the same function
in a  non-direct way.
In the next lemma we  define two weak Jacobi forms of weight $0$
using the quintiple product and the Jacobi forms constructed
in Lemma 1.24.

\proclaim{Lemma 2.5}The functions
$$
\gather
\phi_{0,2}(\tau ,z)=\frac{\phi_{2,2}(\tau ,z)}{\eta(\tau )^4}
={\tsize\frac{1}2} \eta(\tau )^{-4}
\sum_{m\,,n\in \bz}
{(3m-n)}\biggl(\frac{-4}m\biggl)\biggl(\frac{12}n\biggl)
q^{\frac{3m^2+n^2}{24}}r^{\frac{m+n}2}
\\
=\sum_{n\ge 0,\,l\in \bz}f_2(n,l)q^nr^l=
(r+4+r^{-1})+q(r^{\pm 3}-8r^{\pm 2}-r^{\pm 1}
+16)
+q^2(\dots)
\tag{2.18}
\\
\intertext{($r^{\pm l}$ means that we have two summands with
$r^{l}$ and with $r^{-l}$) and}
\phi_{0,3}(\tau ,z)=
\biggl(\frac {\vartheta_{3/2}(\tau ,z)}{\eta(\tau )}\biggr)^2=
\biggl(\frac {\vartheta(\tau ,2z)}{\vartheta(\tau ,z)}\biggr)^2
=\sum_{n\ge 0,\,l}f_3(n,l)q^nr^l\\
\vspace{2\jot}
=r^{-1}
\biggl(\prod_{n\ge 1}(1+q^{n-1}r)(1+q^{n}r^{-1})(1-q^{2n-1}r^2)
(1-q^{2n-1}r^{-2})\biggr)^2\\
{}=(r+2+r^{-1})+q(-4r^{\pm 3}-4r^{\pm 2}+2r^{\pm 1}+4)
+q^2(\dots)
\tag{2.19}
\endgather
$$
are weak Jacobi forms of weight $0$
and index $2$ and $3$ respectively. Moreover, their Fourier coefficients
satisfy the property: if\  $4tn-l^2<0$ and $f_t(n,l)\ne 0$,
then $4tn-l^2=-1$ and $f_t(n,l)=1$.
\endproclaim
\demo{Proof}According to \thetag{2.13} one has to calculate the $q^0$-part
of the Fourier expansions of the given Jacobi forms.
We also remark that for a prime $p$  Fourier coefficient
$f_p(n,l)$ of a Jacobi form of index $p$ depends only on the norm
$4np-l^2$. Thus \thetag{2.15} gives  the Fourier coefficients
$f_2(n,l)$ with $8n-l^2=-1$, $4$, $7$, $8$ and
\thetag{2.17} gives
$f_3(n,l)$ with $12n-l^2=-1$, $3$, $8$, $11$, $12$.
\newline
\qed
\enddemo

The modular form $\Delta_1$ with respect to $\Gamma_3$ and
$\Delta_2$ with respect to  $\Gamma_2$ (see Example  1.14) are
analogous to $\Delta_5$ and they also satisfy
a generalized Euler-type  identity
(for $\Delta_2$ it was given in \cite{GN1}).
Using the modular forms $\Delta_1$, $\Delta_2$ and
$\Delta_5$,  we can define  the  automorphic corrections of
symmetric generalized Cartan matrices of elliptic type of rank three
with a lattice Weyl vector of Theorem 1.3.1 in Part I.
Comparing the divisors of $\Delta_2(Z)$ and $\Delta_1(Z)$
with the divisors of $\hbox{Exp-Lift}(\phi_{0,t})(Z)$
for $t=2$ and $3$ respectively
we obtain

\proclaim{Theorem 2.6}The following identities are valid:
$$
\align
\Delta_1(Z)&=\sum_{M\ge 1}
\sum
\Sb
m >0,\,l\in \bz\\
\vspace{0.5\jot} n,\,m\equiv 1\,mod\,6\\
\vspace{0.5\jot}
4nm-3l^2=M^2
\endSb
\hskip-4pt
\biggl(\dsize\frac{-4}{l}\biggr)
\biggl(\dsize\frac{12}{M}\biggr)
\sum\Sb a|(n,l,m)\endSb \biggl(\dsize\frac{6}{a}\biggr)
q^{n/6}r^{l/2}s^{m/2}\\
{}&=
q^{\frac {1}6}r^{\frac {1}2}s^{\frac {1}2}
\prod
\Sb n,\,l,\,m\in \Bbb Z\\
\vspace{0.5\jot}
(n,l,m)>0
\endSb
\bigl(1-q^n r^l s^{3m}\bigr)^{f_{3}(nm,l)}
\in
\frak N_{1}(\Gamma_3^+, \chi_6)
\tag{2.20}
\endalign
$$
where the character $\chi_6: \Gamma_3^+\to \root{6}\of{1}$
is induced
by $v_\eta^4\times v_H$ and
$$
\align
\Delta_2(Z)&=\sum_{N\ge 1}\
\sum\Sb m>0,\,l\in \bz\\
\vspace{0.5\jot}
 n,\, m\equiv 1\, mod\, 4\\
\vspace{0.5\jot}
2mn-l^2=N^2
\endSb
N\biggl(\frac {-4}{Nl}\biggr)
\sum_{a\,|\,(n,l,m)} \biggl(\frac {-4}{a}\biggr)
\, q^{n/4} r^{l/2} s^{m/2}\\
{}&=
q^{\frac {1}4}r^{\frac {1}2}s^{\frac {1}2} \prod
\Sb n,\,l,\,m\in \Bbb Z\\
\vspace{0.5\jot}
(n,l,m)>0\,\endSb
\bigl(1-q^n r^l s^{2m}\bigr)^{f_{2}(nm,l)}
\in
\frak N_{2}(\Gamma_2^+, \chi_4)
\tag{2.21}
\endalign
$$
where $\chi_4: \Gamma_2^+\to \root{4}\of{1}$ is induced
by $v_\eta^6\times v_H$.
Moreover the divisor of these modular forms is
the irreducible Humbert surface $H_1$
$$
\hbox{Div}_{\Cal A_3^+}(\Delta_1(Z))=H_1,\qquad
\hbox{Div}_{\Cal A_2^+}(\Delta_2(Z))=H_1.
$$
\endproclaim

One can check (one only  needs to calculate the $q^0$-part
of the corresponding Fourier expansion) that the weak Jacobi forms
used  above as the data for the exponential lifting
of all $\Delta_k$-forms
are related by   simple  relations
$$
\phi_{0,1}(\tau,2z)=\phi_{0,2}^2(\tau,z)-8\phi_{0,4}(\tau,z)=
\phi_{0,1}(\tau,z)\cdot \phi_{0,3}(\tau,z)-12\phi_{0,4}(\tau,z).
\tag{2.22}
$$
The form $\phi_{0,3}$ is the square of a weak Jacobi form of index
$\frac {3}2$ with character $\hbox{id}_{SL_2}\times v_H$
$$
\xi_{0,\frac{3}2}(\tau,z)
=\frac{\vartheta(\tau,2z)}{\vartheta(\tau,z)}=
r^{-\frac{1}2}
\prod_{n\ge 1}(1+q^{n-1}r^{1})(1+q^{n}r^{-1})
(1-q^{2n-1}r^2)(1-q^{2n-1}r^{-2})
$$
(or, equivalently,
$\xi_{0,\frac{3}2}(\tau,z)
=\eta(\tau)^{-1}\vartheta_{3/2}(\tau,z))
$.
We set
$$
\xi_{0,6}(\tau,z)=\xi_{0,\frac{3}2}(\tau,2z)=
(r+r^{-1})-q(r^5+r+r^{-1}+r^{-5})+q^2(\dots)\in J_{0,6}^{weak}.
\tag{2.23}
$$
The weak Jacobi form  $\xi_{0,6}$ has only
two orbits of the Fourier coefficients with negative norm of their indices.
They are  $r$ and $-qr^5$.
It is clear that $\eta(\tau)\xi_{0,6}(\tau,z)$ is holomorphic.
The weak Jacobi  $\phi_{0,2}$ also has an expression in terms of
Jacobi theta-series
$$
\phi_{0,2}\cdot \phi_{0,4}-
\phi_{0,3}^2=\xi_{0,6}
\quad\text{or}\quad
\phi_{0,2}(\tau,z)=\frac{\vartheta(\tau,2z)^2}
{\vartheta(\tau,z)\vartheta(\tau,3z)}+
\frac{\vartheta(\tau,4z)\vartheta(\tau,z)}
{\vartheta(\tau,2z)\vartheta(\tau,3z)}.
\tag{2.24}
$$
We obtain similar formulae  for $\phi_{0,1}(\tau,z)$ below
(see \thetag{3.16},
\thetag{3.34} and \thetag{3.35}).
\head
\S 3. Liftings and  Hecke correspondence
\endhead

In this section we shall prove that the exponential lifting
commutes with some Hecke correspondence.
We consider the  multiplicative Hecke operators used in
\cite{GN4} and the multiplicative analogue of the
operator of symmetrisation studied in \cite{G2}.

\subhead
3.1.  Multiplicative symmetrisation
\endsubhead

In \cite{G2} we defined a Hecke operator  which transforms modular forms
with respect  to $\Gamma_t$ into modular forms with respect to
$\Gamma_{tp}$
$$
\hbox{Sym}_{t,p}: \frak M_{k}(\Gamma_t)\to \frak M_{k}(\Gamma_{tp}),
\qquad
\hbox{Sym}_{t,p}: F\mapsto \sum_{M\in (\Gamma_t\cap \Gamma_{tp})
\setminus \Gamma_{tp}} F|_k M.
\tag{3.1}
$$
We call this operator {\it the operator of $p$-symmetrisation}.
It can be represented as an action  of an  element from the Hecke ring
$H(\Gamma_{\infty,t})$ of the parabolic subgroup
$\Gamma_{\infty,t}\subset \Gamma_t$.
To see this, we take a system of representatives
$$
(\Gamma_t\cap \Gamma_{tp})
\setminus \Gamma_{tp}=
\bigl\{J_{tp},\ \ \nabla({\frac {b}{tp}}),\   b\,\hbox{mod}\, p
\bigr\}
\tag{3.2}
$$
where
$$
J_t=
\Biggl(\smallmatrix
o&0&1&0\\
0&0&0&t^{-1}\\
-1&0&1&0\\
0&-t&0&0
\endsmallmatrix\Biggr), \qquad
\nabla(a)=
\Biggl(\smallmatrix
1&0&0&0\\
0&1&0&a\\
0&0&1&0\\
0&0&0&1
\endsmallmatrix\Biggr).
$$
It is valid
$J_tJ_{tp}=\hbox{diag}(1,p,1,p^{-1})$, thus
for any $F\in \frak M_{k}(\Gamma_t)$ we have
$$
\hbox{Sym}_{t,p}(F)=F|_k \bigl(
\Lambda_p+\sum_{b\, mod\, p}\nabla(\frac {b}{tp})\bigr).
$$
The last operator is defined by the
following  element in the parabolic Hecke ring
$$
\hbox{Sym}_p=\Lambda_p+\nabla_{t,p}\in H(\Gamma_{\infty,t}),
\quad\text{where}\quad
\nabla_{t,p}=\sum_{b\, mod\, p}\Gamma_{\infty,t}\nabla(\frac {b}{tp})
\tag{3.3}
$$
and $\Lambda_p$ is the element \thetag{1.10}.
The operator of $p$-symmetrisation  is injective
if $(t,p)=1$ and it commutes with the arithmetic
lifting. If $\phi\in J_{k,t}$, then
$$
\hbox{Sym}_{t,p}(\hbox{Lift}(\phi))=p^{3-k}\,
\hbox{Lift}(\phi|_k T_{-}(p)).
\tag{3.4}
$$
(see  Satz 2.10 and Satz 3.1 in \cite{G2}).
Let us define the multiplicative analogue of the $p$-symmetrisation.

\definition{Definition 3.1}Let $F\in \frak M_k(\Gamma_t, \chi)$
where $k\in \bz/2$.
Then for a prime $p$ we define  {\it the  operator of multiplicative
symmetrisation}
$$
\ms_p:
F\mapsto p^{-k}\,\overline{\chi}(J_t)\prod_{M_i\in \Gamma_t\cap \Gamma_{tp}
\setminus \Gamma_{tp}}
F|_kM_i.
\tag{3.5}
$$
(The additional  constant $ p^{-k}\overline{\chi}(J_t)$ makes
formulae simpler.)
\enddefinition

\proclaim{Lemma 3.2}Let $p$ be an arbitrary prime
and $F\in \frak M_k(\Gamma_t, \chi)$.
Then
$$
\ms_p(F)(Z)\in \frak M_{k(p+1)}(\Gamma_{tp}, \chi^{(p)})
$$
where $\chi^{(p)}$ is a character of $\Gamma_{tp}$.
Moreover if  the modular form $F$ is zero along a
Humbert surface $H_l\subset \Cal A_t$
of  discriminant $D$, then $\ms_p(F)$ is zero along
$\ms_p^*(H_l)$ which is a sum (with some multiplicities)
of Humbert surfaces  with  discriminant
$D$ and $p^2D$ in $\Cal A_{tp}$.
\endproclaim
\demo{Proof}Using \thetag{3.2} and interpretation \thetag{3.3}
of the  $p$-symmetrisation as a Hecke operator,
we get
$$
\ms_p(F)(Z)=
F(\pmatrix \tau&pz\\pz&p^2\omega\endpmatrix) \prod_{b\,mod\, p}
F(\pmatrix \tau&z\\z&\omega+\frac{b}{tp}\endpmatrix)
\in M_{k(p+1)}(\Gamma_{tp}, \chi^{(p)})
\tag{3.6}
$$
where $\chi^{(p)}$ is the character with the properties
$\chi^{(p)}|_{SL_2(\bz)}=(\chi|_{SL_2(\bz)})^{p+1}$
and $\chi^{(p)}|_{H(\bz)}=(\chi|_{H(\bz)})^{p}$.
The formula \thetag{3.6} explains also
the action of $\ms_p$  on the rational
quadratic divisors.
\newline
\qed
\enddemo
Similar to the arithmetic lifting and  the $p$-symmetrisation, 
the exponential lifting commutes
with the multiplicative symmetrisation. We have the following
analog of \thetag{3.4}.
\proclaim{Theorem 3.3}Let  $\phi\in J_{0,t}^{nh}$ be like in Theorem 2.1.
Then for an arbitrary prime p we have
$$
\ms_p(\hbox{Exp-Lift}(\phi_{0,t}))=
\hbox{Exp-Lift}(\phi_{0,t}| T_-(p)).
$$
\endproclaim
\demo{Proof}To prove the theorem, we use the representation
\thetag{2.8} of the exponential lifting.
The formula \thetag{3.6} shows that $\ms_p$ can be written
as the action of $\hbox{Sym}_{p}$ on the function under the exponent
in \thetag{2.8}.
Let us consider the product of  the formal Dirichlet series
$\sum_{m=1}^{\infty}\,T_{-}(m)m^{-1}$
over the Hecke ring $H(\Gamma_{\infty,t})$ with
$\hbox{Sym}_{p}=\Lambda_p+\nabla_{t,p}$.
According to our definition of the normalizing factor
of  Hecke operators in the case of weight zero (see Sect. 1.2),
we can consider the Hecke ring $H(\Gamma_{\infty,t})$ modulo
its central element $\Delta(p)=\Gamma_{\infty,t}(pE_4)$.
I.e. for any  $X\in H(\Gamma_{\infty,t})$
the Hecke operators $X$ and $\Delta(p)X$ are identical.
We recall that $\Delta(p)\Lambda_p=T_{-}(p,p)$
where  $T_{-}(p,p)$ is the embedding of
$T(p,p)=SL_2(\bz)(pE_2) SL_2(\bz)\in H(SL_2(\bz))$
into the Hecke ring $H(\Gamma_{\infty,t})$.
Using the definition one can check that
$$
T_-(m)\nabla_{t,p}
=\cases
pT_-(m)&\text{if } m\equiv 0\ \hbox{mod } p\\
\nabla_{t,p}\, T_-(m)&\text{if } m\not\equiv 0\hbox{ mod } p.
\endcases
$$
Thus
$$
\multline
\biggl(\sum_{m=1}^{\infty}\,T_{-}(m)\,m^{-1}\biggr)
(T_-(p,p)+\nabla_{t,p})\\
=
\sum_{m\ge 1}\biggl(T_-(mp)+p T_-\bigl(\frac m{p}\bigr)
T_-(p,p)\biggr)\,m^{-1}
+\sum_{(m,p)=1}\nabla_{t,p} T_-(m)\,m^{-1}\\
=
T_-(p)\sum_{m\ge 1} T_-(m)\,m^{-1}
+\nabla_{t,p} \cdot\sum_{(m,p)=1} T_-(m)\,m^{-1},
\endmultline
\tag{3.7}
$$
since
$T(m)T(p)=T(mp)+pT(p,p)T(\frac m{p})$  in $H(SL_2(\bz))$.

Let us consider the representation \thetag{2.8} for
the exponential lifting of $\phi_{0,t}$.
According to \thetag{3.7} we have the following identity for
the main factor in \thetag{2.8}
$$
\exp{\biggl(-\sum_{m\ge 1}
\bigl(\widetilde{\phi}_{0,t}| T_-(m)\bigr)|\hbox{Sym}_p (Z)\biggr)}=
\exp{\biggl(-\sum_{m\ge 1}
\bigl(\widetilde{\phi}_{0,t}| T_-(p)\bigr)| T_-(m) (Z)\biggr)},
$$
since $\nabla_{t,p}$ defines zero operator on the space of
Jacobi functions of weight $0$ and index $t$:
$$
(\widetilde{\phi}_{0,t}|\nabla_{t,p})(Z)=
\widetilde{\phi}_{0,t}(Z)\cdot\sum_{b\,mod \,p}
\exp{(\frac{2\pi i  b}{p})}=0.
$$
Let us consider the second exponent in \thetag{2.8}.
We have the  formal identity
$$
V_t\biggl(\Lambda_p+\sum_{b\,mod \,p}\nabla(\frac{b}{tp})\biggr)
=T_-(p)V_{tp}\,(\sqrt{p}E_4)^{-1}
$$
where we consider $T_-(p)$ as the formal sum of
the left cosets fixed in \thetag{1.11} and $V_t$ is
involution \thetag{1.19}.
The second exponent in \thetag{2.8} is not
$\Gamma_{\infty,t}$-invariant, but it is invariant
with respect to the minimal parabolic subgroup
$\Gamma_{00}$ which  is the intersection
of $\Gamma_{\infty,t}$ with the subgroup of the upper-triangular
matrices in  $\Gamma_t$.
This parabolic subgroup is the semidirect product of
the subgroup of all upper-triangular matrices in $SL_2(\bz)$
with the Heisenberg group. Thus we still can consider $T_-(m)$
as an element of the Hecke ring $H(\Gamma_{00})$
of this minimal parabolic subgroup if we take $T_-(m)$
in the standard form \thetag{1.11}.
(See \cite{G10}, where Hecke rings of  parabolic subgroups
of this type were considered for $GL_n$ over a local field.)
Thus for the second exponent in \thetag{2.8} we get
$$
\multline
\exp{\biggl(-\sum_{m\ge 1}
m^{-1}(\widetilde{\phi}_{0,t}^{0}+{\phi}_{0,t}^{0})| T_-(m)|V_t
|\hbox{Sym}_p\,(Z)\biggr)}\\
=
\exp{\biggl(-\sum_{m\ge 1}
m^{-1}\bigl(\widetilde{\phi}_{0,t}^{0}+{\phi}_{0,t}^{0}\bigr)|
T_-(p)| T_-(m)|V_{tp}\,(Z)\biggr)}.
\endmultline
$$
This finishes the proof.
\newline\qed
\enddemo

The multiplicative symmetrisation
will give us some modular forms of small weight
with simple divisor. We consider below the multiplicative
$p$-symmetrisation
of $\Delta_k$-functions (for $k=1/2$, $1$, $2$, $5$)
for different primes to produce
an infinite product expansion of some functions constructed
in \S 1 as the arithmetic lifting, i.e. we get
new identities between arithmetic and exponential liftings.
We recall that for an arbitrary Jacobi form
$\phi_{0,t}(\tau,z)=\sum_{n,l}f(n,l)q^nr^l$ we have the following
formula for the Fourier coefficient of of the Jacobi form
$\bigl(\phi_{0,t}|T_-(m)\bigr)(\tau,z)=\sum_{n,l}f_m(n,l)q^nr^l$
$$
f_m(N,L)=m\sum_{a|(N,L,m)} a^{-1}f(\frac{Nm}{a^2},\frac{L}a).
\tag{3.8}
$$
In particular,
$f_m(0,0)=\sigma_1(m)f(0,0)$.
We  use the identities (see \thetag{2.10}, \thetag{2.15}
\thetag{2.18}, \thetag{2.19})
$$
\align
\phi_{0,1}|T_-(2)(\tau,z)&=r^{2}+2r^{1}+30+2r^{-1}+r^{-2}+q(\dots),\\
\phi_{0,1}|T_-(3)(\tau,z)&=r^{3}+3r^{1}+40+3r^{-1}+r^{-3}+q(\dots),\\
\phi_{0,2}|T_-(2)(\tau,z)&=r^{2}+2r^{1}+12+2r^{-1}+r^{-2}+q(\dots),\\
\phi_{0,2}|T_-(3)(\tau,z)&=r^{3}+3r^{1}+16+3r^{-1}+r^{-3}+q(\dots),
\tag{3.9}\\
\phi_{0,3}|T_-(2)(\tau,z)&=r^{2}+2r^{1}
+\hphantom{3}6+2r^{-1}+r^{-2}+q(\dots),\\
\phi_{0,3}|T_-(3)(\tau,z)&=r^{3}+3r^{1}+
\hphantom{3}8+3r^{-1}+r^{-3}+q(\dots),\\
\phi_{0,4}|T_-(2)(\tau,z)&=r^{2}+2r^{1}+
\hphantom{3}3+2r^{-1}+r^{-2}+q(\dots).
\endalign
$$
\example{Example 3.4}{\it The modular form $\Delta_{11}$}.
Let us consider $\ms_2(\Delta_5)$.
According to the first identity in \thetag{3.9},  we get
$$
\Delta_{11}(Z)=\frac{\ms_2(\Delta_5(Z))}{\Delta_2(Z)^2}\in
\frak N_{11}(\Gamma_2)
\tag{3.10}
$$
with trivial character and with
$\hbox{Div}_{\Cal A_2^{+}}(\Delta_{11})=H_1+H_4$.
The $\gm_2$-modular form of weight $11$
$\hbox{Lift}(\eta^{21}(\tau)\vartheta(\tau,2z))$
defined in Example  1.15 has zero along $H_1+H_4$
(see Lemma 1.16).
Thus  we have proved
the identity
$$
\gather
\Delta_{11}(Z)=
\hbox{Lift}(\eta^{21}(\tau)\vartheta(\tau,2z))=
\hbox{Exp-Lift}\bigl(\phi_{0,1}| T_-(2)-2\phi_{0,2}\bigr)\\
=
qrs^{2}
\prod
\Sb n,\,l,\,m\in \Bbb Z\\
\vspace{0.5\jot}
(n,l,m)>0
\endSb
\bigl(1-q^n r^l s^{2m}\bigr)^{f_{2}^{(11)}(nm,l)}
\in
\frak N_{11}(\Gamma_2)
\tag{3.11}
\endgather
$$
where
$$
\phi_{0,2}^{(11)}(\tau,z)=\phi_{0,1}| T_-(2)(\tau,z)
-2\phi_{0,2}(\tau,z)=
\sum_{n\ge 0,\, l\in \bz}
f_{2}^{(11)}(n,l)\,q^nr^l.
$$
Below we give another  formula for the form
$\Delta_{11}(Z)$ and for the function $\phi_{0,2}^{(11)}$
using a different multiplicative Hecke correspondence
(see \thetag{3.14} and \thetag{3.33}).
\endexample

Theorem 3.2 is valid for all prime $p$, including the divisors of $t$.
Let us consider $p=t=2$. The Hecke operator
$\Lambda_2$ (see \thetag{1.10})
transforms the $\Gamma_1$-modular form $\Delta_5$
to the $\Gamma_4$-modular form $\Delta_5^{(4)}$
$$
\Delta_5^{(4)}(\pmatrix \tau&z\\z&\omega\endpmatrix)
=\Delta_5(\pmatrix \tau&2z\\2z&4\omega\endpmatrix)
\in \frak N_{5}(\Gamma_4, v_\eta^{12}\times \hbox{id}_H)
\tag{3.12}
$$
with divisor $H_4+H_1$ in $\Cal A_4^+$.
Almost the same function we get taking
the  multiplicative $2$-symmetrisation of $\Delta_2(Z)$.
Comparing  divisors of modular forms,  we get the
relation
$$
\ms_2(\Delta_2(Z))=
\Delta_{5}^{(4)}(Z){\Delta_{1/2}(Z)^2}
\in
\frak N_{6}(\Gamma_4, v_\eta^{18}\times \hbox{id}_H).
\tag{3.13}
$$
(The modular form $\ms_2(\Delta_2(Z))$ has the divisor $H_4+3H_1$
in $\Cal A_4^+$.)
The last identity shows that we can define
$\Delta_2(Z)$ in terms of Siegel theta-constants.

In the next lemma we prove some identities between the main
Jacobi forms $\phi_{0,t}$ for $t=1$, $2$, $3$, $4$ which  have been used
in the exponential liftings
\thetag{2.16}, \thetag{2.20}, \thetag{2.21} and
the Jacobi form defined in  \thetag{3.11}.

\proclaim{Lemma 3.5}The identities
$$
\gather
\phi_{0,2}^{(11)}=\phi_{0,1}| T_-(2)-2\phi_{0,2}
=\phi_{0,1}^2-20\phi_{0,2},
\tag{3.14}\\
\phi_{0,2}| T_-(2)(\tau,z)=\phi_{0,1}(\tau,2z)+2\phi_{0,4}(\tau,z),
\tag{3.15}\\
4\phi_{0,1}=\phi_{0,2}| T_+(2)
\tag{3.16}
\endgather
$$
are valid
where $ T_+(2)=\gi\hbox{diag}(1,1,2,2)\gi\in H(\gi)$.
\endproclaim
\demo{Proof}To prove \thetag{3.14} and \thetag{3.15},
we have to calculate the $q^0$-part of the corresponding
Jacobi forms using \thetag{3.8}.
To prove \thetag{3.16}, we apply the Hecke operator
$$
\Lambda_2^*=\gi\hbox{diag}(1,p^{-1},1,p)\gi=
\sum\Sb \lambda,\mu\,mod \,p\\
\kappa\,mod\, p^2\endSb
\gi\pmatrix
1&0&0&\mu\\
p^{-1}\lambda&p^{-1}&p^{-1}\mu&p^{-1}\kappa\\
0&0&1&-\lambda\\
0&0&0&p
\endpmatrix
\in H(\gi)
\tag{3.17}
$$
to \thetag{3.15}:
$$
\phi_{0,2}| (T_-(2)\cdot\Lambda_2^*)(\tau,z)=
\phi_{0,1}|(\Lambda_2\cdot \Lambda_2^*)(\tau,z)
+2\phi_{0,4}|\Lambda_2^*(\tau,z).
$$
We know that
$T_-(p)\Lambda_2^{+}=p^2T_+(p)$ and $\Lambda_p\Lambda_p^{+}=p^4$
(see \cite{G2}, \cite{G3} about properties of these operators
and for more details).
To finish the proof of \thetag{3.16}, we show that
$\phi_{0,4}|\Lambda_2^*=0$ (it means that
$\phi_{0,4}$ is a ``new" Jacobi form). We have
$$
\phi_{0,4}|\Lambda_2^*(\tau,z)=
4\sum_{\mu,\lambda\,mod\,2}
\phi_{0,4}(\tau, {\tsize\frac{\lambda\tau+z+\mu}2})=
\sum_{\lambda\,mod\,2}\sum_{n,l\in \bz} f_4(n,2l)\,
q^{n+\lambda l+\lambda^2}r^{2\lambda+l}.
$$
This Fourier expansion shows that the Jacobi form 
$\phi_{0,4}|\Lambda_2^*$ of weight $0$ and index $1$
is holomorphic, thus
$\phi_{0,4}|\Lambda_2^*\equiv 0$.
\newline
\qed
\enddemo
\remark{Remark 3.6}
Fourier coefficients of the function  $\phi_{0,1}(\tau,z)$,
which  is one of
the canonical generators  of the  graded ring of weak Jacobi
forms, are equal to the   multiplicities of
roots of the automorphic Kac--Moody Lie algebra defined by
$\Delta_5(Z)$.
The identity \thetag{3.16}  gives us a new expression for
this function which is simpler than the expression
$$
\phi_{0,1}(\tau,z)=(144\Delta(\tau))^{-1}
\bigl(E_4^2(\tau)E_{4,1}(\tau,z)-E_6(\tau)E_{6,1}(\tau,z)\bigr)
$$
from  \cite{EZ}.
To calculate the Fourier coefficients of  $\phi_{0,2}$,
one has to use only $\eta(\tau)^{-4}$
and the theta-function $\phi_{2,2}$ (see Lemma 1.24).
Let us calculate the action of $T_+(2)$ on the Fourier
coefficients of $\phi_{0,2}(\tau,z)$. By definition
$$
T_+(2)=\sum_{a,b,c\,mod \,2}
\gi\Biggl(\smallmatrix
1&0&a&b\\
0&1&b&c\\
0&0&2&0\\
0&0&0&2\endsmallmatrix\Biggr)+
\sum_{\lambda,\kappa\,mod \,2}
\gi\Biggl(\smallmatrix
2&0&0&0\\
\lambda&1&0&\kappa\\
0&0&1&-\lambda\\
0&0&0&2
\endsmallmatrix\Biggr).
\tag{3.18}
$$
Thus
$$
\multline
\phi_{0,2}| T_+(2)(\tau,z)\\=2
\sum_{a,b\,mod \,2}\phi_{0,2}(\tau+\frac a{2},z+\frac b{2})+
2\sum_{\lambda\,mod \,2}\phi_{0,2}(2\tau,\lambda \tau+z)
\exp{(2\pi i (\lambda^2\tau+2\lambda z))}\\
=8\sum_{n,l}f_2(2n,2l)\,q^nr^l+
2\sum_{\lambda\,mod \,2}
f_2({\tsize\frac {n+\lambda l+\lambda^2}2},\, l+2\lambda)\,q^nr^l.
\endmultline
$$
It gives us the following formula for  the Fourier coefficients
$f_1(N)=f_1(n,l)$ ($N=4n-l^2$) of the Jacobi form $\phi_{0,1}(\tau,z)$
in terms of
Fourier coefficients $f_2(M)=f_2(n,l)$ ($M=8n-l^2$) of
$\phi_{0,2}(\tau,z)$
$$
f_1(N)=2f_2(4N)+
2^{-1}\bigl(\biggl(\frac {-N}{2}\biggr)+1\bigr)f_2(N)
\tag{3.19}
$$
where
$$
\biggl(\frac {-N}{2}\biggr)=\biggl(\frac {-N}{8}\biggr)=
\cases
\hphantom{-}1 &\text{if } -N\equiv 1\hbox{ mod } 8\\
-1 &\text{if } -N\equiv 5\hbox{ mod } 8\\
\hphantom{-}0 &\text{if } -N\equiv 0\hbox{ mod } 4.
\endcases
$$
\endremark

We continue the list of  Euler type  identities between infinite
sums and infinite products.

\example{Example 3.7}{\it Weight $3$ modular form with respect to $\Gamma_6$.}
For the case of the multiplicative $2$-symmetrisation
of $\Delta_1(Z)$ we obtain
a cusp form constructed in Example  1.20
$$
\ms_2(\Delta_1)=\ml(\phi_{0,3}|T_-(2))
=\hbox{Lift}(\vartheta(\tau,2z)\vartheta(\tau,z)^2\eta(\tau)^3)
\in
\frak N_3(\Gamma_6, v_\eta^{12}\times \hbox{id}_H)
\tag{3.20}
$$
with $\dv_{\Cal A_6^{+}}(F_3^{(6)})=H_4+3H_1$.
Therefore $\ms_2(\Delta_1)$ defines a canonical differential
form with known divisor on the double covering of the moduli
space of Abelian surfaces with polarization of type $(1,6)$.
For more information on this subject see 
\cite{GH2}.
\endexample

We  may construct an interesting modular form using $\Delta_{1/2}(Z)$.
Its $2$-symmetrisation is a modular form of half-integral
weight
$$
F_{3/2}^{(8)}=\ms_2(\Delta_{1/2})=\ml(\phi_{0,4}|T_-(2))
\in
\frak M_{3/2}(\Gamma_8, v_\eta^{9}\times \hbox{id}_H)
\tag{3.21}
$$
with $\dv_{\Cal A_8^{+}}(F_{3/2}^{(8)})=H_4(2)+3H_1$, where
$H_4(2)=\pi_8(\{2z=1\})\subset \Cal A_8^+$.
The square of $F_{3/2}^{(8)}(Z)$ is a modular form of weight 3.

\medskip
The multiplicative $3$-symmetrisation of $\Delta_1(Z)$, $\Delta_2(Z)$
and $\Delta_5(Z)$ produces modular forms with divisor $H_9$.
We obtain the following three functions
$$
\align
F_{16}^{(3)}(Z)&=
\frac{\ms_3(\Delta_5(Z))}{\Delta_1(Z)^4}=
\hbox{Exp-Lift}\bigl(\phi_{0,1}| T_-(3)-4\phi_{0,3}\bigr)
\in
\frak N_{16}(\Gamma_3, v_{\eta}^{8}\times \hbox{id}_H)
\tag{3.22}\\
F_{8}^{(6)}&=
\ms_3(\Delta_2)=
\hbox{Exp-Lift}\bigl(\phi_{0,2}| T_-(3)\bigr)
\in
\frak N_{8}(\Gamma_6)\\
F_4^{(9)}&=\ms_3(\Delta_1)=\ml(\phi_{0,3}|T_-(3))
\in
\frak N_4(\Gamma_9, v_\eta^{16}\times \hbox{id}_H)
\tag{3.23}
\endalign
$$
with the divisors
$$
\dv_{\Cal A_3^{+}}(F_{16}^{(3)})=H_9,
\quad\dv_{\Cal A_6^{+}}(F_{8}^{(6)})=H_9+4H_1,
\quad\dv_{\Cal A_9^{+}}(F_4^{(9)})=H_9(3)+4H_1
$$
where $H_9(3)=\pi_9(\{3z=1\})\subset \Cal A_9^+$.
(We remark that there are two non-equivalent irreducible
Humbert surfaces of discriminant $9$ in $\Cal A_9^+$.)
\medskip

\subhead
3.2. Multiplicative Hecke operators
\endsubhead
\smallskip

Let $F\in \frak M_k(\Gamma_t,\chi)$ be a modular form of integral
weight and
$$
X=\Gamma_t M \Gamma_t=\sum_i\Gamma_tg_i\in
\Cal H_*(\Gamma_t)=\bigotimes_{(p,t)=1} \Cal H_p(\Gamma_t)
\cong \bigotimes_{(p,t)=1} \Cal H_p(\gm_1)
$$
be a Hecke element with a good reduction modulo all primes
dividing $t$. Then we can define {\it a multiplicative
Hecke operator}
$$
[F]_X: =\prod_i F|_k g_i.
\tag{3.24}
$$
This is again a $\Gamma_t$-modular form. We call it
{\it the Hecke product of $F$ defined by $X$}.
In Theorem A.7 of  \cite{GN4, Appendix A}
we proved that the exponential lifting commutes with
multiplicative Hecke operators. More exactly,
for arbitrary $\phi\in J_{0,t}^{nh}$ and
$X\in \Cal H_*(\Gamma_t)$ the identity
$$
[\ml(\phi_{0,t})]_X=
c\cdot \ml\bigl(\phi| \Cal J_0^{(t)}(X)\bigr)
\tag{3.25}
$$
is valid,
where
${\Cal J}^{(t)}_0$ is a natural projection
of the Hecke ring  $\Cal H_*(\Gamma_t)$ into the Hecke-Jacobi ring
of the parabolic subgroup of $\Gamma_t$
(see \cite{G2}, \cite{G6}, \cite{G8})
and $c$ is a constant.
The operator ${\Cal J}^{(t)}_0$ is the same one which appears
in the commutative relation between the arithmetic lifting
and Hecke operators (see \cite{G8}).

We consider below  \thetag{3.25} for the Hecke operator
$T(p)=\gm_t\hbox{diag}(1,1,p,p)\gm_t$
in the case of good and bad reduction.
If $T(p)\in H(\Gamma_t)$ ($(p,t)=1$) we have
$$
\align
T(p)=\gm_t&\Biggl(\smallmatrix
p&0&0&0\\
0&p&0&0\\
0&0&1&0\\
0&0&0&1
\endsmallmatrix\Biggr)
+\sum_{a_1,a_2,a_3\,mod\,p}
\gm_t\Biggl(\smallmatrix
1&0&a_1&a_2\\0&1&a_2&a_3/t\\0&0&p&0\\0&0&0&p
\endsmallmatrix\Biggr)
\\
\vspace{2\jot}
{}&+
\sum_{a\,mod\,p}
\gm_t\Biggl(\smallmatrix
1&0&a&0\\0&p&0&0\\0&0&p&0\\0&0&0&1
\endsmallmatrix\Biggr)+
\sum_{b_1,b_2\,mod\,p}
\gm_t\Biggl(\smallmatrix
p&0&0&0\\-b_1&1&0&b_2/t\\0&0&1&b_1\\0&0&0&p
\endsmallmatrix\Biggr)
\tag{3.26}
\endalign
$$
and
$$
\Cal J_0(T(p))=T_0(p)+p^2+p
\tag{3.27}
$$
where $T_0(p)$ is the Hecke-Jacobi operator
$$
T_0(p)=\sum\Sb b\,mod \,p\\
\vspace{0.5\jot}
c \,mod \,p^2 \endSb
\gi\Biggl(\smallmatrix
1&0&c&b\\
0&p&pb&0\\
0&0&p^2&0\\
0&0&0&p\endsmallmatrix\Biggr)
+\sum\Sb a,b\,mod \,p\\
\vspace{0.5\jot}a\ne 0 \endSb
\gi\Biggl(\smallmatrix
p&0&a&ab\\
0&p&ab&ab^2\\
0&0&p&0\\
0&0&0&p\endsmallmatrix\Biggr)+
\sum_{\lambda\,mod \,p}
\gi\Biggl(\smallmatrix
p^2&0&0&0\\
p\lambda &p&0&0\\
0&0&1&-\lambda\\
0&0&0&p
\endsmallmatrix\Biggr).
\tag{3.28}
$$
The element $T_0(p)$ defines an operator
on the space of Jacobi forms which does not  change the index
$$
T_0(p): J_{k,t}\to J_{k,t}.
$$
It coincides (up to a constant) with the  Hecke operator
$T_p$ defined in \cite{EZ, \S 4}.

We want to use \thetag{3.25} with $X=T(p)$
for the modular forms $\Delta_5(Z)$,
$\Delta_2(Z)$, $\Delta_1(Z)$ and $\Delta_{1/2}(Z)$.
To write down the right hand side of \thetag{3.25},
we need a formula for the action of $T_{0}(p)$
on the  Jacobi forms of index $t$ for a divisor $p$ of $t$.
Let $\phi_{0,t}(\tau,z)$ be an arbitrary Jacobi form
of weight zero and index $t$ with Fourier expansion
$$
\widetilde{\phi_{0,t}}(Z)=\sum_{n,l\in \bz}g(n,l)\,q^nr^ls^t=
\sum_{N} g(N)\exp{(2\pi i\, \hbox{tr}(NZ))}
$$
where $N=\bigl(\smallmatrix n&l/2\\l/2&t
\endsmallmatrix\bigr)\in M_2(\bz)$
($t$ is fixed).
In accordance with \thetag{3.28} we get
$$\gather
\widetilde{\phi_{0,t}}| T_0(p)(Z)=\\
\vspace{2\jot}
\sum_{M_i\in \gi\setminus T_0(p)}
\widetilde{\phi_{0,t}}| M_i (Z)
=
p^3
\sum\Sb l\equiv 0\,mod \,p\\
\vspace{0.5\jot}n\equiv 0 \,mod\,p^2\endSb
g\bigl(\bigl(\smallmatrix n&l/2\\l/2&t\endsmallmatrix\bigr)\bigr)
\exp{\bigl(2\pi i \,\hbox{tr}
(p^{-2}N
\bigl[\bigl(
\smallmatrix
1&0\\0&p
\endsmallmatrix\bigr)\bigr]
Z)\bigr)}\\
\vspace{2\jot}
+
\sum_{n,l\in \bz} G_p(N)
g(N)
\exp{\bigl(2\pi i\, \hbox{tr}(NZ)\bigr)}
+
\sum_{n,l\in \bz}g(N)
\sum_{\lambda\,mod \,p}
\exp{\bigl(2\pi i\, \hbox{tr}
(N
\bigl[\bigl(\smallmatrix p&0\\
\lambda&1\endsmallmatrix\bigr)\bigr]
Z)\bigr)}
\endgather
$$
where $N[M]={}^tMNM$ and we denote
by $G_p(N)$  the Gauss sum
$$
G_p(n,l,t)=-p+\sum\Sb a,b\,mod \,p\endSb
\exp{\biggl(2\pi i \,\frac {na+lab+tab^2}{p}\biggr)}.
$$
Changing $N$ to $N[M^{-1}]$ in the first and third sum,
we get the formula for the Fourier coefficient
$g_p(n,l)$ of Jacobi form $\phi_{0,t}| T_0(p)(\tau,z)$
of index $t$
$$
g_p(n,l)=
p^3 g(p^2 n, pl)+
G_p(n,l,t)g(n,l)+
\sum_{\lambda\,mod \,p}
g({\tsize\frac{n+\lambda l +\lambda^2 t}{p^2}},\,
{\tsize\frac{l +2\lambda t}{p}})
\tag{3.29}
$$
where we put $g(n,l)=0$ if $n\not\in \bz$ or $l\not\in \bz$.
In the case of a good reduction, when $(p,t)=1$,
the $G_p(N)$ is given by the formula
$$
G_p(n,l,t)=p\biggl(\frac {-(4nt-l^2)}p\biggr),\qquad(p,t)=1
$$
(see \thetag{3.19} about  definition of the Kronecker symbol for
$p=2$). If $p|t$, we can represent the formula for
$G_p$ in the form useful for exact calculations
$$
G_p(n,l,t)=
p\,\cases
\hphantom{p-,}0&\text{ if }\  l\not\equiv 0\ \hbox{mod }p\\
p-1&\text{ if }\  (n,l)\equiv (0,0)\ \hbox{mod }p\\
\hphantom{p}-1&\text{ if }\  n\not\equiv 0\ \hbox{mod }p,\
\text{ and }\
l\equiv 0\ \hbox{mod }p.
\endcases
$$
For the Fourier coefficients
depending on $\lambda$ in \thetag{3.29} we have
$\hbox{det}\bigl( N
\bigl[\bigl(\smallmatrix p&0\\
\lambda&1\endsmallmatrix\bigr)\bigr]\bigr)=
p^2\hbox{det}( N)$. If $(t,p)=1$, then
$n+\lambda l +\lambda^2 t$ is a full square mod $p$
for $4nt-l^2\equiv 0\ \hbox{mod\,}p^2$.
Thus  there exists only
one $\lambda$ mod $t$ which gives us a non-trivial
term in the third summand in \thetag{3.29}.
Therefore we prove

\proclaim{Lemma 3.8} Let us suppose that the Fourier coefficient
$g(n,l)$
of the Jacobi form $\phi_{0,t}$ of weight zero and index $t$
depends only on the norm $N=4nt-l^2\in \bz$. We denote
$g(N)=g(n,l)$. For any prime $p$ such that
$(t,p)=1$, the Fourier coefficients of
$\phi_{0,t}| T_0(p)$ are given by the formula
$$
g_p(N)=p^3g(p^2N)+ p\biggl(\frac {-N}p\biggr)g(N)+
g\biggl(\frac {N}{p^2}\biggr)
$$
where we set $g(\frac {N}{p^2})=0$ if $\frac {N}{p^2}\not\in \bz$.
\endproclaim

\example{Example 3.9}{\it The Igusa modular form $\Delta_{35}(Z)$.}
The Igusa modular form $\Delta_{35}(Z)$ is
the first Siegel modular form of odd weight with respect to
$\gm_1=Sp_4(\bz)$. It has weight $35$ (see \cite{Ig1}).
We defined this modular form in \cite{GN4} as a
Hecke product of $\Delta_5(Z)$.
In this example we recall this construction.
Let us take the modular form $\Delta_5(Z)$ which has
the divisor $H_1$ in $\Cal A_1$. Using
the system  of representatives $T(p)$,  we then get
$$
\align
[\Delta_5(Z)]_{T(2)}=
\kern-10pt\prod_{a,b,c\,mod\, 2}&\kern-8pt
\Delta_5 ({\tsize \frac{z_1+a}2 ,
\tsize\frac{z_2+b}2, \tsize\frac{z_3+c}2})
\prod_{a \,mod\, 2}
\Delta_5 ({\tsize\frac{z_1+a}2,z_2,2z_3})\,
\Delta_5 ({\tsize 2z_1, z_2, \frac{z_3+a}2})\\
\vspace{2\jot}
{}&\times\Delta_5 ({\tsize 2z_1, 2z_2, 2z_3})
\prod_{b\,mod\, 2}
\Delta_5 ({\tsize 2z_1, -z_1+z_2, \frac{z_1-2z_2+z_3+b}2}).
\tag{3.30}
\endalign
$$
One can check that
$\hbox{div}_{\Cal A_1}\bigl([\Delta_5(Z)]_{T(p)}
\bigr)=(p+1)^2H_1+H_{p^2}$.
Thus
$$
\Delta_{35}(Z)=\frac{[\Delta_5(Z)]_{T(2)}}{\Delta_5(Z)^8}
=\hbox{Exp-Lift}
\bigl(\phi_{0,1}| (T_0(2)-2) \bigr)\in
\frak N_{35}(\Gamma_1)
$$
and
$$
\Delta_{35}(Z)=
q^2 r s^2\,(q-s)
\prod\Sb n,\,l,\,m\in \Bbb Z\\
\vspace{0.5\jot}
(n,l,m)>0\,\endSb
\bigl(1-q^n r^l s^m\bigr)^{f_1^{(2)}(4nm-l^2)}
\tag{3.31}
$$
where $f_1(4n-l^2)=f_1(n,l)$ are the Fourier coefficients
of $\phi_{0,1}(\tau,z)$
and
$$
f_1^{(2)}(N)=8f_1(4N) +2(\biggl(\frac{-N}2\biggl)-1)f_1(N)
+f_1(\frac {N}4)
$$
according to Lemma 3.8.
Remark that we cannot construct $\Delta_{35}(Z)$
as an arithmetic lifting of a holomorphic Jacobi form.
Nevertheless \thetag{3.29} gives us $\Delta_{35}(Z)$
as a finite Hecke product of the lifted form $\Delta_5(Z)$.
\endexample

\example{Example 3.10}{\it The modular form $D_{6}(Z)$.}
Using \thetag{3.25}--\thetag{3.26} and Lemma 3.8, we get
modular forms with divisor $H_{p^2}$ in $\Cal A_2^{+}$ and
$\Cal A_3^{+}$  respectively
$$
\align
F_p^{(2)}(Z)&=c_2\frac{[\Delta_2(Z)]_{T(p)}}{\Delta_2(Z)^{(p+1)^2}}
=\hbox{Exp-Lift}
\bigl(\phi_{0,2}| (T_0(p)-p-1)\bigr)\in
\frak N_{2p(p^2-1)}(\Gamma_2),\\
F_p^{(3)}(Z)&=c_3\frac{[\Delta_1(Z)]_{T(p)}}{\Delta_1(Z)^{(p+1)^2}}
=\hbox{Exp-Lift}
\bigl(\phi_{0,3}| (T_0(p)-p-1) \bigr)\in
\frak N_{p(p^2-1)}(\Gamma_3,\chi_3^{(p)})
\endalign
$$
where $p\ne 2$  for $\Delta_2$
and  $p\ne 3$ for $\Delta_1$.
The character $\chi_2^{(p)}$
is trivial or has order two.
\newline
(The $c_p$ are constants which can be  easily calculated.)
In particular  we get the modular form
$$
D_6(Z)=\frac{2^{22}{[\Delta_1(Z)]_{T(2)}}}{\Delta_1(Z)^{9}}
=\hbox{Exp-Lift}
\bigl(\phi_{0,3}| (T_0(2)-3)\bigr)
$$
which we have constructed  as arithmetic lifting
of the  Jacobi form  $\eta^{11}\vartheta_{3/2}$
in Example 1.17.
Thus according  to \thetag{3.9}
we get
$$
D_6(Z)=\hbox{Lift}(\eta^{11}\vartheta_{3/2})
=q^{\frac {1}2}r^{\frac {1}2}s^{\frac {3}2} \prod
\Sb n,\,l,\,m\in \Bbb Z\\
\vspace{0.5\jot}
(n,l,m)>0 \endSb
\bigl(1-q^n r^l s^{3m}\bigr)^{g_{3}(nm,l)}
\in
\frak N_6(\Gamma_3^+, v_\eta^{12}\times v_H)
\tag{3.32}
$$
where
$$
\phi_{0,3}^{(6)}(\tau,z)=\phi_{0,3}| (T_0(2)-3)(\tau,z)
=\sum_{n,l}g_3(n,l)q^nr^l
=r^2-r+12-r^{-1}+r^{-2}+q(\dots).
$$
\endexample

\example{Example 3.11}{\it Hecke product for $p=t=2$
and $p=t=3$}.
Using \thetag{3.29} we can consider the cases of bad reduction
$p=t=3$ and $p=t=3$.
The Hecke operator
$T^+(2)=\gm_2^+\hbox{diag}(1,1,2,2)\gm_2^+$ from the
Hecke ring $H(\Gamma_2^+)$ of the maximal normal extension
$\Gamma_2^+$ contains $18$ left cosets:
the $15$ cosets from \thetag{3.26} and
$$
\sum\Sb a,b \,mod\, 2\\
\vspace{0.5\jot}(a,b)\ne (0,0)\endSb
\gm_2^+
\pmatrix
-a&2&b&0\\
1&0&0&b/2\\
0&0&0&1\\
0&0&2&a
\endpmatrix.
$$
Therefore the modular form
$[\Delta_2]_{T^+(2)}$ of weight $36$ has divisor $2H_4+9H_1$
and we obtain the identity
$$
\align
\Delta_{11}(Z)^2&=
\bigl(\hbox{Lift}(\eta^{21}(\tau)\vartheta(\tau,2z)\bigr)^2
=\bigl(\ml(\phi_{0,2}^{(11)}(Z))\bigr)^2\\
{}&=c[\Delta_2]_{T^+(2)}(Z)/\Delta_2(Z)^7
=\ml(\phi_{0,2}| T_0(2))(Z)
\endalign
$$
because
$$
\phi_{0,2}| T_0(2)(\tau,z)=
2r^{2}+44+2r^{-2}+q(\dots).
$$
As a corollary we can continue  the identity \thetag{3.14}
between the basic Jacobi forms
(see Lemma 3.5)
$$
\phi_{0,2}^{(11)}=
\phi_{0,1}| T_-(2)-2\phi_{0,2}=\phi_{0,1}^2-20\phi_{0,2}
=\frac{1}2\phi_{0,2}| T_0(2).
\tag{3.33}
$$

Similarly to the case $p=2$  we obtain  for $p=t=3$ the Jacobi form
$
\phi_{0,3}| T_0(3)(\tau,z)=
2r^{3}+68+2r^{-3}+q(\dots)
$.
Therefore the Hecke product for $p=3$ of $\Gamma_3$-modular form
$\Delta_1$ is related with
$$
F_{16}^{(3)}(Z)=\ml{\bigl(\frac{1}2\phi_{0,3}| (T_0(3)-2)
\bigr)}\in
\frak N_{16}(\Gamma_3, v_\eta^{8}\times \hbox{id}_H)
$$
with divisor $H_9$.
We have constructed this modular form in \thetag{3.22}
using multiplicative
$3$-symmetrisation of $\Delta_5$. In terms of Jacobi forms it is
equivalent to the relation
$$
\phi_{0,3}|T_0(3)=2\phi_{0,1}|T_-(3)-6\phi_{0,3}.
$$
Recall that the Jacobi form $\phi_{0,3}(\tau,z)$ is
the square of an infinite product (see \thetag{2.19}).
\endexample

\example{Example 3.12}{\it $\Gamma_4$-modular forms with
$H_{p^2}$-divisors}.
Our next examples are connected with  the Jacobi form $\phi_{0,4}$.
In the case of good reduction ($p\ne 2$), we can define the same
function as in Example 3.10.
Since
$$
\phi_{0,4}| T_0(p)(\tau,z)=
r^p+pr+(p^3+1)+r^{-p}+pr^{-1}+q(\dots) \qquad(p\ne 2)
$$
we obtain a modular form
$$
F_p^{(4)}(Z)=\ml\bigl(\phi_{0,4}| (T_0(p)-p-1)\bigr)
\in
\frak N_{p(p^2-1)/2}(\Gamma_4,v_{\eta}^{p(p^2-1)}
\times \hbox{id}_H)
$$
with divisor $H_{p^2}$.
For instance for $p=3$ we get a modular form of weight $12$ with
divisor $H_9$ in $\Cal A_4^+$.
Using \thetag{3.29} for $p=2$, we obtain a very nice identity
$$
\phi_{0,4}| T_0(2)(\tau,z)=\phi_{0,1}(\tau,2z).
\tag{3.34}
$$
It gives us the second new formula for the generator $\phi_{0,1}(\tau,z)$
(see \cite{3.16}).
We can represent it in an equivalent form using
the operator $\Lambda_2^*$ (see \thetag{3.17}).
One can  check  (see \cite{G2}) that
$$
\align
p^{-1}T_0(p)\cdot\Lambda_2^*=T_+(1,p^2)&=
\sum\Sb a,b,c\,mod \,p^2 \endSb
\gi\Biggl(\smallmatrix
1&0&a&pb\\
0&1&b&c\\
0&0&p^2&0\\
0&0&0&p^2\endsmallmatrix\Biggr)
+\sum\Sb a,b\,mod \,p\\
\vspace{0.5\jot}
 c\,mod\, p^2\\
\vspace{0.5\jot}a\ne 0 \endSb
\gi\Biggl(\smallmatrix
p&0&a&pb\\
0&1&b&c\\
0&0&p&0\\
0&0&0&p^2\endsmallmatrix\Biggr)\\
{}&+
\sum\Sb a\,mod \,p\\ c\,mod\,p^2\endSb
\gi\Biggl(\smallmatrix
p^2&0&0&0\\
-a &1&0&c\\
0&0&1&-a\\
0&0&0&p
\endsmallmatrix\Biggr).
\endalign
$$
Thus \thetag{3.34} is equivalent to
$$
8\phi_{0,1}(\tau,z)=\phi_{0,4}| T_+(1,4)(\tau,z).
\tag{3.35}
$$
We recall that $\phi_{0,4}(\tau,z)$ is given by
an infinite product (see \thetag{2.10}).
One can find a formula for the action of $T_+(1,4)$ on Fourier
coefficients of Jacobi forms similar to
the formula \thetag{3.19} for $T_+(2)$.
\endexample

\head
\S~4. Reflective modular forms
\endhead

In this section we construct some other modular forms
with known divisors and prove new identities
between arithmetic and exponential liftings.
Our main aim is to construct reflective  modular forms,
i.e. forms  with
divisors   defined by reflections in the corresponding
orthogonal group. We are especially interesting in
modular forms of divisors with multiplicity one.

\subhead
4.1. The modular form $D_2(Z)$
\endsubhead
In what follows we shall  use new  weak Jacobi  forms
in order to prove that some modular forms constructed
in \S 1 as the arithmetic lifting have infinite product expansion.
Using the basic weak Jacobi form from \S 2
(see \thetag{2.15}, \thetag{2.19}, \thetag{2.23}),
we define a weak Jacobi form
of weight 0 and index $9$
$$
\multline
\phi_{0,9}(\tau,z)=
\phi_{0,1}(\tau,3z)+7\phi_{0,3}(\tau,z)\xi_{0,6}(\tau,z)-
\phi_{0,3}(\tau,z)^3 \\
=(r^2-r+4-r^{-1}+r^{-2})+\dots+q^2(r^9+\dots).
\endmultline
\tag{4.1}
$$
Since $\eta^6\phi_{0,9}$ is holomorphic, the  part of the Fourier
expansion written above contains all Fourier coefficients
of negative norm. Using this Jacobi form we can construct the function
$D_2(Z)$ (see \thetag{1.30}) from Example 1.17
as exponential lifting.

\proclaim{Theorem 4.1}The following identity is valid:
$$
\align
D_2(Z)&=\sum_{N\ge 1}
\sum\Sb m >0,\,l\in \bz\\
\vspace{0.5\jot} n,\,m\equiv 1\,mod\,6\\
\vspace{0.5\jot} 4nm-l^2=N^2\endSb
\hskip-4pt
N\biggl(\dsize\frac{-4}{N}\biggr)\biggl(\dsize\frac{12}{l}\biggr)
\sum\Sb a|(n,l,m)\endSb \biggl(\dsize\frac{6}{a}\biggr)
q^{n/6}r^{l/2}s^{3m/2}\\
{}&=
q^{\frac {1}6}r^{\frac {1}2}s^{\frac {3}2} \prod
\Sb n,\,l,\,m\in \Bbb Z\\
\vspace{0.5\jot}
(n,l,m)>0\,\endSb
\bigl(1-q^n r^l s^{9m}\bigr)^{f_{9}(nm,l)}
\in \frak N_{2}(\gm_{9}, v_{\eta}^4\times v_H).
\tag{4.2}
\endalign
$$
Moreover $D_2(Z)$ has the  divisor
$$
\hbox{Div}_{\Cal A_9^+}(D_2(Z))=H_4(2)+H_9(9)
$$
(see \thetag{1.23}).
\endproclaim
\demo{Proof}The lifted form $D_2(Z)$ has zero along $H_4(2)$
according to Lemma 1.16.
To prove the identity  we have to check that
$D_2(Z)$ has zero along $H_9(9)$. This Humbert surface
is defined by the element
$\ell_9=2f_2+\frac 1{2}f_3+f_2\in \widehat {L}_9$
(see \thetag{1.14} and \thetag{1.22}). In the basis \thetag{1.14}
the matrix of  the reflection  $\sigma_{\ell_9}$  has
the form
$$
(\sigma_{\ell_9})
=\pmatrix
9&-72&16\\
2&-17&8\\
4&-36&9
\endpmatrix.
$$
Thus, if we consider $D_2$ as a modular form with respect
to the  orthogonal group of lattice $L_9$, we get
$$
\gather
 D_{2}(\sigma_{l_9}(\frak z))=
\sum_{n,l, m}
a(n,l,m)
\exp{\bigl(2\pi i (\tsize\frac {9n-24l+16m}{6} z_1+
\tsize\frac{12n-17l+24m}2 z_2
+\tsize\frac{2n-6l+9m }{6} z_3))}\\
=-D_{2}(\frak z),
\endgather
$$
because the Fourier coefficient $a(n,l,m)$
has the  factor $\biggl(\dsize\frac{12}{l}\biggr)$.
\newline
\qed
\enddemo
Similarly to $\phi_{0,36}(\tau,z)$ (see \thetag{2.12})
we introduce a weak Jacobi form  of index $12$
$$
\xi_{0,12}(\tau,z)=\frac {\vartheta_{3/2}(\tau,3z)}
{\vartheta_{3/2}(\tau,z)}=
\frac {\vartheta(\tau,6z)\vartheta(\tau,z)}
{\vartheta(\tau,3z)\vartheta(\tau,2z)}.
\tag{4.3}
$$
Using the same arguments as in the case of
the weak Jacobi form $\phi_{0,36}$, we see that one needs to calculate
the Fourier coefficients of $\phi_{0,12}(\tau,z)$ for  $q^n$
with $n\le 3$. It gives us:
\proclaim{Lemma 4.2}The weak Jacobi form
$$
\align
\xi_{0,12}(\tau,z)&=
(r-1+r^{-1})\prod_{n\ge 1}(1-q^nr^{\pm 1}+q^{2n}r^{\pm 2})
(1+q^{2n-1}r^{\pm 2}+q^{4n-2}r^{\pm 4})\\
{}&
(1+q^{3n-2}r^{\pm 3})(1+q^{3n-1}r^{\pm 3})
(1-q^{6n-1}r^{\pm 6})(1-q^{6n-5}r^{\pm 6})
\endalign
$$
($\pm m$ means that there are two factors in the product:
 with $-m$ and with $+m$)
has a non-zero Fourier coefficient $f_{12}(n,l)$
of negative norm only for $48n-l^2=-1$ ($r$ and $-qr^7$) and
$48n-l^2=-4$ ($-q^2r^{10}$) and its first Fourier coefficients are
$$
\xi_{0,12}(\tau,z)
=(r-1+r^{-1})-q(r^7+\dots)-q^2(r^{10}+\dots)
-q^3(r^{12}+\dots)+\dots\, .
$$
In particular, $\eta^2\xi_{0,12}$ is holomorphic.
\endproclaim
Using $\xi_{0,12}$ we define
$$
\multline
\phi_{0,18}(\tau,z)=\xi_{0,12}(\tau,z)\cdot
\xi_{0,6}(\tau,z)=
\frac {\vartheta(\tau,6z)\vartheta(\tau,4z)\vartheta(\tau,z)}
{\vartheta(\tau,3z)\vartheta(\tau,2z)^2}
\\
=(r^2-r+2-r^{-1}+r^{-2})+\dots+q^3(r^{15}+\dots)
+q^4(2r^{17}+\dots)+\dots\, .
\endmultline
\tag{4.4}
$$
Since $\eta^3\phi_{0,18}$ is holomorphic, one needs  to know
only the Fourier coefficients with norm $\ge -9$.
All Fourier coefficients of this type are written in
the Fourier expansion above.
We give, without proof, the identity for the modular
form $D_1(Z)$ (see \thetag{1.29}) from Example 1.17:
$$
\align
D_1(Z)&=\sum_{M\ge 1}
\sum\Sb m >0,\,l\in \bz\\
\vspace{0.5\jot} n,\,m\equiv 1\,mod\,12\\
\vspace{0.5\jot} 4nm-l^2=M^2 \endSb
\hskip-4pt
\biggl(\dsize\frac{12}{Ml}\biggr)
\sum\Sb a|(n,l,m)\endSb
\biggl(\dsize\frac{-4}{a}\biggr)
q^{n/12}r^{l/2}s^{3m/2}\\
{}&=
q^{\frac {1}{12}}r^{\frac {1}2}s^{\frac {3}2}
\prod
\Sb n,\,l,\,m\in \Bbb Z\\
\vspace{0.5\jot}
(n,l,m)>0\,\endSb
\bigl(1-q^n r^l s^{18m}\bigr)^{f_{18}(nm,l)}
\in \frak N_{1}(\gm_{18}, v_{\eta}^2\times v_H).
\tag{4.5}
\endalign
$$
The weak Jacobi forms introduced above provide another formula
for the Jacobi form $\phi_{0,36}$
(see \thetag{2.12}) used in
the exponential lifting of the singular modular form $D_{1/2}(Z)$
$$
\phi_{0,36}(\tau,z)=\phi_{0,4}(\tau,3z)-\xi_{0,6}(\tau,2z)
\xi_{0,12}(\tau,z).
\tag{4.6}
$$

\subhead
4.2. Anti-symmetric modular forms
\endsubhead

The arithmetic lifting provides us with modular forms
which are  invariant with respect to the main exterior
involution $V_t$ of the group $\Gamma_t$
(see \thetag{1.19}).
Using the  exponential lifting, one can construct
anti-invariant modular forms, i.e. forms satisfying
$F(V_t(Z))=-F(Z)$. For example for  $t=1$ we have the anti-invariant
form  $F(Z)=\Delta_{35}(Z)$.
In this subsection we construct anti-invariant
modular forms for arbitrary $\Gamma_t$ ($t>1$).
For  $\Gamma_2$, $\Gamma_3$ and
$\Gamma_4$  these modular forms  have only the Humbert surface
$H_{4t}(0)=\pi_t\{\tau-t\omega=0\}$ as their  divisor.
(We remark that for $t=4$ there are two
irreducible  Humbert surfaces with discriminant $16$.)

Let us consider the function
$
\psi_{0,t}(\tau,z)=\Delta(\tau)^{-1}E_{12,t}(\tau,z)
$
where
$$
\Delta(\tau)=q\prod_{n\ge 1}(1-q^n)^{24}=
q-24q^2+253q^3+\dots
$$
and $E_{12,t}(\tau,z)$ is
a Jacobi--Eisenstein series of weight $12$ and index $t$.

There exists a formula for Fourier coefficients of $E_{k,1}$
in terms of H. Cohen's numbers (see \cite{EZ, \S 2}).
One can find  the table of the values of  Fourier
coefficients of $E_{4,1}(\tau,z)$ and $E_{6,1}(\tau,z)$ in
\cite{EZ, \S 1}.
Using the basic Jacobi forms
$\phi_{0,2}(\tau,z)=r^{\pm 1}+4+\dots\ $,
$\phi_{0,3}(\tau, z)=r^{\pm 1}+2+\dots\ $
and the forms
$\phi_{0,2}^{(11)}=r^{\pm 2}+22+\dots$ and
$\phi_{0,3}^{(6)}=
r^{\pm 2}-r^{\pm 1}+12+\dots$,
which are the data  for the exponential liftings
\thetag{3.11} and \thetag{3.32},
we define
$$
\align
\psi_{0,2}(\tau,z)&=\Delta(\tau)^{-1}E_{6,1}(\tau,z)^2
-2\phi_{0,2}^{(11)}(\tau,z)+176\phi_{0,2}(\tau,z)\\
{}&=\sum_{n\ge 0,\, l\in \bz}c_2(n,l)q^nr^l=
q^{-1}+24+q(\dots),
\\
\psi_{0,3}(\tau,z)&=\Delta(\tau)^{-1}E_{4,1}(\tau,z)^3
-3\phi_{0,3}^{(6)}(\tau,z)-171\phi_{0,3}(\tau,z)\\
{}&=\sum_{n\ge 0,\, l\in \bz}c_3(n,l)q^nr^l=
q^{-1}+24+q(\dots).
\endalign
$$
The Jacobi forms $\psi_{0,p}$ ($p=2$, $3$) contain the only
type of  Fourier coefficients with indices of negative norm.
This is $q^{-1}$ of norm $-4p$.
Thus we can use both functions to produce the exponential
liftings
$$
\align
\Psi_{12}^{(2)}(Z)&=\ml(\psi_{0,2})=
q\prod\Sb n,l,m\in \bz\\
\vspace{0.5\jot}
(n,l,m)>0\endSb
 (1-q^nr^ls^{2m})^{c_2(nm,l)}
\in \frak M_{12}(\Gamma_2),
\tag{4.7}
\\
\Psi_{12}^{(3)}(Z)&=\ml(\psi_{0,3})=
q\prod\Sb n,l,m\in \bz\\
\vspace{0.5\jot}
(n,l,m)>0\endSb
 (1-q^nr^ls^{3m})^{c_3(nm,l)}
\in \frak M_{12}(\Gamma_3).
\tag{4.8}
\endalign
$$
According to Theorem 2.1
$$
\Psi_{12}^{(p)}(V_p<Z>)=-\Psi_{12}^{(p)}(Z)\quad (p=2,\,3)
\quad\text{and }\
\dv_{\Cal A_p}(\Psi_{12}^{(p)})=
\cases H_8&\text{for } p=2\\
       H_{12}&\text{for } p=3.
\endcases
$$
The  Fourier-Jacobi expansion of $\Psi_{12}^{(p)}$
starts with coefficients
$$
\Psi_{12}^{(p)}(Z)=\Delta_{12}(\tau)-\Delta_{12}(\tau)\psi_{0,p}(\tau,z)
\exp{(2\pi i p \omega)}+\dots .
$$
Therefore the constructed modular forms $\Psi_{12}^{(p)}(Z)$ 
($p=2$, $3$) are not cusp forms.

If we do the same for $t=4$, we get a Jacobi form
we used to construct $\Delta_{35}(Z)$. Let us take
the  Jacobi form
$$
\phi_{0,1}|(T_0(2)-2)(\tau, 2z)=q^{-1}+(r^4+70+r^{-4})+q(\dots).
$$
Its exponential lifting is zero along two Humbert
surfaces with discriminant $16$. To delete the second  component,
we consider the additional
Jacobi--Eisenstein series which has the constant term equals zero
(such a series exists if the index contains a perfect square).
For $t=4$ this Jacobi--Eisenstein series is
the eight power of the Jacobi theta-series $\vartheta(\tau,z)$.
\remark{Remark 4.3}It is easy to modify the calculation
of Fourier coefficients in \cite{EZ} to find
an exact formula for the Fourier expansion
of the additional Eisenstein series with
a perfect square index. It gives us, for example,
a formula for the eight power of the triple product
in terms of H.Cohen numbers. We hope to present
such formulae later.
\endremark

Using $\vartheta(\tau,z)^8$, we define
$$
\align
\psi_{0,4}(\tau,z)&=
\bigl(\phi_{0,1}|(T_0(2)+26)\bigr)(\tau, 2z)-
\Delta(\tau)^{-1}E_4(\tau)\vartheta(\tau,z)^8-
8\bigl(\phi_{0,4}|(T_0(3)+4)\bigr)(\tau,z)\\
\vspace{2\jot}
{}&=\sum_{n\ge 0,\, l\in \bz}c_4(n,l)q^nr^l=
q^{-1}+24+q(\dots).
\endalign
$$
Similarly to $\psi_{0,2}$ and $\psi_{0,3}$
the Jacobi form $\psi_{0,4}$ contains  only
the  Fourier coefficients of type $q^{-1}$
with index of negative norm.
Taking its exponential lifting we obtain
the $\Gamma_4$-modular form of weight $12$
$$
\Psi_{12}^{(4)}(Z)=\ml(\psi_{0,4})(Z)=
q\prod\Sb n,l,m\in \bz\\
\vspace{0.5\jot}
(n,l,m)>0\endSb
 (1-q^nr^ls^{4m})^{c_4(nm,l)}
\in \frak M_{12}(\Gamma_4).
\tag{4.9}
$$
According to Theorem 2.1, $\Psi_{12}^{(4)}(Z)$
is anti-invariant and $\dv_{\Cal A_4}(\Psi_{12}^{(4)})=
H_8(0)$.

\remark{Remark 4.4} The fact that all three modular forms
constructed above have weight $12$ has the following explanation
using the   Borcherds singular modular form
$\Phi_{12}(\tau,\frak z,\omega)$
with respect to the  orthogonal group of the lattice $L$ of signature
$(26,2)$, which is  the orthogonal sum of two
unimodular hyperbolic plains $U$ and the Leech lattice $\Cal L$.
See \cite{Bo5} and  also \cite{Bo4}, \cite{Bo2}.
The modular form $\Phi_{12}(\tau,\frak z,\omega)$
is anti-invariant, i.e.
$$
\Phi_{12}(\tau,\frak z,\omega)=-\Phi_{12}(\omega,\frak z,\tau).
$$
Let us take a primitive $\ell \in \Cal L$ such that
$\ell^2=2t$ ($t>1$).
Then the orthogonal group
of the lattice $2U\oplus<2t>$ (see the notation in Sect. 1.3)
of signature $(3,2)$ acts on the tree dimensional
subdomain $\bh_t^+$ of the tube  homogeneous domain of dimension
$26$.
The  $\bh_t^+$ is isomorphic to the Siegel
upper half plane  and the restriction $\Phi|_{\bh_t^+}$
is a Siegel modular form of weight $12$
with respect to the group $\Gamma_t$ (see \thetag{1.21} and Lemma 1.9).
Thus  {\it for arbitrary $t>1$ there exists
an anti-symmetric modular form of weight $12$
with respect to $\Gamma_t$.}
Since $\Psi_{12}^{(t)}(Z)$ ($t=2$, $3$, $4$) have weight $12$
and the divisor $H_{4t}$, they coincide with
$c\Phi|_{\bh_t^+}$.
It also proves that {\it the divisor of $\Phi|_{\bh_t^+}$
is exactly equal to $H_{4t}$ for} $t=2$, $3$, $4$.
For arbitrary $t$ the restriction
$\Phi|_{\bh_t^+}$ may have additional divisors.
\endremark

\subhead
4.3. The reflective modular forms with divisors
$H_1$ and $H_4$
\endsubhead
\example{Example 4.5} {\it The first cusp form for} $\Gamma_5$.
In \thetag{1.32} we defined a cusp form of weight $5$
with respect to $\Gamma_5$ with trivial character.
Let us construct this function
as the exponential lifting. We set
$$
\phi_{0,5}(\tau,z)=\phi_{0,2}(\tau,z)\phi_{0,3}(\tau,z)=
(r^{2}+6r^{1}+10+6r^{-1}+r^{-2})+q(-3r^4+\dots)+\dots\, .
$$
Its Fourier expansion does not contain the summand $qr^5$ with
$N=20n^2-l^2=-5$. Thus
$$
\hbox{Div}_{\Cal A_5^+}(\hbox{Exp-Lift}(\phi_{0,5}))
=7H_1+H_4.
$$
According to Lemma 1.16
$$
F_5^{(5)}(Z)
=\hbox{Lift}(\eta(\tau)^3\vartheta(\tau,z)^6 \vartheta(\tau,2z))=
\hbox{Exp-Lift}(\phi_{0,5})
\in \frak N_5(\Gamma_5).
\tag{4.10}
$$
We can define another weak Jacobi form of index $5$
$$
\phi_{0,5}^{(2)}(\tau,z)=2\phi_{0,2}(\tau,z)\phi_{0,3}(\tau,z)
-\phi_{0,1}(\tau,z)\phi_{0,4}(\tau,z)=
(r^{\pm 2}+r^{\pm 1}+8)+q(r^5+\dots)+\dots\, .
$$
According to Theorem 2.1
$$
\hbox{Div}_{\Cal A_5^+}(\hbox{Exp-Lift}(\phi^{(2)}_{0,5}))
=2H_1+H_4+H_5.
$$
One can prove that
$$
G_4^{(5)}(Z)
=\hbox{Lift}(\eta(\tau)^6\vartheta(\tau,z) \vartheta(\tau,2z))
=\hbox{Exp-Lift}(\phi_{0,5}^{(2)})
\in \frak N_4(\Gamma_5, v_\eta^{12}\times v_H).
\tag{4.11}
$$
\endexample

\example{Example 4.6}{\it The cusp forms for} $\Gamma_6$.
If the polarization $t$ contains more than one of primes,
the structure of the  graded ring of  modular forms
reflects the properties of rings of modular forms
for different polarizations.
In \cite{G2} it was proved injectivity of
the $p$-symmetrisation
$
\hbox{Sym}_{t,p}: \frak M_{k}(\Gamma_t)\to \frak M_{k}(\Gamma_{tp})
$
if $(t,p)=1$ (see \thetag{3.1}).
In this example we present infinite product expansions
of  the $\Gamma_6$-cusp forms of weight $3$ constructed
in Example 1.15 and 1.20.

According to the dimension formula of the space of Jacobi forms
(see \cite{EZ}, \cite{SkZ})
we have
$
\hbox{dim}(J_{6,6}^{cusp})=1
$.
A square root of a generator of this space was defined
in Example 1.20:
$
\phi_{6,6}(\tau,z)
=(\eta(\tau)^3\vartheta(\tau,z)^2\vartheta(\tau,2z))^2.
$
Its arithmetic lifting is
$$
F_3^{(6)}(Z)
=\hbox{Lift}(\eta(\tau)^3\vartheta(\tau,z)^2\vartheta(\tau,2z))
\in \frak N_3(\Gamma_6^+, v_\eta^{12}\times \hbox{id}_H).
$$
For $(1,6)$-polarization there are two
Humbert surfaces with discriminant $1$
in $\Cal A_6^+$:
$$
H_1(1)=\pi_6(\{Z\in \bh_2\,|\, z=0\}),\qquad
H_1(5)=\pi_6(\{Z\in \bh_2\,|\, \tau+5z+6\omega=0\}),
$$
and one Humbert surface $H_4$ of discriminant $4$.
According to Lemma 1.16
$\dv_{\Cal A_6^+}(F_3^{(6)})$ contains $3H_1(1)+H_4$.
One can prove that
$$
F_3^{(6)}(Z)=\hbox{Lift}\bigl(
\eta(\tau)^3\vartheta(\tau,z)^2\vartheta(\tau,2z)\bigr)
=\hbox{Exp-Lift}
(3\phi_{0,3}^2-2\phi_{0,2}\phi_{0,4})
\tag{4.12}
$$
and
$$
\hbox{Div}_{\Cal A_6^+}(F_3^{(6)}(Z))
=3H_1(1)+2H_1(5)+H_4.
$$

In Example 1.15 we defined $(-1)$-Lift of the Jacobi form
$\eta(\tau)^5\vartheta(\tau,2z)$.
One can prove that
$$
\align
\hbox{Lift}(\eta(\tau)^5\vartheta(\tau,2z))&
=\hbox{Exp-Lift}
(5\phi_{0,3}^2-4\phi_{0,2}\phi_{0,4})
\in \frak N_3(\Gamma_6^+, v_\eta^{8}\times \hbox{id}_H),
\tag{4.13}\\
\hbox{Lift}_{-1}(\eta(\tau)^5\vartheta(\tau,2z))
&=\hbox{Exp-Lift}(\phi_{0,3}^2)
\in \frak N_3(\Gamma_6^+, v_\eta^{16}\times \hbox{id}_H),
\tag{4.14}
\endalign
$$
and their divisors are equal to
$2H_1(1)+4H_1(5)+H_4$ and
$5H_1(1)+H_4$ respectively.
\endexample

\example{Example 4.7} {\it The first cusp form for} $\Gamma_7$.
For the $\Gamma_7$-modular form of weight $2$
from Example  1.20 we have
$$
F_2^{(7)}(Z)=\hbox{Lift\,}(\vartheta(\tau,z)^3\vartheta(\tau,2z))
=
\hbox{Exp-Lift}(\phi_{0,3}\cdot\phi_{0,4})(Z)
\in \frak N_{2}(\Gamma_{7},v_\eta^{12}\times v_H)
\tag{4.15}
$$
with
$$
\hbox{Div}_{\Cal A_7^+}(F_2^{(7)})
=4H_1+H_4.
$$
This divisor is a part of the divisor of the arithmetic lifting
(see Lemma 1.16). Since $\eta^4\phi_{0,3}\phi_{0,4}$
is holomorphic, we have to check in the Fourier expansion
of $\phi_{0,3}\phi_{0,4}$ only Fourier coefficients
with norm $\ge -4$. We see that
$$
\phi_{0,3}(\tau,z)\phi_{0,4}(\tau,z)=
(r^2+3r+4+3r^{-1}+r^{-2})+\dots\,,
$$
thus its exponential lifting has the divisor mentioned above.
This finishes the proof of \thetag{4.15}.

The first non-zero cusp form for $\Gamma_7$ has weight $4$
(see \cite{G2, \S 3}).
The Jacobi form
$\phi_{4,7}(\tau,z)=\vartheta(\tau,z)^6\vartheta(\tau,2z)^2$
is the cusp form of weight $4$ and index $7$, and there exists
the only one such form.
Moreover,
$$
F_4^{(7)}=
\hbox{Lift\,}(\phi_{4,7})=(F_2^{(7)})^2
$$
is the first cusp form with respect to  $\Gamma_7$
with trivial character.
\endexample
\example{Example 4.8}{\it The cusp form of weight $1$ for } $\Gamma_{10}$.
We set
$$
\align
\phi_{0,10}(\tau,z)&=\phi_{0,4}(\tau,z)\xi_{0,6}(\tau,z)=
\frac{\vartheta(\tau, 4z)\vartheta(\tau, 3z)}
{\vartheta(\tau, 2z)\vartheta(\tau, z)}\\
{}&=(r^2+r+2+r^{-1}+r^{-2})+\dots+q^2(r^9+\dots).
\tag{4.16}
\endalign
$$
The Jacobi forms $\eta^2\phi_{0,4}$ and
$\eta\xi_{0,6}$ are holomorphic at infinity. Thus
in \thetag{4.16} we have all types of Fourier coefficients of
$\phi_{0,10}(\tau,z)$ of negative norm.
One can prove that the exponential lifting of $\phi_{0,10}$
coincides with the modular form of weight $1$ for $\Gamma_{10}$
constructed in Example 1.20
$$
F_1^{(10)}(Z)=
\hbox{Lift\,}(\vartheta(\tau,z)\vartheta(\tau,2z))
=
\hbox{Exp-Lift}(\phi_{0,10})
\in \frak N_{1}(\Gamma_{10}, v_\eta^{6}\times v_H)
\tag{4.17}
$$
and $
\hbox{Div}_{\Cal A_{10}^+}(F_1^{(10)})
=2H_1(1)+H_1(9)+H_4
$.
\endexample

\head
5. Some applications
\endhead

In this section we shortly describe some applications of automorphic
forms constructed in \S~1--\S~4. We hope to give more detailed
and advanced applications in further papers.

\subhead
5.1. Denominator identities for automorphic Lorentzian Kac--Moody
algebras
\endsubhead

Almost all identities between arithmetic and exponential liftings
proved above  give denominator formulae for
Lorentzian Kac--Moody algebras.
This means that they give Fourier expansions (or the
corresponding $q$-expansions) of
{\it Lie type}. See Sect. 2.5 of Part I and considerations below.

We denote by $M_t=U(4t)\oplus \langle 2 \rangle$ ($t \in \bn /4$) a
hyperbolic lattice with
the bases $f_2, \widehat{f}_3, f_{-2}$ and with the Gram matrix
$$
\left(\matrix
0&0&-4t\\
0&2&0\\
-4t&0&0
\endmatrix
\right).
$$
We denote an element $nf_2+l\widehat{f}_3+mf_{-2}$ of this lattice
(or of $M_t \otimes \bq$) by its coordinates $(n,l,m)$.

All our formulae have the form:
$$
\Phi=\sum_{(n,l,m) \in M_t}{N(n,l,m)q^{\rho_1+n}r^{\rho_2+l}s^{t(\rho_3+m)}}=
$$
$$
=\sum_{w \in W}{\epsilon (w) \hskip-5pt \sum_{a \in \br_{++}{\Cal M}
\cap M_t}
(e^{w(\rho)}+N(a)e^{w(\rho+a)})}
=q^{\rho_1}r^{\rho_2}s^{t\rho_3}\prod_{(n,l,m)>0}{(1-q^nr^ls^{tm})}^{f(n,l,m)}
\tag{5.1}
$$
for some $W$, $\M$, $\rho=(\rho_1, \rho_2, \rho_3)$ which we define
below. The ``exponent''
$e^{(n,l,m)}:=q^nr^ls^{tm}$. The {\it Weyl group} $W$ is a
reflection subgroup $W\subset W(M_t)\subset O(M_t)$, and
$\M_t$ is a fundamental polyhedron of $W$ in the hyperbolic space
$\La (M_t)$ defined by the hyperbolic lattice  $M_t$.
The $W$ and $\M$ are defined by an
{\it acceptable
set} $P(\M)\subset M_t$ of elements orthogonal to faces of $\M$ and
directed outward. The set $P(\M )$
is also called the {\it set of simple real roots}. The main property of
this set is that $(\alpha, \alpha )>0$ and
$(\alpha , \alpha )\,|\,2(M_t, \alpha)$ for any $\alpha \in P(\M)$.
Moreover,
it defines a {\it generalized Cartan matrix}
$$
A=\left( {2(\alpha, \alpha^\prime )\over (\alpha, \alpha)}\right), \ \
\alpha, \alpha^\prime \in P(\M),
$$
which means that all diagonal elements of $A$
are equal to $2$ and all non-diagonal elements are non-positive
integers. The generalized Cartan matrix $A$ is the main invariant of
the Fourier expansion of Lie type.
The Weyl group $W$ is generated by the reflections
$s_\alpha\in O^+(M_t)$, $\alpha \in P(\M)$. We recall that
$s_\alpha (\alpha )=-\alpha$ and
$s_\alpha$ is identical on the $\alpha^\perp$. The
$\epsilon : W\to \{\pm 1\}$ is a character. It is defined by the set
$P(\M)_\0o=\{ \alpha \in P(\M)\,|\,\epsilon (s_\delta)=-1\}$ which is
called {\it the set of even simple real roots}. Respectively, the set
$P(\M )_\1o=P(\M )\setminus P(\M )_\0o=
\{ \alpha \in P(\M)\,|\,\epsilon (s_\delta)=1\}$
is called the set of {\it odd simple real roots}. By definition,
$a=(n,l,m)\in \br_{++}\M$ if $(a, P(\M))\le 0$ and $a \not=0$.
The  {\it lattice Weyl vector}
$\rho=(\rho_1, \rho_2, \rho_3) \in \br_{++}\M\cap (M_t\otimes \bq )$
is defined by the property: $(\rho, \alpha)=-(\alpha, \alpha)/2$ for
any $\alpha \in P(\M )$. It
is true that $(\rho, \rho)\le 0$. The case $(\rho, \rho)<0$ is called
{\it elliptic}. The case $(\rho, \rho)=0$ is called
{\it parabolic}. All {\it Fourier coefficients}
$N(n,l,m)$, $N(a)$ and {\it multiplicities} $f(n,l,m)$
are integers. The inequality $a=(n,l,m)>0$ means that $a \in M_t$ and
either $(a, \rho)<0$ or  $a \in \bq_{++}\rho$ (if $(\rho, \rho)=0)$.
If the multiplicity
$f(n,l,m)\not=0$, then either $a \in W(P(\M)_\0o)$ and $f(n,l,m)=1$
(i.e. $a$ is a {\it positive even real root}),
or $a \in W(P(\M))_\1o$ and $f(n,l,m)=-1$ (i.e. $a$ is
a {\it positive odd real root}), or $a \in 2W(P(\M ))_\1o$ (i.e.
$a$ is a {\it positive
even real root which is multiple to an odd real  root})
and $f(n,l,m)=1$, or  $(a,a)\le 0$ (i.e. $a$ is a {\it positive
imaginary root}).

Below we describe sets $P(\M)$, $P(\M)_\0o$ and the generalized
Cartan matrices $A$ for all our formulae. To describe some of them
we also use the {\it group of symmetries of $P(\M)$} which is
the group
$$
\Sym (P(\M))=\{g \in O^+(M_t)\,|\, g(P(\M))=P(\M),\
 g(P(\M))_\1o=P(\M)_\1o \}.
$$
We also use notation $W^{(k_1,...,k_r)}(K)$ for the reflection subgroup
of a lattice $K$ generated by reflections in all primitive
elements $\delta \in K$ such that
$(\delta, \delta) \in \{k_1,...,k_n\}\subset \bn$.
The group $W(K)$ denote the full reflection group of the lattice $K$
which is generated by reflections in all elements $\delta \in K$ with
positive squares $(\delta, \delta )>0$.
For all cases the semi-direct product
$W\rtimes \Sym(P(\M))$ has finite index in
$O^+(M_t)$ (i.e. $W$ has {\it restricted arithmetic type})
and for  almost all of them $W\rtimes \Sym(P(\M))$
is equal to the full group
$O^+(M_t)$. Here $O^+(M_t)$ denote the subgroup of $O(M_t)$
which fixes the light cone of $M_t$.

\subsubhead
5.1.1. Cases $\Phi=\Phi_{t,j,\mu}$ and $\Phi=\widetilde{\Phi}_{t,I,\mu}$
where $t=1,\,2,\,3,\,4$, $j=0,\,I,\,II$ and $\mu=\0o,\, \1o$
\endsubsubhead
For these cases

\noindent
$\Phi_{1,II,\0o}=\Delta_5$ (formula \thetag{2.16}
and \cite{GN1}, \cite{GN2});

\noindent
$\Phi_{2,II,\0o}=\Delta_2$ (formula \thetag{2.21} and \cite{GN1});

\noindent
$\Phi_{3,II,\0o}=\Delta_1\ $ (formula \thetag{2.20});

\noindent
$\Phi_{4,II,\0o}=\Delta_{1/2}\ $ (formula \thetag{2.11});

\noindent
$\Phi_{1,0,\0o}=\Delta_{35}\ $ (formula \thetag{3.31} and \cite{GN4});

\noindent
$\Phi_{2,0,\0o}=\Delta_{11}\ $ (formula \thetag{3.11});

\noindent
$\Phi_{3,0,\0o}=D_6\cdot\Delta_1\ $ (formula \thetag{3.32});

\noindent
$\Phi_{4,0,\0o}=\Delta_5^{(4)}\ $ (formula \thetag{3.12});

\noindent
$\Phi_{t,I,\0o}=\Phi_{t,0,\0o}/\Phi_{t,II,\0o},\ $
$\ t=1,\,2,\,3,\,4$ (see \cite{GN4});

\noindent
$\widetilde{\Phi}_{t,I,\0o}(Z)=\Phi_{t,0,\0o}(Z)/\Phi_{t,II,\0o}(2Z),\ $
$\ t=1,\,2,\,3,\,4$ (see \cite{GN4});

\noindent
$\Phi_{t,\1o}=1$ for $t=1$ and $\Phi_{t,\1o}=\Psi^{(t)}_{12}\ $ for
$t=2,\,3,\,4$
(formulae \thetag{4.7}, \thetag{4.8}, \thetag{4.9});

\noindent
$\Phi_{t,j,\1o}=\Phi_{t,j,\0o}\cdot\Phi_{t,\1o},\ $ $\ t=1,\,2,\,3,\,4$;

\noindent
$\widetilde{\Phi}_{t,I,\1o}=\widetilde{\Phi}_{t,I,\0o}\cdot\Phi_{t,\1o},\ $
$\ t=1,\,2,\,3,\,4$.

\smallpagebreak

According to our definitions,
$\Phi_{1,j,\1o}=\Phi_{1,j,\0o}$ and
$\widetilde{\Phi}_{1,I,\1o}=\widetilde{\Phi}_{1,I,\0o}$.
We denote the corresponding to these forms $A$, $W$, $\M$  and
$\rho$ by the same indexes and ``tilde''.
We describe  them below.
For all of them we have
$$
W_{t,j,\mu}\rtimes \Sym (P_{t,j,\mu}(\M_{t,j,\mu}))=O^+(M_t),\ \
\widetilde{W}_{t,I,\mu}\rtimes \Sym (\widetilde{P}_{t,I,\mu}
(\widetilde{\M}_{t,I,\mu}))=O^+(M_t).
$$

We have for $(t,I,\1o)$ and $(t,0,\1o)$:
$$
W_{t,I,\1o}=W_{t,0,\1o}=W(M_t)=W^{(2,8,8t)}(M_t)
$$
with the same fundamental polyhedron
$\M_{t,I,\1o}=\M_{t,0,\1o}$.

\smallpagebreak

{\bf Cases $(t,I,\1o)$, $t=1,\,2,\,3,\,4$.} We have
$$
P_{t,I,\1o}(\M_{t,I,\1o})=\{\delta_1=(0,-1,0),\,
\delta_2=(1,2,0),\,  \delta_3=(-1,0,1)\},
\ \, \quad \, P_{t,I,\1o}(\M_{t,I,\1o})_\1o=\{\delta_1\}
$$
with the Gram matrix
$$
G(P_{t,I,\1o}(\M_{t,I,\1o}))=(\delta_i,\,\delta_j)=
\left(\matrix
\hphantom{-}2&-4&\hphantom{-}0\\
-4&\hphantom{-}8&-4t\\
\hphantom{-}0&-4t&\hphantom{-}8t
\endmatrix
\right)
$$
and the generalized Cartan matrix
$$
A_{t,I,\1o}=
\left(
\matrix
\hphantom{-}2 &-4 &\hphantom{-}0\\
-1& \hphantom{-}2 &-t\\
\hphantom{-}0 &-1 & \hphantom{-}2
\endmatrix
\right)
$$
with the set of odd indices $\{1\}$. We have
$$
\rho_{t,I,\1o}=\left( {2t+3\over 2t},\ {1\over 2},\ {3\over 2t}\right).
$$
The matrix $A_{1,I,\1o}=A_{1,I,\0o}$ is equal
to the matrix $r=-59/2$ of the Table 1 in Part I and is twisted
to the symmetric generalized Cartan matrix $A_{1,0}$ of Theorem
1.3.1 in Part I.

\smallpagebreak

{\bf Cases $(t,0,\1o)$, $\ t=1,\,2,\,3,\,4$.} We have:
$$
P_{t,0,\1o}(\M_{t,0,\1o})=\{2\delta_1=(0,-2,0),\
\delta_2=(1,2,0),\  \delta_3=(-1,0,1)\}
\ \ \text{and}\ \  P_{t,0,\1o}(\M_{t,0,\1o})_\1o=\emptyset
$$
with the Gram matrix
$$
G(P_{t,0,\1o}(\M_{t,0,\1o}))=
\left(\matrix
\hphantom{-}8&-8&\hphantom{-}0\\
-8&\hphantom{-}8&-4t\\
\hphantom{-}0&-4t&\hphantom{-}8t
\endmatrix
\right)
$$
and the generalized Cartan matrix
$$
A_{t,0,\1o}=
\left(
\matrix
\hphantom{-}2 &-2 &\hphantom{-}0\\
-2& \hphantom{-}2 &-t\\
\hphantom{-}0 &-1 & \hphantom{-}2
\endmatrix
\right)
$$
with the empty set of odd indices. We have
$$
\rho_{t,0,\1o}=\left({t+2\over t},\ 1,\ {2\over t}\right).
$$
The matrix $A_{1,0,\1o}=A_{1,0,\0o}$ is equal
to the symmetric generalized Cartan matrix $A_{1,0}$ of Theorem
1.3.1 in Part I.

\smallpagebreak

{\bf Cases $\widetilde{(t,I,\1o)}$, $\ t=1,\,2,\,3,\,4$.} We have
$$
\widetilde{W}_{t,I,\1o}=W^{(8,8t)}(M_t)
$$
with the fundamental polyhedron
$$
\widetilde{\M}_{t,I,\1o}=\Sym(\widetilde{P}_{t,I,\1o}
(\widetilde{\M}_{t,I,\1o}))(\M_{t,I,\1o})
$$
where
$\Sym(\widetilde{P}_{t,I,\1o}(\widetilde{\M}_{t,I,\1o}))$
is generated by the reflection
in $\delta_1=(0,-1,0)$. Respectively,
$$
\align
\widetilde{P}_{t,I,\1o}(\widetilde{\M}_{t,I,\1o})&=
\Sym(\widetilde{P}_{t,I,\1o}(\widetilde{\M}_{t,I,\1o}))
(\delta_2, \delta_3)\\
{}&=\{\delta_2=(1,2,0),\,  \delta_3=(-1,0,1),\,  \delta_2^\prime =(1,-2,0)\},
\endalign
$$
$$
\widetilde{P}_{t,I,\1o}(\widetilde{\M}_{t,I,\1o})_\1o=\emptyset.
$$
We have
$$
G(\widetilde{P}_{t,I,\1o}(\widetilde{\M}_{t,I,\1o}))=
\left(
\matrix
\hphantom{-}8 & -4t & -8\\
-4t& \hphantom{-}8t &-4t\\
-8 & -4t& \hphantom{-}8
\endmatrix
\right),
$$
the generalized Cartan matrix
$$
\widetilde{A}_{t,I,\1o}=
\left(
\matrix
\hphantom{-}2 & -t &-2\\
-1&  \hphantom{-}2 &-1\\
-2& -t & \hphantom{-}2
\endmatrix
\right),
$$
with the empty set of odd indices. The lattice Weyl vector
$$
\widetilde{\rho }_{t,I,\1o}=({t+1\over t},\, 0,\, {1\over t}).
$$
The matrix $\widetilde{A}_{1,I,\1o}=\widetilde{A}_{1,I,\0o}$
gives the symmetric generalized Cartan matrix
$A_{1,I}$ of Theorem 1.3.1 in Part I.
The matrix $\widetilde{A}_{4,I,\1o}$ gives the generalized
Cartan matrix corresponding to $r=-5/2$ of Table 1 in Part I.
This matrix is twisted to the matrix $A_{1,II}$.

\smallpagebreak

{\bf Cases $(t,II,\1o)$, $t=2,\,3,\,4$.} The case $(1,II,\1o)$
coincides with $(1,II,\0o)$ and will be considered later together
with all cases $(t,II,\0o)$.
We have
$$
W_{t,II,\1o}=W^{(2,8t)}(M_t),\qquad t=2,\,3,\,4,
$$
with the fundamental polyhedron
$$
\M_{t,II,\1o}=\Sym(P_{t,II,\1o}(\M_{t,II,\1o}))(\M_{t,0,\1o}),
\qquad t=2,\,3,\,4,
$$
where
$\Sym (P_{t,II,\1o}(\M_{t,II,\1o}))$ is generated by the reflection
in $\delta_2=(1,2,0)$ and
$$
P_{t,II,\1o}(\M_{t,II,\1o})
=\Sym (P_{t,II,\1o}(\M_{t,II,\1o}))(2\delta_1,
\delta_3)=
$$
$$
=\{2\delta_1=(0,-2,0),\,\delta_3=(-1,0,1),\,\delta_3^\prime=(t-1,2t,1),
\,2\delta_1^\prime=(2,2,0)\}
$$
with the Gram matrix
$$
G(P_{t,II,\1o}(\M_{t,II,\1o}))=
\left(
\matrix
\hphantom{-}8 & \hphantom{-}0 &-8t &-8\\
\hphantom{-}0 & \hphantom{-}8t&-4t^2+8t&-8t\\
-8t &-4t^2 + 8t &\hphantom{-}8t &\hphantom{-}0\\
-8 &-8t &\hphantom{-}0 &\hphantom{-}8
\endmatrix
\right)
$$
and the generalized Cartan matrix
$$
A_{t,II,\1o}=
\left(
\matrix
\hphantom{-}2 &\hphantom{-}0 &-2t &-2\\
\hphantom{-}0 &\hphantom{-}2 &-t + 2 &-2\\
-2 &-t + 2 &\hphantom{-}2 &\hphantom{-}0\\
-2 &-2t &\hphantom{-}0 &\hphantom{-}2
\endmatrix
\right)
$$
with the empty set of odd indices. The lattice Weyl vector
$$
\rho_{t,II,\1o}=({t+1\over t},\ 1,\ {1\over t}).
$$
The matrix $A_{4,II,\1o}$ gives the third generalized Cartan matrix
for $r=-2$ of Table 1 in Part I. It is twisted to the generalized
Cartan matrix $A_{2,I}$ of Theorem 1.3.1 in Part I.

\smallpagebreak

{\bf Cases $(t,\1o)$, $\ t=2,\,3,\,4$.} For these cases
$$
W_{t,\1o}=W^{(8t)}(M_t),\qquad t=2,\,3,\,4,
$$
with the fundamental polyhedron
$$
\M_{t,\1o}=\Sym(P_{t,\1o}(\M_{t,\1o}))(\M_{t,0,\1o}),
\qquad t=2,\,3,\,4,
$$
where (we consider only $t=2,\,3,\,4$) the group
$\Sym (P_{t,\1o}(\M_{t,\1o}))$ is generated by the reflections in
$\delta_1=(0,-1,0)$ and $\delta_2=(1,2,0)$. We have
$$
\rho_{t,\1o}=(1,0,0)
$$
and for $\delta_3=(-1,0,1)$
$$
P_{t,\1o}(\M_{t,\1o})=P_{t,\1o}(\M_{t,\1o})_\0o=
\Sym (P_{t,\1o}(\M_{t,\1o}))(\delta_3)=
$$
$$
=\{\text{primitive\ } \delta \in M_t\ |\
(\delta,\delta)=8t,\ 4t|(\delta,M_t),\
(\delta,\rho)=-4t\}.
$$
These case is parabolic, the group $\Sym (P_{t,\1o}(\M_{t,\1o}))$
is the affine reflection group on a line (it is isomorphic to $\bz$ up to
index two). The Gram matrix is
$$
G(P_{t,\1o}(\M_{t,\1o}))=(\alpha, \alpha^\prime),\ \ \alpha,\,\alpha^\prime
\in P_{t,\1o}(\M_{t,\1o}),
$$
and the generalized Cartan matrix is
$$
A_{t,\1o}=\left(
{(\alpha, \alpha^\prime)\over 4t}
\right),\ \ \alpha,\,\alpha^\prime
\in P_{t,\1o}(\M_{t,\1o}),
$$
with the empty set of odd indices. It is symmetric of parabolic type
and has a lattice Weyl vector. It is interesting that for $t=2,\,3$
this matrix does not have a parabolic submatrix of
the type $\widetilde{{\Bbb A}_1}$.

\smallpagebreak

{\bf Case $(t,0,\0o)$, $\ t=1,\,2,\,3,\,4$.} For this case
$$
W_{t,0,\0o}=W^{(2,8)}(M_t).
$$
and
$$
\M_{t,0,\0o}=\Sym (P_{t,0,\0o}(\M_{t,0,\0o}))(\M_{t,0,\1o}).
$$
The group $\Sym (P_{1,0,\0o}(\M_{1,0,\0o}))$ is trivial and
the case $(1,0,\0o)$
coincides with the $(1,0,\1o)$ which has been considered before.
Below we assume that $t=2,3,4$. We have
$$
P_{t,0,\0o}(\M_{t,0,\0o})=P_{t,0,\0o}(\M_{t,0,\0o})_\0o
=\Sym (P_{t,0,\0o}(\M_{t,0,\0o}))(2\delta_1,\delta_2)=
$$
$$
=\{2\delta_1=(0,-2,0),\,\delta_2=(1,2,0),\,\delta_2^\prime=(0,2,1)\}
$$
where $\Sym (P_{t,0,\0o}(\M_{t,0,\0o}))$ is generated by the
reflection in $\delta_3=(-1,0,1)$. We have
$$
\rho_{t,0,\0o}=\left({2\over t},\ 1,\ {2\over t}\right),
$$
$$
G(P_{t,0,\0o}(\M_{t,0,\0o})))=
\left(
\matrix
\hphantom{-}8 & -8 & -8\\
-8&  \hphantom{-}8 &-4t+8\\
-8&-4t+8&   \hphantom{-}8
\endmatrix
\right)
$$
and
$$
A_{t,0,\0o}=
\left(
\matrix
\hphantom{-}2 & -2 & -2\\
-2& \hphantom{-} 2 &-t+2\\
-2&-t+2&   \hphantom{-}2
\endmatrix
\right)
$$
with the empty set of odd indices. The matrices $A_{t,0,\0o}$, $t=1,2,3$,
coincide with the symmetric generalized Cartan matrices $A_{t,0}$ of Theorem
1.3.1 of Part I. The matrix $A_{4,0,\0o}$ coincides with the matrix
$A_{1,II}$ of this theorem.

\smallpagebreak

{\bf Case $(t,I,\0o)$, $t=1,2,3,4$.}
For this case
$$
W_{t,I,\0o}=W^{(2,8)}(M_t)
$$
and
$$
\M_{t,I,\0o}=\Sym (P_{t,I,\0o}(\M_{t,I,\0o}))(\M_{t,I,\1o}).
$$
The group $\Sym (P_{1,I,\0o}(\M_{1,I,\0o}))$ is trivial and
the case $(1,I,\0o)$
coincides with the $(1,I,\1o)$ which has been considered before.
Below we assume that $t=2,\,3,\,4$. We have
$$
P_{t,I,\0o}(\M_{t,I,\0o})
=\Sym (P_{t,I,\0o}(\M_{t,I,\0o}))(\delta_1,\delta_2)=
$$
$$
=\{\delta_1=(0,-1,0),\,\delta_2=(1,2,0),\,\delta_2^\prime=(0,2,1)\}
$$
and
$$
P_{t,I,\0o}(\M_{t,I,\0o})_\1o=\{\delta_1=(0,-1,0)\}
$$
where $\Sym (P_{t,I,\0o}(\M_{t,I,\0o}))$ is generated by the
reflection in $\delta_3=(-1,0,1)$. We have
$$
\rho_{t,I,\0o}=\left({3\over 2t},\ {1\over2},\ {3\over 2t}\right),
$$
$$
G(P_{t,I,\0o}(\M_{t,I,\0o})))=
\left(
\matrix
{\hphantom{-}{2}}&{{-4}}&{{-4}}\\
{{-4}}&{\hphantom{-}{8}}&{{-4}t + {8}}\\
{{-4}}&{{-4}t + {8}}&{\hphantom{-}{8}}
\endmatrix
\right)
$$
and
$$
A_{t,I,\0o}=
\left(
\matrix
{\hphantom{-}{2}}&{{-4}}&{{-4}}\\
{-1}&{\hphantom{-}{2}}&{-t + {2}}\\
{-1}&{-t + {2}}&{\hphantom{-}{2}}
\endmatrix
\right)
$$
with the set $\{1\}$ of odd indices.
In the Table 1 of Part I, the matrix $A_{2,I,\0o}$ corresponds to
$r=-17/2$, the matrix $A_{3,I,\0o}$ corresponds to the second matrix with
$r=-11/2$,
the matrix $A_{4,I,\0o}$ corresponds to the third matrix with
$r=-4$. These matrices
are twisted to the matrices $A_{2,0}$, $A_{3,0}$, $A_{1,II}$ of Theorem
1.3.1 in Part I respectively.

\smallpagebreak

{\bf Cases $\widetilde{(t,I,\0o)}$, $\ t=1,\,2,\,3,\,4$.}
For these cases
$$
\widetilde{W}_{t,I,\0o}=W^{(8)}(M_t)
$$
with the fundamental polyhedron
$$
\widetilde{\M}_{t,I,\0o}=\Sym(\widetilde{P}_{t,I,\0o}
(\widetilde{\M}_{t,I,\0o}))(\M_{t,I,\1o}).
$$
The case $\widetilde{(1,I,\0o)}$ coincides with $\widetilde{(1,I,\1o)}$
and has been considered. Below we suppose that $t=2,3,4$.
For $t=2,3,4$, the group
$\Sym(\widetilde{P}_{t,I,\0o}(\widetilde{\M}_{t,I,\0o}))$
is generated by the reflections
in $\delta_1=(0,-1,0)$ and $\delta_3=(-1,0,1)$, and is isomorphic to
$D_2$. We have
$$
\widetilde{P}_{t,I,\0o}(\widetilde{\M}_{t,I,\0o})=
\widetilde{P}_{t,I,\0o}(\widetilde{\M}_{t,I,\0o})_\0o=
\Sym(\widetilde{P}_{t,I,\0o}(\widetilde{\M}_{t,I,\0o}))(\delta_2)=
$$
$$
=\{\delta_{21}=\delta_2=(1,2,0),\, \delta_{22}=(0,2,1),\,
\delta_{23}=(0,-2,1), \,\delta_{24}=(1,-2,0)\},
$$
with the Gram matrix
$$
G(\widetilde{P}_{t,I,\0o}(\widetilde{\M}_{t,I,\0o}))=
\left(
\matrix
{\hphantom{-}{8}}&{{-4}t + {8}}&{{-4}t - {8}}&{{-8}}\\
{{-4}t + {8}}&{\hphantom{-}{8}}&{{-8}}&{{-4}t - {8}}\\
{{-4}t - {8}}&{{-8}}&{\hphantom{-}{8}}&{{-4}t + {8}}\\
{{-8}}&{{-4}t - {8}}&{{-4}t + {8}}&{\hphantom{-}{8}}
\endmatrix
\right),
$$
and the generalized Cartan matrix
$$
\widetilde{A}_{t,I,\0o}=
\left(
\matrix
{\hphantom{-}{2}}&{-t + {2}}&{-t - {2}}&{{-2}}\\
{-t + {2}}&{\hphantom{-}{2}}&{{-2}}&{-t - {2}}\\
{-t - {2}}&{{-2}}&{\hphantom{-}{2}}&{-t + {2}}\\
{{-2}}&{-t - {2}}&{-t + {2}}&{\hphantom{-}{2}}
\endmatrix
\right),
$$
with the empty set of odd indices. The lattice Weyl vector
$$
\widetilde{\rho }_{t,I,\0o}=({1\over t},\ 0,\ {1\over t}).
$$
For $t=1,\,2,\,3$, the generalized Cartan matrices
$\widetilde{A}_{t,I,\0o}$
give the symmetric generalized Cartan matrices
$A_{t,I}$, and $\widetilde{A}_{4,I,\0o}$ gives the symmetric
generalized Cartan matrix $A_{2,II}$ of Theorem 1.3.1 in Part I.

\smallpagebreak

{\bf Cases $(t,II,\0o)$, $\ t=1,\,2,\,3,\,4$.}
For these cases
$$
W_{t,II,\0o}=W^{(2)}(M_t)
$$
with
$$
P_{t,II,\0o}(\M_{t,II,\0o})
=\Sym (P_{t,II,\0o}(\M_{t,II,\0o}))(\M_{t,0,\1o})
$$
where the group $\Sym (P_{t,II,\0o}(\M_{t,II,\0o}))$ is generated by
the reflections in $\delta_2=(1,2,0)$ and $\delta_3=(-1,0,1)$.
We can take
$$
\rho_{t,II,\0o}=\left({1\over 2t},\ {1\over 2},\ {1\over 2t} \right).
$$
Then
$$
P_{t,II,\0o}(\M_{t,II,\0o})=P_{t,II,\0o}(\M_{t,II,\0o})_\0o=
\Sym (P_{t,II,\0o}(\M_{t,II,\0o}))(\delta_1),\ \
\delta_1=(0,-1,0).
$$
We have
$$
P_{1,II,\0o}(\M_{1,II,\0o})=\{\delta_{11}=(0,-1,0),\,
\delta_{12}=(1,1,0),\,\delta_{13}=(0,1,1)\},
$$
$$
P_{2,II,\0o}(\M_{2,II,\0o})=\{\delta_{11}=(0,-1,0),\,\delta_{12}=(1,1,0),
\,\delta_{13}=(1,3,1),\, \delta_{14}=(0,1,1)\},
$$
$$
\split
P_{3,II,\0o}(\M_{3,II,\0o})=\{&\delta_{11}=(0,-1,0),\, \delta_{12}=(1,1,0),\\
&\delta_{13}=(2,5,1),\,
\delta_{14}=(2,7,2),\,\delta_{15}=(1,5,2),\,\delta_{16}=(0,1,1)\},
\endsplit
$$
and
$$
\split
P_{4,II,\0o}(\M_{4,II,\0o})&=
\Sym(P_{4,II,\0o}(\M_{4,II,\0o}))((0,-1,0))\\
&=\{\delta\in M_4 |\ (\delta , \delta)=2,\
(\delta, \rho)=-1\}.
\endsplit
$$
It follows that $G(P_{t,II,\0o}(\M_{t,II,\0o}))=A_{t,II,\0o}$ where
for $t=1,\,2,\,3$ the generalized Cartan matrix
$A_{t,II,\0o}$ is equal to the generalized Cartan matrix $A_{t,II}$ of
Theorem 1.3.1 in Part I, and the
$A_{4,II,\0o}$ is the symmetric parabolic generalized Cartan matrix
(it is infinite)
$$
A_{4,II,\0o}=\bigl((\alpha, \alpha^\prime)\bigr),\ \
\alpha, \alpha^\prime \in P_{4,II,\0o}(\M_{4,II,\0o}).
$$

\subsubhead
5.1.2. The case $D_2(Z)$
\endsubsubhead
The automorphic form $D_2(Z)$ was considered in Theorem 4.1.
For this case $t=9$ and the Weyl group is
$$
W=W^{(2,8,18)}(M_9).
$$
The full group $W(M_9)=W^{(2,8,18,72)}(M_9)$ has the fundamental
polyhedron $\M_0$ with the set of primitive orthogonal vectors
$$
P(\M_0)=\{\alpha_1=(0,-1,0),\, \alpha_2=(1,2,0), \,\alpha_3=(2,9,1),
\,\delta=(-1,0,1)\}
$$
with the Gram matrix
$$
G(P(\M_0))=
\left(
\matrix
\hphantom{-}{2}&{-4}&{-18}&\hphantom{-}{0}\\
{-4}&\hphantom{-}{8}&\hphantom{-}{0}&{-36}\\
{-18}&\hphantom{-}{0}&\hphantom{-}{18}&{-36}\\
\hphantom{-}{0}&{-36}&{-36}&\hphantom{-}{72}
\endmatrix
\right).
$$
We have
$$
\rho=\left({1\over 6},\ {1\over 2},\ {1\over 6}\right),
$$
The group $\Sym(P(\M))$ is generated by the reflection in $\delta=(-1,0,1)$,
and
$$
\split
&P(\M)=\Sym(P(\M))(\alpha_1, \alpha_2, \alpha_3)=
\{(0,-1,0),\,(1,2,0),\,(2,9,1),\,(1,9,2),\,(0,2,1)\},\\
&P(\M)_\1o=\{ (0,-1,0)\}.
\endsplit
$$
The Gram matrix
$$
G(P(\M))=
\left(
\matrix
\hphantom{-}{2}&{-4}&{-18}&{-18}&{-4}\\
{-4}&\hphantom{-}{8}&\hphantom{-}{0}&{-36}&{-28}\\
{-18}&\hphantom{-}{0}&\hphantom{-}{18}&{-18}&{-36}\\
{-18}&{-36}&{-18}&\hphantom{-}{18}&\hphantom{-}{0}\\
{-4}&{-28}&{-36}&\hphantom{-}{0}&\hphantom{-}{8}
\endmatrix
\right)
$$
and the generalized Cartan matrix
$$
A(P(\M))=
\left(\matrix
\hphantom{-}{2}&{-4}&{-18}&{-18}&{-4}\\
{-1}&\hphantom{-}{2}&\hphantom{-}{0}&{-9}&{-7}\\
{-2}&\hphantom{-}{0}&\hphantom{-}{2}&{-2}&{-4}\\
{-2}&{-4}&{-2}&\hphantom{-}{2}&\hphantom{-}{0}\\
{-1}&{-7}&{-9}&\hphantom{-}{0}&\hphantom{-}{2}
\endmatrix\right)
$$
with the set $\{1\}$ of odd indices. This is the generalized Cartan
matrix of the case $r=-3/2$ of Table 1 in Part I. It is
twisted to the symmetric generalized Cartan
matrix $A_{1,III}$ of Theorem 1.3.1 in Part I.

\subsubhead
5.1.3. The case $D_{1/2}(Z)$
\endsubsubhead
The final formula for the automorphic form $D_{1/2}(Z)$  is given
in \thetag{2.14}.
For this case $t=36$ and the Weyl group is
$$
W=W^{(2,8,18,32)}(M_{36}).
$$
The group $W(M_{36})=W^{(2,8,18,32,72,288)}(M_{36})$ has the
fundamental polyhedron $\M_0$ with the set of primitive orthogonal
vectors
$$
\split
P(\M_0)=\{\alpha_1=(0,-1,0),\, \alpha_2=(1,2,0),\,&
\alpha_3=(7,32,1), \,
\alpha_4 =(5,27,1),\\
&\delta_1=(2,18,1),\, \delta_2=(-1,0,1)\}
\endsplit
$$
with the Gram matrix
$$
G(P(\M_0))=
\left(\matrix
\hphantom{-}{2}&{-4}&{-64}&{-54}&{-36}&\hphantom{-}{0}\\
{-4}&\hphantom{-}{8}&{-16}&{-36}&{-72}&{-144}\\
{-64}&{-16}&\hphantom{-}{32}&\hphantom{-}{0}&{-144}&{-864}\\
{-54}&{-36}&\hphantom{-}{0}&\hphantom{-}{18}&{-36}&{-576}\\
{-36}&{-72}&{-144}&{-36}&\hphantom{-}{72}&{-144}\\
\hphantom{-}{0}&{-144}&{-864}&{-576}&{-144}&\hphantom{-}{288}
\endmatrix\right).
$$
The fundamental polyhedron $\M$ for the Weyl group $W$ is
$$
\split
\M=\Sym (P(\M))(\M_0),\ \  &P(\M)=\Sym(P(\M))(\alpha_1,\alpha_2,
\alpha_3,\alpha_4),\\
&P(\M)_\1o=\Sym(P(\M))(\alpha_1).
\endsplit
$$
The lattice Weyl vector is
$$
\rho=\left({1\over 24},\ {1\over 2},\ {1\over 24}\right).
$$
The symmetry group $\Sym (P(\M))$ is generated by the reflections in
$\delta_1=(2,18,1)$ and $\delta_2=(-1,0,1)$. It is infinite, and
this case is parabolic with the parabolic generalized Cartan matrix
$$
A(P(\M))=\left({2(\alpha,\alpha^\prime)\over (\alpha,\alpha)}
\right),\ \ \alpha,\alpha^\prime \in P(\M).
$$

\subhead
5.2. Discriminants of moduli of K3 surfaces and Mirror Symmetry
\endsubhead
\smallskip
A lattice $T$ with signature $(n,2)$ is called $2$-reflective if
there exists an automorphic form $\Phi$ on $\Omega (T)$
with respect to a subgroup of $O(T)$ of finite index such that
the divisor of $\Phi$ is union of quadratic divisors orthogonal to
elements of $T$ with norm $2$. The automorphic form $\Phi$ is
called {\it $2$-reflective} (for $T$).
Obviously, any overlattice $T\subset T^\prime$ of
finite index of a $2$-reflective lattice $T$ is also
$2$-reflective. Thus, one $2$-reflective lattice $T$
gives many other $2$-reflective lattices. A $2$-reflective lattice $T$
is called {\it strongly $2$-reflective} if there exists a $2$-reflective
automorphic form $\Phi$ for $T$ such that the divisor of $\Phi$ is union
of {\it all} quadratic divisors orthogonal to elements of
$T$ with norm $2$.
The automorphic form $\Phi$ is then called {\it strongly $2$-reflective}.
Conjecturally (see \cite{N6} and Part I) number of $2$-reflective
lattices $T$ of rank $\ge 5$ is finite. $2$-reflective lattices and
strongly $2$-reflective lattices are important for moduli of K3-surfaces.
By the global Torelli Theorem (see \cite{P-\u S--\u S})
and epimorphicity of period map (see \cite{Ku}),
moduli $\M_S$ of K3 surfaces with a condition $S\subset L_{K3}$
on Picard lattice have discriminant given by an automorphic form equation
if the lattice $T=S^\perp_{L_{K3}}(-1)$ is strongly $2$-reflective.
Here $S$ is a hyperbolic lattice of signature $(1,k)$, which is
primitively embedded into the unimodular even lattice $L_{K3}$ of
signature $(3,19)$. This equation is given by a strongly reflective
automorphic form $\Phi$ of the lattice $T$ which is then called {\it
the discriminant automorphic form}. The zero divisor of a
reflective automorphic form $\Phi$ gives a part of the discriminant.

We get a list of strongly $2$-reflective lattices $T$ of rank $5$.

\proclaim{Theorem 5.2.1} The lattices $T$,
$$
T=U^2\oplus \langle 2k \rangle,\qquad k=2,\, 3,\, 4,
$$
and
$$
T=U(k)^2\oplus \langle 2 \rangle,\qquad
k=1,\,2,\dots,\,7,\,8,\, 10,\, 12,\, 16,
$$
are strongly $2$-reflective. We have the following list of
strongly $2$-reflective automorphic forms $\Phi$ for these lattices (after
some isomorphism of $\Omega (T)$ with ${\Bbb H}_t^+$ which is
described in the proof below).
$$
\split
&U^2\oplus \langle 2 \rangle\ \ \ \Delta_{35}=\Phi_{1,0,\0o}\ \ \
\text{(see \thetag{3.31} and \cite{GN4})};\\
&U^2\oplus \langle 2t \rangle\ \ \ \Psi^{(t)}_{12}=\Phi_{t,\1o},
\ \ t=2,3,4\ \ \text{(see \thetag{4.7} --- \thetag{4.9})};\\
&U(2)^2\oplus \langle 2 \rangle\ \ \ \Delta_{11}=\Phi_{2,0,\0o}\ \ \
\text{(see \thetag{3.11})};\\
&U(3)^2\oplus \langle 2 \rangle\ \ \ D_6\cdot\Delta_{1}=\Phi_{3,0,\0o}\ \ \
\text{(see \thetag{3.32})};\\
&U(4)^2\oplus \langle 2 \rangle\ \ \ \Delta_{5}=\Phi_{1,II,\0o}\ \ \
\text{(see \thetag{2.16} and \cite{GN1}, \cite{GN2})};\\
&U(5)^2\oplus \langle 2 \rangle\ \ \ F_5^{(5)}\ \ \
\text{(see Example 4.5, \thetag{4.10})};\\
&U(6)^2\oplus \langle 2 \rangle\ \ \ F_3^{(6)}\ \ \
\text{(see Example 4.6 \thetag{4.12})};\\
&U(7)^2\oplus \langle 2 \rangle\ \ \ F_2^{(7)}\ \ \
\text{(see Example 4.7 \thetag{4.15})};\\
&U(8)^2\oplus \langle 2 \rangle\ \ \ \Delta_{2}=\Phi_{2,II,\0o}\ \ \
\text{(see \thetag{3.11} and \cite{GN1})};\\
&U(10)^2\oplus \langle 2 \rangle\ \ \ F_1^{(10)}\ \ \
\text{(see Example 4.8 \thetag{4.17})};\\
&U(12)^2\oplus \langle 2 \rangle\ \ \ \Delta_{1}=\Phi_{3,II,\0o}\ \ \
\text{(see \thetag{2.20})};\\
&U(16)^2\oplus \langle 2 \rangle\ \ \ \Delta_{1/2}=\Phi_{4,II,\0o}\ \ \
\text{(see \thetag{2.11})}.
\endsplit
$$
\endproclaim

\demo{Proof} All automorphic forms $\Phi$ we use in this theorem  are
automorphic with respect to the paramodular group
$\widehat{O}(L_t)$ where $L_t=U^2\oplus \langle 2t \rangle$,
$t \in \bn$ (see Sect. 1.3),
and have the divisors which are unions of quadratic divisors $\Ha_\delta$,
where $\delta$ are some primitive roots of $L_t$ completely
described by the construction of $\Phi$. Here $\delta \in L_t$ is
a root of $L_t$ if the reflection in $\delta$ belongs to $O(L_t)$.

For each automorphic form $\Phi$ of the theorem,
we replace $L_t$ by a canonical (with respect to $L_t$) sublattice
$T\subset (L_t\otimes \bq )(m)$, $m \in \bq$. Here
{\it canonical} means that
$O(T)=O(L_t)$ in $O(L_t\otimes \bq)$. Then a primitive root $\delta$
in $L_t$ is replaced by a primitive root $\tilde{\delta} \in T$,
$\tilde{\delta} \in \bq_{++}\delta$. Obviously, then $\Phi$ will be
automorphic with respect to the corresponding
subgroup of finite index of $O(T)$. Below, for each $\Phi$,
we give a procedure which changes $L_t$ by $T$ in such a way that
any quadratic divisor $\Ha_\delta$, $\delta \in L_t$,
of $\Phi$ will be replaced by $\Ha_{\tilde{\delta}}$, $\tilde{\delta}\in
T$ and $(\tilde{\delta}, \tilde{\delta})=2$ in the lattice $T$.
Moreover, any element $\tilde{\delta}\in T$ with the norm 2
gives a rational quadratic divisor of $\Phi$. Thus, the lattice $T$
will be strongly $2$-reflective with the strongly $2$-reflective
automorphic form $\Phi$.

For $\Delta_{35}$ ($t=1$) and $\Psi_{12}^{(t)}$, $t=2,\,3,\,4$, we put
$T=L_t=U^2\oplus \langle 2t \rangle$.

For $\Delta_5$, $\Delta_2$, $\Delta_1$, $\Delta_{1/2}$ which
correspond to $t=1,2,3,4$ respectively, one should put
$T=(L_t)^\ast(4t)$.
Then $T=U(4t)^2\oplus \langle 2 \rangle$.

For $D_6\cdot\Delta_1$, $F_5^{(5)}$, $F_2^{(7)}$ which
correspond to odd $t=3$, $5$, $7$ respectively, we take
$T\otimes \bz_p=(L_t\otimes \bz_p)^\ast(t)$ for $p\not=2$, and
$T\otimes \bz_2=(L_t\otimes \bz_2)(t)$. Here $\bz_p$ denote
$p$-adic integers.
Then $T=U(t)^2\oplus \langle 2 \rangle$.

For $F^{(6)}_3$ and $F^{(10)}_1$ which correspond to even
$t=6$, $10$ respectively, with odd $t/2$, we first replace $L_t$
by its overlattice of index $2$ (it is unique). We get
an odd lattice $\tilde L_t=U^2\oplus \langle t/2 \rangle$. The
lattice $L_t$ is its maximal even sublattice. Then we consider
$T=(\tilde L_t)^\ast (t)=U(t)^2 \oplus \langle 2 \rangle$.

Using properties of the paramodular orthogonal group
$O(L_t)$ described in Sect. 1.3, one can  check
the properties of $T$ we claimed. It finishes the proof.
\enddemo

The automorphic forms $\Phi$ of Theorem 5.2.1 can be used to
continue series of examples of the variant of
Mirror Symmetry for K3 surfaces we suggested in
\cite{GN3}. We denote the automorphic form $\Phi$ corresponding to
the lattice $T$ from Theorem 5.2.1
by $\Phi_T$. We remark that for the most part of cases
the form $\Phi_T$ has zeros of multiplicity one and is the unique
strongly reflective automorphic form for $T$
with this property (by Koecher principle, see \cite{Ba}, \cite{F1}
for example). All other forms have very small multiplicities of
zeros and are automorphic with respect to natural groups of small
indices in the full orthogonal groups.
It follows that they are very exceptional.

The automorphic forms $\Phi_T$ may be considered as ``algebraic functions'' on
moduli $\M_{T^\perp}$ of K3 surfaces with the condition
$T^\perp:=T(-1)^\perp_{L_{K3}}$ on
the Picard lattice (B-model). For any $T$, the hyperbolic
lattices $T^\perp$ may be easily
computed using the discriminant forms technique (see \cite{N1}) and
give very interesting moduli of K3 surfaces. We hope to
study equations of the corresponding K3 surfaces
in further publications.

Let us consider an isotropic element $c$ in the first summand
$U(k)$ of the lattice $T$. We then get a hyperbolic lattice
$$
S=(c^\perp_T/[c])(-1)=
\cases
U\oplus \langle -2t \rangle &\text{for $T=U^2\oplus \langle 2t \rangle$,}\\
U(k)\oplus \langle -2 \rangle &\text{for $T=U(k)^2\oplus \langle 2 \rangle$,}
\endcases
$$
and can consider moduli $\M_S$ of K3
surfaces with condition $S$ on the Picard lattice (A-model).
A general K3 surface $X$ from this moduli has Picard lattice $S$. The Fourier
expansion of $\Phi_T$ at the cusp $c$ which we used in this paper
is then related with the geometry of
non-singular rational curves on these K3 surfaces $X$ (see \cite{GN3}).

More exactly it means (see details in \cite{GN3} and Part I)
that replacing $M_t$ by the lattice $S$,
we have the similar to \thetag{5.1} Fourier expansion
related with the group $W=W^{(-2)}(S)$, some its quadratic character
$\epsilon$, a generalized lattice Weyl vector $\rho\in S\otimes \bq$,
and a fundamental polyhedron $\M$ for $W$. For this case,
we can interpret $\br_{++}\M\cap S$ as the set $\text{NEF}(S)$ of
numerically effective (or irreducible with non-negative square)
elements of the Picard lattice
$S$ of $X$ and the set $P(\M)$ as the set of irreducible
non-singular rational curves on $X$. We can interpret
$\alpha >0$, $\alpha \in S$ (for the product part of
\thetag{5.1}), as the
set $\text{EF}(S)$ of effective elements  of the Picard lattice
$S$ for the K3 surfaces $X$. We  get the identity
$$
\split
&\Phi_T (z)=\\
&=\sum_{w \in W^{(-2)}(S)}{\epsilon (w)
\left(\exp{(2\pi i (w(\rho)\cdot z))}-
\sum_{a \in \text{NEF}(S)}{N(a)\exp{(2\pi i (w(\rho+a)\cdot z))}}\right)}\\
&=\exp{\left(2\pi i (\rho\cdot z)\right)}
=\prod_{\alpha \in \text{EF}(S)}{\left(1-
\exp{(2\pi i(\alpha\cdot z))}
\right)^{\mult~\alpha}}
\endsplit
\tag{5.2}
$$
where $z \in \Omega (V^+(S))=S\otimes \br+i V^+(S)$ and $V^+(S)$ is
the light cone of effective elements of $S$.
For the lattices $T=U(4)^2\oplus \langle 2 \rangle$ and
$U(8)^2\oplus \langle 2 \rangle$ these formulae were given in \cite{GN3}.

Considering linear system $|k\rho|$, one gets very canonical and
beautiful projective models of $X$.
The case $T=U(12)^2\oplus \langle 2 \rangle$ corresponding to $\Delta_1$
is especially nice. For this case the
Picard lattice $S=U(12)\oplus \langle -2 \rangle$.
The linear system $|6\rho|$ gives the embedding of $X$ in ${\Bbb P}^4$ as
intersection of quadric and cubic. The surface $X$ has exactly $6$ non-singular
rational curves, all of them have degree $6$.
The set of their classes is  $P_{3,II,\0o}(\M_{3,II,\0o})\subset S$
for the $\rho$ from Sect. 5.1.1.

\smallpagebreak

We recall (see Part I) that a lattice $T$ with
two negative squares is called {\it reflective} if
there exists an automorphic form $\Phi$ on $\Omega (T)$ with respect to
a subgroup of finite index of $O(T)$ such that the zero divisor of $\Phi$
is union of quadratic divisors orthogonal to roots of $T$. The
automorphic form $\Phi$ is called {\it reflective} for $T$. Conjecturally
(see Part I) the set of reflective lattices $T$ is finite up to
multiplication of the form of $T$. It is an
interesting problem to find all reflective lattices $T$ and their
reflective automorphic forms.

In this paper we found many reflective lattices $T$ and reflective
automorphic forms for $\rk T=5$. They are lattices
$$
T=U^2\oplus \langle 2t \rangle,
\qquad t=1,\, 2,\dots,\,9,\,
10,\, 12,\, 16,\, 18,\, 36,
$$
related with the paramodular groups. Using canonical with respect to
$T$ sublattices of $(T\otimes \bq) (k)$
(defined in the proof of Theorem 5.2.1), one
can construct many other reflective lattices.
For example, we get lattices of Theorem 5.2.1.
Without any doubt
these lattices are also very important for the theory of Lorentzian
Kac--Moody algebras, geometry of K3 surfaces and mirror symmetry
(see \cite{GN3} and \cite{GN6}).

\enddocument
\end